\documentclass[11pt,letterpaper]{article}
\usepackage[T1]{fontenc}
\usepackage[margin=1in]{geometry}
\usepackage{amsmath,amsthm,amssymb,mathtools,booktabs,array}
\usepackage{lmodern}
\usepackage{microtype}
\usepackage{enumitem,longtable}
\usepackage{needspace}
\usepackage[authoryear,round]{natbib}
\usepackage[unicode,hypertexnames=false]{hyperref}
\hypersetup{
  colorlinks=true,
  linkcolor=blue,
  citecolor=blue,
  urlcolor=blue,
  filecolor=blue,
  runcolor=blue,
  pdftitle={Necessary and sufficient conditions for the self-normalized central limit theorem for independent triangular arrays},
  pdfauthor={Xiang Li, Song-Hao Liu, Qi-Man Shao, Hao Shi, Jing-Yu Xu},
  pdfkeywords={Self-normalized sums, central limit theorem, independent triangular arrays, maximal-term negligibility}
}
\numberwithin{equation}{section}
\theoremstyle{plain}
\newtheorem{theorem}{Theorem}[section]
\newtheorem{proposition}[theorem]{Proposition}
\newtheorem{lemma}[theorem]{Lemma}
\newtheorem{prooflemma}{Lemma}[section]
\newtheorem{corollary}[theorem]{Corollary}
\theoremstyle{remark}
\newtheorem{remark}{Remark}[section]

\newcommand{\E}{\mathbb E}
\renewcommand{\P}{\mathbb P}
\newcommand{\R}{\mathbb R}
\newcommand{\N}{\mathrm N}
\newcommand{\ind}{\mathbf1}
\newcommand{\dto}{\xrightarrow{\mathrm d}}
\newcommand{\pto}{\xrightarrow{\mathbb P}}

\makeatletter
\newcommand{\proofanchor}{\phantomsection\protected@edef\@currentlabel{\thesection}}
\makeatother
\allowdisplaybreaks[1]
\providecommand{\bmhead}[1]{\section*{#1}}

\title{Necessary and sufficient conditions for the self-normalized central limit theorem for independent triangular arrays}
\author{Xiang Li$^{1}$ \and Song-Hao Liu$^{2}$ \and Qi-Man Shao$^{3}$ \and Hao Shi$^{1}$ \and Jing-Yu Xu$^{1}$}
\date{}

\begin{document}
\maketitle
\vspace{-1.5em}
\begin{center}
\small
\begin{minipage}{0.98\textwidth}
\centering
$^{1}$Department of Statistics and Data Science,\\
Southern University of Science and Technology,\\
Shenzhen, Guangdong, 518055, China\\[0.4em]
$^{2}$School of Mathematical Sciences, Dalian University of Technology,\\
Dalian, Liaoning, China\\[0.4em]
$^{3}$Department of Statistics and Data Science,\\
Shenzhen International Center for Mathematics,\\
Southern University of Science and Technology,\\
Shenzhen, Guangdong, 518055, China\\[0.7em]
\href{mailto:lixiang3@sustech.edu.cn}{\texttt{lixiang3@sustech.edu.cn}};\quad
\href{mailto:liusonghao@dlut.edu.cn}{\texttt{liusonghao@dlut.edu.cn}}\\
\href{mailto:shaoqm@sustech.edu.cn}{\texttt{shaoqm@sustech.edu.cn}};\quad
\href{mailto:Hao.Shi@link.cuhk.edu.hk}{\texttt{Hao.Shi@link.cuhk.edu.hk}}\\
\href{mailto:12131253@mail.sustech.edu.cn}{\texttt{12131253@mail.sustech.edu.cn}}
\end{minipage}
\end{center}

\begin{abstract}
For a rowwise independent triangular array $\{X_{n,i}:1\le i\le n\}$,
put $S_n=\sum_{i=1}^nX_{n,i}$ and $V_n^2=\sum_{i=1}^nX_{n,i}^2$.
We prove that $S_n/V_n\dto\N(0,1)$ implies
$\max_{i\le n}|X_{n,i}|/V_n\pto0$.
The same assumption gives $V_n^2/b_n^2\pto1$ for any choice of medians
$b_n^2$ of $V_n^2$.
Combined with the results of \citet{Shao2018}, this yields necessary and
sufficient conditions for the self-normalized central limit theorem
and resolves his conjecture for independent sequences.
The proof identifies the joint limits of the sum and its sum of squares,
then uses the standard normal ratio to remove every non-Gaussian component.
\end{abstract}

\medskip

\section{Introduction}
\label{main:strategy}
Let $X_1,X_2,\ldots$ be independent real random variables, and put
\[
 S_n=\sum_{i=1}^nX_i,\qquad V_n^2=\sum_{i=1}^nX_i^2.
\]
Ratios with denominator $V_n$ are taken to be zero when $V_n=0$.
We seek necessary and sufficient conditions for $S_n/V_n\dto\N(0,1)$.

For nondegenerate i.i.d. variables, \citet{GGM1997} proved that
$S_n/V_n\dto\N(0,1)$ is equivalent to $\E X_1=0$ and membership
in the domain of attraction of the normal law.
This settled the conjecture of \citet{LMRS1973} and extended the
symmetric result of \citet{GM1991}.
By \citet{OBrien1980}, the domain condition is equivalent to
$\max_{i\le n}|X_i|/V_n\pto0$.
For independent variables symmetric about zero, \citet{Egorov1996}
proved that this maximal-term condition characterizes
$S_n/V_n\dto\N(0,1)$, provided $\P(V_n>0)\to1$.

For general independent sequences, \citet{Shao2018} gave three
sufficient conditions for $S_n/V_n\dto\N(0,1)$:
\begin{align}
 &\frac{\max_{i\le n}|X_i|}{V_n}\pto0,
 \tag{C1}\label{main:shao-C1}\\
 &\sum_{i=1}^n
   \left\{\E\left(\frac{X_i}{V_n}\right)\right\}^{2}
   \longrightarrow0,
 \tag{C2}\label{main:shao-C2}\\
 &\E\left(\frac{S_n}{V_n+\alpha_n}\right)\longrightarrow0.
 \tag{C3}\label{main:shao-C3}
\end{align}
Here, for all sufficiently large $n$, $\alpha_n>0$ is a deterministic
solution of
\begin{equation}\label{main:shao-scale}
 \sum_{i=1}^n\E\frac{X_i^2}{\alpha_n^2+X_i^2}=1.
\end{equation}
Under \eqref{main:shao-C1}, Shao also proved necessity of
\eqref{main:shao-C2} and~\eqref{main:shao-C3}, and conjectured that
all three conditions are necessary.
The remaining implication is from the normal limit to \eqref{main:shao-C1}.

Condition~\eqref{main:shao-C1} compares the observed maximum with the
random denominator $V_n$.
For centered, rowwise independent variables $\xi_{n,i}$ with finite
variances, the deterministic variance scale is
$s_n^2=\sum_{i=1}^n\E\xi_{n,i}^2>0$.
The classical Lindeberg--Feller theorem states that, under Feller's
condition $\max_{i\le n}\E\xi_{n,i}^2/s_n^2\to0$,
the convergence $\sum_{i=1}^n\xi_{n,i}/s_n\dto\N(0,1)$ holds if and only if
\[
 \frac1{s_n^2}\sum_{i=1}^n
 \E\bigl[\xi_{n,i}^2\ind_{\{|\xi_{n,i}|>\varepsilon s_n\}}\bigr]
 \longrightarrow0\qquad(\varepsilon>0);
\]
see \citet{Lindeberg1922,Feller1936}.
The classical converse assumes negligibility of individual variances.
For the random scale $V_n$, we prove that the normal ratio itself
implies \eqref{main:shao-C1} and concentration under median normalization.

This proves Shao's conjecture: \eqref{main:shao-C1}--\eqref{main:shao-C3}
are necessary and sufficient.
More precisely, let $b_n^2$ be any median of $V_n^2$, with $b_n\ge0$.
If $S_n/V_n\dto\N(0,1)$, then $b_n>0$ eventually and
\begin{subequations}\label{intro:median-normalization}
\begin{gather}
 V_n^2/b_n^2\pto1,\qquad S_n/b_n\dto\N(0,1),
 \label{intro:normalization}\\
 \sum_{i=1}^n\P(|X_i|>\varepsilon b_n)\longrightarrow0
 \qquad\text{for every }\varepsilon>0.
 \label{intro:tail}
\end{gather}
\end{subequations}
These conclusions extend to rowwise independent triangular arrays.

Under the normal limit, median normalization makes
$(S_n/b_n,V_n^2/b_n^2)$ tight.
We retain coordinate limits in the joint-limit approach studied for
infinitesimal arrays by \citet{Mason2005,Mason2013}.
Exponential integrability and entropy comparisons under quadratic tilts
show that the Gaussian component is nondegenerate.
Integrating over this component excludes an extra deterministic term in the square sum.
The normal ratio gives lower bounds tending to one for two signed
exponential transforms along the same sequence of parameters tending to infinity.
A nonzero coordinate or Poisson component gives an asymptotic upper bound
strictly below one for their product, which yields a contradiction.
The zero mean of the normal ratio then eliminates the remaining drift,
and the median property fixes the variance at one.
A further subsequence argument and Laplace factorization then give
\eqref{intro:tail}.

\textbf{Notation} The symbols $\dto$ and $\pto$ denote convergence in
distribution and in probability.
Finite positive constants $C,c$ may change between occurrences;
dependence on fixed quantities is stated where needed.

\section{Main results}
\label{main:results}
For each $n\ge1$, let $X_{n,1},\ldots,X_{n,n}$ be independent,
almost surely finite real random variables. Put
\[
 S_n=\sum_{i=1}^{n}X_{n,i},\qquad
 V_n=\left(\sum_{i=1}^{n}X_{n,i}^2\right)^{1/2},
\]
and define
\[
 T_n=\begin{cases}S_n/V_n,&V_n>0,\\0,&V_n=0.\end{cases}
\]
Every ratio with denominator $V_n$ is taken to be zero when $V_n=0$.

Our main theorem shows that a standard normal limit for $T_n$ implies
maximal-term negligibility and yields a deterministic normalization.

\Needspace{15\baselineskip}
\begin{theorem}
\label{main:array}
Under the assumptions above, suppose that
\begin{equation}\label{hist:cond:N}
 T_n\dto\N(0,1).
\end{equation}
Then
\begin{equation}\label{main:maximal}
 \frac{\max_{i\le n}|X_{n,i}|}{V_n}\pto0.
\end{equation}
More precisely, let $b_n^2$ be any median of $V_n^2$, and let $b_n$ denote its nonnegative square root. Thus
\[
 \P(V_n^2\le b_n^2)\ge\tfrac12,
 \qquad \P(V_n^2\ge b_n^2)\ge\tfrac12.
\]
Then $b_n>0$ for all sufficiently large $n$, and
\begin{subequations}\label{main:strong}
\begin{gather}
 V_n^2/b_n^2\pto1,\qquad S_n/b_n\dto\N(0,1),
 \label{main:normalization}\\
 \sum_{i=1}^{n}\P(|X_{n,i}|>\varepsilon b_n)\longrightarrow0
 \qquad\text{for every }\varepsilon>0.
 \label{main:tail}
\end{gather}
\end{subequations}
\end{theorem}

The theorem holds for every deterministic choice of medians.
In particular, $V_n/b_n\pto1$, so the same normal limit holds with
deterministic normalization. Reindexing, finite row repetition, and zero
padding extend the results to arrays with arbitrary finite row lengths.

Combining Theorem~\ref{main:array} with the results of
\citet[Theorems~1.1 and~1.2, p.~1742]{Shao2018}
gives the following necessary and sufficient conditions for the
self-normalized central limit theorem.

\begin{theorem}
\label{main:shao-characterization}
Under the rowwise independence assumptions above,
\eqref{hist:cond:N} holds if and only if
\eqref{main:shao-C1}--\eqref{main:shao-C3} hold with $X_i$ replaced
by $X_{n,i}$ in row $n$. The same replacement is made in
\eqref{main:shao-scale}.
\end{theorem}

For the array $X_{n,i}=X_i$, the normal limit \eqref{hist:cond:N}
and maximal-term negligibility \eqref{main:maximal} also force
$V_n\uparrow\infty$ almost surely.

\begin{corollary}
\label{cor:fixed-sequence}
Let $(X_i)_{i\ge1}$ be a fixed sequence of independent, almost surely
finite real random variables. Define $S_n,V_n,T_n$ from its first
$n$ coordinates, with the same zero-denominator convention. If
$T_n\dto\N(0,1)$ along the full sequence, then
\[
 \frac{\max_{i\le n}|X_i|}{V_n}\pto0,
 \qquad \sum_{i=1}^\infty X_i^2=\infty\quad\text{almost surely}.
\]
In particular, if $\sum_{i=1}^\infty X_i$ converges almost surely
and $0<\sum_{i=1}^\infty X_i^2<\infty$ almost surely, then
\begin{equation}\label{cor:finite-series-ratio}
 \frac{\sum_{i=1}^\infty X_i}
      {\left(\sum_{i=1}^\infty X_i^2\right)^{1/2}}
 \not\sim\N(0,1).
\end{equation}
\end{corollary}

\Needspace{16\baselineskip}
\section{Proofs of the main results}
\label{sec:return-array}

\subsection{Proof of Theorem~\ref{main:array}}
\label{sec:proof-main-array}
The first step is to control the square sum at a deterministic scale.
The lemma below shows that a uniform positive lower bound for
$\P(V_n^2\le d_n^2)$ gives tightness of $V_n^2/d_n^2$.

\proofanchor
\label{hist:sec:compact}

\begin{prooflemma}
\label{hist:lem:small-ball}
Suppose that $T_n\dto T$ and $\P(T=1)=\P(T=-1)=0$.
Let $d_n>0$ be deterministic numbers such that, for some fixed $c>0$,
\begin{equation}\label{hist:eq:small-ball}
 \P(V_n^2\le d_n^2)\ge c
 \qquad\text{for every }n.
\end{equation}
Then the family $V_n^2/d_n^2$ is tight, and
\begin{equation}\label{hist:eq:aggregate-large}
 \lim_{R\to\infty}\limsup_{n\to\infty}
       \sum_{i\le n}\P(|X_{n,i}|>Rd_n)=0.
\end{equation}
\end{prooflemma}

\begin{proof}
Write
\[
 S_{n,-i}=\sum_{\substack{j\le n\\j\ne i}}X_{n,j},\qquad
 V_{n,-i}^2=\sum_{\substack{j\le n\\j\ne i}}X_{n,j}^2.
\]
By \eqref{hist:eq:small-ball},
\[
 \P(V_{n,-i}^2\le d_n^2)\ge c,\qquad
 \P(|X_{n,i}|\le d_n)\ge c.
\]
For $K>1$, the events
$\{V_{n,-i}^2\le d_n^2,\ |S_{n,-i}|>Kd_n\}$ and
$\{|X_{n,i}|\le d_n\}$ are independent.
On their intersection,
$|S_n|>(K-1)d_n$ and $V_n^2\le2d_n^2$, so
$|T_n|>(K-1)/\sqrt2$. Hence
\begin{equation}\label{hist:eq:remove-coordinate}
 \P(V_{n,-i}^2\le d_n^2,\ |S_{n,-i}|>Kd_n)
 \le c^{-1}\P\!\left(|T_n|>\frac{K-1}{\sqrt2}\right).
\end{equation}
By tightness of $(T_n)$, choose a fixed $K>1$ so that the right-hand
side is at most $c/2$ for every $n$. Thus
\[
 G_{n,i}=\{V_{n,-i}^2\le d_n^2,\ |S_{n,-i}|\le Kd_n\},
 \qquad \P(G_{n,i})\ge c/2,
\]
for every $n$ and $1\le i\le n$.

For $R>1$, set $E_{n,i}=\{|X_{n,i}|>Rd_n\}\cap G_{n,i}$.
If $E_{n,i}$ occurs, then $|X_{n,j}|\le d_n$ for $j\ne i$;
hence the events $E_{n,i}$ are disjoint. On $E_{n,i}$,
\begin{align*}
 \left|T_n-\operatorname{sgn}X_{n,i}\right|
 &\le\frac{|S_{n,-i}|}{\sqrt{X_{n,i}^2+V_{n,-i}^2}}
       +1-\frac{|X_{n,i}|}{\sqrt{X_{n,i}^2+V_{n,-i}^2}}
 \le\frac KR+\frac1{2R^2}=: \eta_R,
\end{align*}
where the second inequality uses $1-(1+u)^{-1/2}\le u/2$ for $u\ge0$.
Independence of $X_{n,i}$ and $G_{n,i}$ gives
\begin{equation}\label{hist:eq:one-jump-bound}
 \frac c2\sum_i\P(|X_{n,i}|>Rd_n)
 \le\sum_i\P(E_{n,i})
 \le\P\{\operatorname{dist}(T_n,\{-1,1\})\le\eta_R\}.
\end{equation}
For fixed $R$, the set $\{t:\operatorname{dist}(t,\{-1,1\})\le\eta_R\}$
is closed. The Portmanteau theorem gives
\[
 \limsup_n\sum_i\P(|X_{n,i}|>Rd_n)
 \le\frac2c\P\{\operatorname{dist}(T,\{-1,1\})\le\eta_R\}.
\]
Since $\P(T=1)=\P(T=-1)=0$, letting $R\to\infty$ with $K$ fixed
proves \eqref{hist:eq:aggregate-large}.

Put $f_{n,i}=\E e^{-X_{n,i}^2/d_n^2}\in(0,1]$.
By \eqref{hist:eq:small-ball} and independence,
\[
 \prod_i f_{n,i}=\E e^{-V_n^2/d_n^2}\ge c/e,
 \qquad
 \sum_i(1-f_{n,i})\le-\log\prod_i f_{n,i}\le1-\log c.
\]
For $u\ge0$, the inequality
\begin{equation}\label{hist:eq:laplace-truncation}
 1-e^{-u}\ge(1-e^{-1})\min(u,1)
\end{equation}
therefore yields
\begin{equation}\label{hist:eq:truncated-energy}
 \sum_i\E\min(X_{n,i}^2/d_n^2,1)
 \le C_c:=\frac{1-\log c}{1-e^{-1}}.
\end{equation}
For $R>1$ and $L>0$, on $\{\max_i|X_{n,i}|\le Rd_n\}$,
$V_n^2=\sum_iX_{n,i}^2\ind_{\{|X_{n,i}|\le Rd_n\}}$.
The union bound, Markov's inequality and
$y^2\ind_{\{|y|\le R\}}\le R^2\min(y^2,1)$ give
\begin{equation}\label{hist:eq:q-tight-bound}
 \P(V_n^2>L^2d_n^2)
 \le\sum_i\P(|X_{n,i}|>Rd_n)+\frac{R^2C_c}{L^2}.
\end{equation}
Fix $\varepsilon>0$. By \eqref{hist:eq:aggregate-large}, we can choose $R>1$
such that
\[
 \limsup_{n\to\infty}\sum_{i\le n}\P(|X_{n,i}|>Rd_n)
 <\varepsilon/4.
\]
With this $R$ fixed, choose $L$ large enough that
$R^2C_c/L^2<\varepsilon/4$.
By \eqref{hist:eq:q-tight-bound},
\[
 \limsup_{n\to\infty}\P(V_n^2>L^2d_n^2)<\varepsilon/2.
\]
Hence there is an $N$ such that
$\P(V_n^2>L^2d_n^2)<\varepsilon/2$ for every $n\ge N$.
For each of the finitely many indices $n<N$, the ratio $V_n^2/d_n^2$
is almost surely finite, so its tail probability tends to zero as
$L\to\infty$. Increasing $L$ if necessary makes the same bound hold
for all the first $N$ indices. Thus
\begin{equation}\label{hist:eq:scale-tightness}
 \sup_{n\ge1}\P(V_n^2>L^2d_n^2)\le\varepsilon/2<\varepsilon,
\end{equation}
which proves tightness of the whole family.
\end{proof}

The next proposition applies Lemma~\ref{hist:lem:small-ball} at the
median scale and at smaller scales. These two applications give joint
tightness and keep the normalized square sum away from zero.

\begin{proposition}
\label{hist:prop:median}
Assume \eqref{hist:cond:N}.  For each $n$, let $b_n^2$ be any median
of $V_n^2$.  Then $b_n>0$ eventually, and
\begin{align}
 \lim_{M\to\infty}\limsup_{n\to\infty}
       \P(V_n^2>Mb_n^2)&=0,
       \label{hist:eq:median-upper}\\
 \lim_{\delta\downarrow0}\limsup_{n\to\infty}
       \P(V_n^2<\delta b_n^2)&=0.
       \label{hist:eq:median-lower}
\end{align}
After discarding the finite initial segment on which a median might
be zero, the pairs $(S_n/b_n,V_n^2/b_n^2)$ are tight.  Every
subsequential joint limit $(A,Q)$ satisfies
\begin{equation}\label{hist:eq:joint-positive}
 0<Q<\infty\quad\text{almost surely},\qquad
 A/\sqrt Q\sim\N(0,1).
\end{equation}
\end{proposition}

\begin{proof}[Proof of Proposition~\ref{hist:prop:median}]
By the convention defining $T_n$, the event $\{V_n^2=0\}$ is contained
in $\{T_n=0\}$. The Portmanteau theorem gives
\[
 \limsup_{n\to\infty}\P(V_n^2=0)
 \le\limsup_{n\to\infty}\P(T_n=0)
 =0,
\]
because the standard normal law has no atom at zero.
A median $b_n^2$ satisfies
\[
 \P(V_n^2\le b_n^2)\ge\tfrac12,\qquad
 \P(V_n^2\ge b_n^2)\ge\tfrac12.
\]
Since $\P(V_n^2=0)<1/2$ eventually, $b_n^2>0$ for all sufficiently
large $n$. Discard the finite initial segment and take $b_n>0$.
Then \eqref{hist:eq:small-ball} holds with $d_n=b_n$ and $c=1/2$ and
Lemma~\ref{hist:lem:small-ball} gives \eqref{hist:eq:median-upper}.

To prove \eqref{hist:eq:median-lower}, suppose that it fails.
There are a subsequence $n_j$, numbers $\delta_j\downarrow0$,
and $\rho>0$ such that
\[
 \P(V_{n_j}^2\le\delta_jb_{n_j}^2)\ge\rho
 \qquad\text{for every }j.
\]
Apply Lemma~\ref{hist:lem:small-ball} to the rows $n_j$ with
$d_{n_j}=\sqrt{\delta_j}b_{n_j}$ and $c=\rho$.
It gives tightness of $V_{n_j}^2/(\delta_jb_{n_j}^2)$.
The median property gives
\[
 \P\left(\frac{V_{n_j}^2}{\delta_jb_{n_j}^2}\ge\delta_j^{-1}\right)
 =\P(V_{n_j}^2\ge b_{n_j}^2)\ge\tfrac12.
\]
Since $\delta_j^{-1}\to\infty$, this contradicts tightness and proves
\eqref{hist:eq:median-lower}.

The identity
\[
 \frac{S_n}{b_n}=T_n\frac{V_n}{b_n}
\]
holds also when $V_n^2=0$, since both sides then vanish.  Tightness of
the two factors implies tightness of their product. Together with
\eqref{hist:eq:median-upper}, this proves joint tightness of
$(S_n/b_n,V_n^2/b_n^2)$.

Consider any convergent subsequence with limit $(A,Q)$.  Its second
coordinate is nonnegative and finite almost surely.  For every
$\delta>0$, the open-set part of the Portmanteau theorem gives
\[
 \P(Q<\delta)
 \le\liminf_j\P(V_{n_j}^2/b_{n_j}^2<\delta)
 \le\limsup_n\P(V_n^2/b_n^2<\delta).
\]
Letting $\delta\downarrow0$ and using
\eqref{hist:eq:median-lower} yields $\P(Q=0)=0$.
The map $(a,q)\mapsto a/\sqrt q$ is continuous whenever $q>0$;
define it to be zero at $q=0$ for the prelimit pairs.  The continuous
mapping theorem therefore identifies its limiting law with that of
$A/\sqrt Q$.  Since the same prelimit ratios are $T_{n_j}$ and
converge to $\N(0,1)$, \eqref{hist:eq:joint-positive} follows.

\end{proof}

With tightness established, the next step is to identify the possible
joint limits.
Rowwise independence yields a representation separating the Gaussian,
individual-coordinate, and Poisson contributions, together with
deterministic terms.

\begin{proposition}
\label{array:joint-representation}
Let $\{Y_{n,i}:1\le i\le n\}$ be a triangular array of real random
variables, independent within each row. Set
\[
 A_n=\sum_{i=1}^nY_{n,i},\qquad
 \widehat V_n=\left(\sum_{i=1}^nY_{n,i}^2\right)^{1/2}.
\]
Suppose that
\[
 (A_n,\widehat V_n^2)\dto(A,Q),\qquad
 Q<\infty\quad\text{almost surely}.
\]
Then the joint limit admits the representation
\begin{equation}\label{array:joint-canonical-form}
 (A,Q)\overset d=
 \left(b+\sigma Z+U_{\mathrm{ord}}+U_{\mathrm{poi}},
       \;\sigma^2+\kappa+R_{\mathrm{ord}}+R_{\mathrm{poi}}\right),
\end{equation}
where $b\in\mathbb R$ and $\sigma,\kappa\ge0$ are constants,
$Z\sim\N(0,1)$, and
\begin{equation}\label{def:components}
 \begin{split}
 U_{\mathrm{ord}}
 &=\sum_j\left\{
      Y_j\ind_{\{|Y_j|\le1\}}
      -\E[Y_j\ind_{\{|Y_j|\le1\}}]
    \right\}
   +\sum_jY_j\ind_{\{|Y_j|>1\}},\\
 R_{\mathrm{ord}}&=\sum_jY_j^2,\\
 U_{\mathrm{poi}}
 &=\int_{|x|\le1}x(\mathcal N-\nu)(dx)
   +\int_{|x|>1}x\mathcal N(dx),\\
 R_{\mathrm{poi}}&=\int x^2\mathcal N(dx).
 \end{split}
\end{equation}
Here $Y_1,Y_2,\ldots$ are real-valued random variables and
$\mathcal N$ is a Poisson random measure with intensity $\nu$.
The variables $Y_1,Y_2,\ldots,Z$ and the random measure $\mathcal N$
are mutually independent, and
\begin{equation}\label{array:component-summability}
 \sum_j\E(1\wedge Y_j^2)<\infty,\qquad
 \int(1\wedge x^2)\nu(dx)<\infty,\qquad
 \nu(\{0\})=0.
\end{equation}
\end{proposition}

\begin{remark}
\label{rem:representation-convergence}
By \eqref{array:component-summability}, the first series in
$U_{\mathrm{ord}}$ converges almost surely and in $L^2$.
The second sum has only finitely many nonzero terms almost surely.
The compensated integral in $U_{\mathrm{poi}}$ is an $L^2$ limit.
It also converges almost surely along nested truncations away from zero.
Both $R_{\mathrm{ord}}$ and $R_{\mathrm{poi}}$ are finite almost surely.
\end{remark}

\proofanchor
\label{array:joint-representation-section}

For the infinitesimal part, we use the following joint limit
representation of the sum and square sum. It is based on the classical
array theorem of \citet[Chapter~2, Theorem~4.7, p.~61]{AG1980}, as used by
\citet[pp.~854--855]{Mason2005}, read with the correction \citep{Mason2013}.

Put
\begin{equation}
\label{def:truncations}
 w(x)=\frac{x^2}{1+x^2},\qquad
 \tau(x)=\frac{x}{1+x^2},\qquad
 \tau_0(x)=x\ind_{\{|x|\le1\}}.
\end{equation}
\begin{prooflemma}
\label{app:infinitesimal-remainder}
Let $\{Z_{n,i}:1\le i\le n\}$ have independent, almost surely finite
real entries within each finite row, and assume
\[
 \max_{i\le n}\P(|Z_{n,i}|>\epsilon)\longrightarrow0
 \qquad\text{for every }\epsilon>0.
\]
With $w,\tau,\tau_0$ as in \eqref{def:truncations}, put
\[
 \mu_n(dx)=w(x)\sum_i\mathcal L(Z_{n,i})(dx),\qquad
 \kappa_n=\sum_i[\E\tau(Z_{n,i})]^2.
\]
Here $\mathcal L(X)$ denotes the law of $X$.
Suppose that $\mu_n\Rightarrow\mu$ weakly as finite measures on $\R$
and that $\kappa_n\to\kappa<\infty$. Define
\[
 \nu(dx)=\ind_{\{x\ne0\}}w(x)^{-1}\mu(dx).
\]
Then $0\le\kappa\le \mu(\{0\})$, $\nu$ is a L\'evy measure, and, with
$\sigma^2=\mu(\{0\})-\kappa$ and $\gamma_\nu=\int(\tau_0-\tau)\,d\nu$,
\begin{equation}\label{app:remainder-limit}
 \left(\sum_i\{Z_{n,i}-\E\tau(Z_{n,i})\},\ \sum_iZ_{n,i}^2\right)
 \dto(\sigma Z+U_{\mathrm{poi}}+\gamma_\nu,\ \sigma^2+\kappa+R_{\mathrm{poi}}).
\end{equation}
Here $Z\sim\N(0,1)$ is independent of a Poisson random measure
$\mathcal N$ with intensity $\nu$; $U_{\mathrm{poi}}$ and
$R_{\mathrm{poi}}$ are defined in \eqref{def:components}.
The integral defining $\gamma_\nu$ is absolutely convergent.
\end{prooflemma}

The proof is given in Supplementary Section~\ref{app:infinitesimal}.

\begin{proof}[Proof of Proposition~\ref{array:joint-representation}]
We pass to further subsequences below and retain the index $n$
throughout the proof.

\medskip
\noindent\textbf{Step 1.}Let us first show that centering preserves
tightness of the row sums. Tightness of $\widehat V_n^2$ allows us, after discarding finitely many rows, to
choose $M_0<\infty$ such that $\inf_n\P(\widehat V_n^2\le M_0)\ge1/2$.
Write $p_{n,i}=\E(1-e^{-Y_{n,i}^2})$. Independence gives
\[
 \tfrac12e^{-M_0}\le\E e^{-\widehat V_n^2}
 =\prod_i(1-p_{n,i})
 \le\exp\left(-\sum_i p_{n,i}\right).
\]
The inequality $e^{x^2}\ge1+x^2$ gives
$w(x)\le1-e^{-x^2}$. Hence
\begin{equation}\label{array:G1}
 C:=\sup_n\sum_i\E w(Y_{n,i})\le  \sup_n \sum_{i}p_{n,i}\le M_0+\log 2<\infty.
\end{equation}
Also, whenever $\P(\widehat V_n^2\le L^2)>0$, independence and
$p\le-\log(1-p)$ yield
\begin{equation}\label{array:G2}
 \sum_i\P(|Y_{n,i}|>L)
 \le -\sum_{i=1}^{n}\log\P\left(|Y_{n,i}|\le L\right)
 =-\log\P(\max_i|Y_{n,i}|\le L)
 \le-\log\P(\widehat V_n^2\le L^2).
\end{equation}
Tightness of $\widehat V_n^2$ gives
$\inf_n\P(\widehat V_n^2\le L^2)\to1$ as $L\to\infty$.
Thus $\sup_n\sum_i\P(|Y_{n,i}|>L)\to0$ as $L\to\infty$.

{
Put $a_{n,i}=\E\tau(Y_{n,i})$. For every deterministic
subset $I_n\subseteq\{1,\ldots,n\}$,
}{\begin{equation}\label{array:centered-sum-decomposition}
\begin{aligned}
 \sum_{i\in I_n}(Y_{n,i}-a_{n,i})
 &=\sum_{i\in I_n}\{\tau(Y_{n,i})-a_{n,i}\}
 +\sum_{i\in I_n}
 (Y_{n,i}-\tau(Y_{n,i}))\ind_{\{|Y_{n,i}|\le1\}}\\
 &\quad+\sum_{i\in I_n}
 (Y_{n,i}-\tau(Y_{n,i}))\ind_{\{|Y_{n,i}|>1\}}.
\end{aligned}
\end{equation}
}
Rowwise independence, $\tau^2\le w$, and
\eqref{array:G1} give
{\begin{equation}\label{array:truncated-sum-second-moment}
\begin{aligned}
 \E\left|\sum_{i\in I_n}\{\tau(Y_{n,i})-a_{n,i}\}\right|^2
 =\sum_{i\in I_n}\operatorname{Var}(\tau(Y_{n,i}))
 \le\sum_{i\in I_n}\E w(Y_{n,i})
 \le C.
\end{aligned}
\end{equation}
}For $|x|\le1$, $|x-\tau(x)|\le w(x)$, so
{\begin{equation}\label{array:small-remainder-first-moment}
 \E\left[
   \sum_{i\in I_n}|Y_{n,i}-\tau(Y_{n,i})|
   \ind_{\{|Y_{n,i}|\le1\}}
 \right]
 \le\sum_{i\in I_n}\E w(Y_{n,i})
 \le C.
\end{equation}
}For $|x|>1$, $|x-\tau(x)|\le|x|$ and the
Cauchy--Schwarz inequality yield
{\begin{equation}\label{array:large-remainder-bound}
 \sum_{i\in I_n}|Y_{n,i}-\tau(Y_{n,i})|
       \ind_{\{|Y_{n,i}|>1\}}
 \le\sqrt{N_n\widehat V_n^2},
 \qquad
 N_n:=\sum_{i=1}^n\ind_{\{|Y_{n,i}|>1\}}.
\end{equation}
}Since $\ind_{\{|x|>1\}}\le2w(x)$, \eqref{array:G1} implies
\[
 \sup_n\E N_n
 \le2\sup_n\sum_{i=1}^n\E w(Y_{n,i})
 \le2C.
\]
Thus $N_n$ is tight by Markov's inequality. Together with tightness
of $\widehat V_n^2$, this implies tightness of $N_n\widehat V_n^2$.

{Apply Chebyshev's inequality to \eqref{array:truncated-sum-second-moment}
and Markov's inequality to \eqref{array:small-remainder-first-moment}.
Bound the large terms using \eqref{array:large-remainder-bound}.
The decomposition \eqref{array:centered-sum-decomposition} and
the union bound then give, for every $M>0$,
}\[
 \P\!\left(
   \left|\sum_{i\in I_n}(Y_{n,i}-a_{n,i})\right|>3M
 \right)
 \le \frac{C}{M^2}+\frac{C}{M}
      +\P\!\left(\sqrt{N_n\widehat V_n^2}>M\right).
\]
This upper bound does not depend on $I_n$ and tends to zero uniformly
in $n$ as $M\to\infty$. Hence
\begin{equation}\label{array:G3}
 \left\{\sum_{i\in I_n}(Y_{n,i}-a_{n,i}):n\ge1\right\}
 \quad\text{is tight uniformly over all deterministic choices of }I_n.
\end{equation}

Taking $I_n=\{1,\ldots,n\}$ and writing
$a_n=\sum_{i=1}^n a_{n,i}$ gives tightness of $A_n-a_n$.
Since $A_n$ is tight by assumption, their difference
$a_n=A_n-(A_n-a_n)$ is also tight. Being deterministic, the sequence
$(a_n)$ is therefore bounded.

\medskip
\noindent\textbf{Step 2.} We will next split each row sum into a sum of
larger individual terms and an infinitesimal remainder, then study
their limiting distributions separately. If $Y_{n,i_n}\dto Y$ is prescribed, put $Y_{n,i_n}$ first in each
row. Otherwise, choose any coordinate first. These permutations leave
$A_n$ and $\widehat V_n^2$ unchanged. Order the remaining coordinates
by decreasing $\E w(Y_{n,j})$ for $j\ge2$. Append zeros beyond
the end of each row. For each fixed $j$, the family $(Y_{n,j})_n$
is tight because
\[
 \P(|Y_{n,j}|>L)\le\P(\widehat V_n^2>L^2).
\]
Choose nested infinite index sets
$\mathcal I_1\supseteq\mathcal I_2\supseteq\cdots$ so that each
of the first $j$ coordinates has a distributional limit along
$\mathcal I_j$. Let $n_k$ be the $k$th index in $\mathcal I_k$.
Then $Y_{n_k,j}$ converges in distribution for every fixed $j$.
Let $Y_1,Y_2,\ldots$ be independent random variables with these
limiting laws. If a limit $Y$ was prescribed, the first coordinate
still converges to $Y$, so $Y_1\overset d=Y$. Independence within
each row gives
\[
 (Y_{n_k,1},\ldots,Y_{n_k,m})
 \dto(Y_1,\ldots,Y_m)
 \qquad\text{for every fixed }m.
\]
 Since $w$ is bounded and continuous, Fatou's lemma for the
sum over $j$ gives
\begin{equation}\label{array:G4}
 \sum_j\E w(Y_j)\le C,
 \qquad \sum_j\E(1\wedge Y_j^2)\le2C.
\end{equation}
To prove convergence of the ordinary-coordinate contribution, put
$a_j=\E\tau(Y_j)$ and split its partial sums into a martingale
and a correction:
\[
 \begin{aligned}
 \sum_{j=1}^m(Y_j-a_j)&=H_m+\sum_{j=1}^m d_j,\quad 
 H_m:=\sum_{j=1}^m\{\tau(Y_j)-a_j\},
 \quad d_j:=Y_j-\tau(Y_j).
 \end{aligned}
\]
By independence and centering, $(H_m)$ is a martingale for the filtration
$\{\sigma(Y_1,\ldots,Y_m)\}_{m\ge1}$. Since $\tau^2\le w$,
\[
 \E H_m^2
 =\sum_{j=1}^m\operatorname{Var}(\tau(Y_j))
 \le\sum_{j=1}^m\E w(Y_j)\le C.
\]
The $L^2$-bounded martingale convergence theorem therefore gives
$H_m\to H_\infty$ almost surely and in $L^2$.

{For $d_j=Y_j-\tau(Y_j)$, $|x-\tau(x)|\le w(x)$ when $|x|\le1$, so Tonelli's
theorem and \eqref{array:G4} give
}\[
 \E\!\left[\sum_{j=1}^\infty |d_j|\ind_{\{|Y_j|\le1\}}\right]
 =\sum_{j=1}^\infty\E\bigl[|d_j|\ind_{\{|Y_j|\le1\}}\bigr]
 \le C.
\]
Thus this nonnegative sum is finite almost surely. Moreover,
{\begin{equation}\label{array:ordinary-large-value-summability}
 \sum_{j=1}^\infty\P(|Y_j|>1)
 \le2\sum_{j=1}^\infty\E w(Y_j)\le2C.
\end{equation}
}{By the first Borel--Cantelli lemma, only finitely many terms
$d_j\ind_{\{|Y_j|>1\}}$ are nonzero almost surely. Each term is finite, so
$\sum_j|d_j|<\infty$ almost surely. Therefore
}\begin{equation}\label{array:G5}
 F:=\sum_{j=1}^\infty(Y_j-a_j)
   =H_\infty+\sum_{j=1}^\infty d_j\in\mathbb R
 \quad\text{almost surely}.
\end{equation}
Changing the centering from $\tau$ to $\tau_0$ expresses $F$ as
$U_{\mathrm{ord}}$ plus a deterministic correction.
The bound
\[
 \sum_j\operatorname{Var}(\tau_0(Y_j))
 \le\sum_j\E[Y_j^2\ind_{\{|Y_j|\le1\}}]\le2C
\]
and the $L^2$-bounded martingale convergence theorem show that
$\sum_j\{\tau_0(Y_j)-\E\tau_0(Y_j)\}$ converges almost surely and in
$L^2$. Together with the finitely many large values, this defines
$U_{\mathrm{ord}}$ in \eqref{def:components}. Moreover, $|\tau_0-\tau|\le w$
implies
\[
 \sum_j\left|\E[(\tau_0-\tau)(Y_j)]\right|
 \le\sum_j\E w(Y_j)\le C.
\]
For each $m$, the identity
\[
 \begin{aligned}
 \sum_{j=1}^m(Y_j-a_j)
 &=\sum_{j=1}^m\{\tau_0(Y_j)-\E\tau_0(Y_j)\}
   +\sum_{j=1}^mY_j\ind_{\{|Y_j|>1\}}
+\sum_{j=1}^m\E[(\tau_0-\tau)(Y_j)]
 \end{aligned}
\]
therefore gives, on letting $m\to\infty$,
{\begin{equation}\label{array:ordinary-centering-identity}
 F=U_{\mathrm{ord}}+\sum_j\E[(\tau_0-\tau)(Y_j)].
\end{equation}
}

The bound $\sum_j\E[Y_j^2\ind_{\{|Y_j|\le1\}}]\le2C$ also implies
that the small-value square sum is finite almost surely.
By \eqref{array:ordinary-large-value-summability}, only finitely many
$Y_j^2\ind_{\{|Y_j|>1\}}$ are nonzero almost surely.
Hence $R_{\mathrm{ord}}=\sum_jY_j^2<\infty$ almost surely.
For every fixed positive integer $m$, bounded continuity of $\tau$ and
finite-dimensional convergence give
\[
 P_{n,m}:=\left(\sum_{j\le m}(Y_{n,j}-a_{n,j}),
                   \sum_{j\le m}Y_{n,j}^2\right)
 \dto P_m:=\left(\sum_{j\le m}(Y_j-a_j),
                   \sum_{j\le m}Y_j^2\right).
\]
Moreover, $P_m\to(F,R_{\mathrm{ord}})$ almost surely as $m\to\infty$.
Let $d$ be the bounded Lipschitz metric on probability laws on $\mathbb R^2$.
To approximate $(F,R_{\mathrm{ord}})$ by a growing initial block,
choose strictly increasing integers $n_m$ such that, for every $n\ge n_m$,
$d(\mathcal L(P_{n,m}),\mathcal L(P_m))\le1/m$.
For all sufficiently large $n$, set
\[
 m_n=\max\{m\le n:n_m\le n\}.
\]
The sequence $m_n$ is nondecreasing, and for every fixed $M$,
we have $m_n\ge M$ whenever $n\ge\max\{M,n_M\}$.
Hence $m_n\uparrow\infty$, and
\[
 \begin{aligned}
 d\bigl(\mathcal L(P_{n,m_n}),\mathcal L(F,R_{\mathrm{ord}})\bigr)
 &\le m_n^{-1}
   +d\bigl(\mathcal L(P_{m_n}),\mathcal L(F,R_{\mathrm{ord}})\bigr)
 \longrightarrow0.
 \end{aligned}
\]
Therefore
\begin{equation}\label{array:G6}
 P_{n,m_n}\dto(F,R_{\mathrm{ord}}).
\end{equation}

{
For $j\ge2$, the ordering of coordinates $2,3,\ldots$ and
\eqref{array:G1} give
\[
 (j-1)\E w(Y_{n,j})
 \le\sum_{k=2}^{j}\E w(Y_{n,k})\le C.
\]
Thus $\E w(Y_{n,j})\le C/m_n$ for $j>m_n$.
For every $\epsilon>0$,
}\begin{equation}\label{array:G7}
 \P(|Y_{n,j}|>\epsilon)
 \le\frac{1+\epsilon^2}{\epsilon^2}\,\E w(Y_{n,j})
 \le\frac{1+\epsilon^2}{\epsilon^2}\frac{C}{m_n}
 \longrightarrow0.
\end{equation}
Boundedness of $(a_n)$, \eqref{array:G3}, and tightness of
$\widehat V_n^2$ give a further subsequence on which
$a_n\to a\in\mathbb R$ and
\begin{equation}\label{array:G8}
 (D_n,E_n):=
 \left(\sum_{j>m_n}(Y_{n,j}-a_{n,j}),
             \sum_{j>m_n}Y_{n,j}^2\right)
 \dto(D,E).
\end{equation}
For each $n$, $(D_n,E_n)$ and $P_{n,m_n}$ depend on disjoint sets
of independent coordinates. Their limits $(D,E)$ and
$(F,R_{\mathrm{ord}})$ are therefore independent. Using
\eqref{array:G6}, \eqref{array:G8}, and $a_n\to a$, we pass
to the limit in
\[
 (A_n,\widehat V_n^2)
 =(a_n,0)+P_{n,m_n}+(D_n,E_n)
\]
to obtain
\begin{equation}\label{array:G9}
 (A,Q)\overset d=(a+F+D,R_{\mathrm{ord}}+E).
\end{equation}

\medskip
\noindent\textbf{Step 3.} It remains to identify the limit of the
infinitesimal remainder.
{
Set $Z_{n,i}=Y_{n,i}\ind_{\{i>m_n\}}$ for $1\le i\le n$.
The sums of $Z_{n,i}-\E\tau(Z_{n,i})$ and $Z_{n,i}^2$ give the two
coordinates of $(D_n,E_n)$.
Define
}\[
 \mu_n(dx)=w(x)\sum_i\mathcal L(Z_{n,i})(dx),
 \qquad \kappa_n=\sum_i(\E\tau(Z_{n,i}))^2.
\]
Since $\tau^2\le w$, Jensen's inequality and \eqref{array:G1} give
\[
 \kappa_n=\sum_i(\E\tau(Z_{n,i}))^2
 \le\sum_i\E\tau(Z_{n,i})^2
 \le\int\mu_n(dx)\le C.
\]
Also, by \eqref{array:G2}, for sufficiently large $L$,
\[
 \sup_n\mu_n(\{|x|>L\})
 \le\sup_n\sum_i\P(|Z_{n,i}|>L)
 \le-\log\inf_n\P(\widehat V_n^2\le L^2)
 \longrightarrow0.
\]
The total masses are bounded by $C$, so $(\mu_n)$ has a weakly
convergent subsequence. Since $0\le\kappa_n\le C$, pass to a further
subsequence on which $\mu_n\Rightarrow\mu$ and $\kappa_n\to\kappa$.
By \eqref{array:G7}, the array $\{Z_{n,i}\}$ is infinitesimal.
Together with \eqref{array:G8} and the limits $\mu_n\Rightarrow\mu$,
$\kappa_n\to\kappa$, Lemma~\ref{app:infinitesimal-remainder} gives
\begin{equation}\label{array:G11}
 (D,E)\overset d=(\sigma Z+U_{\mathrm{poi}}+\gamma_\nu,
                         \sigma^2+\kappa+R_{\mathrm{poi}}).
\end{equation}
The parameters in \eqref{array:G11} are
\[
 \begin{gathered}
 \nu(dx)=\ind_{\{x\ne0\}}w(x)^{-1}\mu(dx),\qquad
 \sigma^2=\mu(\{0\})-\kappa\ge0,\\
 \gamma_\nu=\int(\tau_0-\tau)\,d\nu.
 \end{gathered}
\]
Since $1\wedge x^2\le2w(x)$ and $|\tau_0(x)-\tau(x)|\le w(x)$,
\[
 \int(1\wedge x^2)\nu(dx)\le2\mu(\mathbb R\setminus\{0\}),
 \qquad
 \int|\tau_0-\tau|\,d\nu\le\mu(\mathbb R\setminus\{0\})<\infty.
\]
Use independent $Z$ and $\mathcal N$ to realize \eqref{array:G11}.
Choose them independent of $(Y_j)_{j\ge1}$. This preserves
the law in \eqref{array:G9}.

{Substituting \eqref{array:G11} and \eqref{array:ordinary-centering-identity}
into \eqref{array:G9}, and setting
}\[
 b=a+\gamma_\nu+\sum_j\E[(\tau_0-\tau)(Y_j)],
\]
gives \eqref{array:joint-canonical-form}.
\end{proof}

The next two propositions give the main steps that reduce the limit
representation \eqref{array:joint-canonical-form} to its Gaussian part.
Their proofs are given in Sections~\ref{sec:rigidity-start}
and~\ref{array:gaussian-exclusion}, respectively. Some supporting lemma
proofs are given in Supplementary Sections~\ref{supp:tail-band-tools}
and~\ref{app:positive-gaussian-energy}.

\begin{proposition}
\label{main:first-eliminations}
Under the assumptions of Proposition~\ref{array:joint-representation},
suppose that
\[
 0<Q<\infty\quad\text{almost surely},\qquad
 A/\sqrt Q\sim\N(0,1).
\]
Then
\[
 \sigma>0,\qquad \kappa=0.
\]
\end{proposition}

\begin{proposition}
\label{arr:rigidity}
Let $(A,Q)$ and its representation be as in
Proposition~\ref{array:joint-representation}. Suppose that
\[
 0<Q<\infty\quad\text{almost surely},\qquad
 A/\sqrt Q\sim\N(0,1),\qquad \sigma>0,\quad\kappa=0.
\]
Then, in this representation,
\[
 Y_j=0\quad\text{almost surely for every }j,\qquad
 \nu=0,\qquad b=0.
\]
Consequently,
\[
 (A,Q)\overset d=(\sigma Z,\sigma^2),\qquad Z\sim\N(0,1).
\]
\end{proposition}

\begin{proof}[Proof of Theorem~\ref{main:array}]
Assume \eqref{hist:cond:N} and choose any medians $b_n^2$ of $V_n^2$.
Proposition~\ref{hist:prop:median} gives $b_n>0$ eventually.
Discard the finite initial segment on which $b_n$ may vanish, and put
\[
 Y_{n,i}=X_{n,i}/b_n,\qquad
 A_n=S_n/b_n,\qquad \widehat V_n=V_n/b_n.
\]
By Proposition~\ref{hist:prop:median}, every subsequence has a further
subsequence along which
\[
 (A_n,\widehat V_n^2)\dto(A,Q),\qquad
 0<Q<\infty\quad\text{almost surely},\qquad
 A/\sqrt Q\sim\N(0,1).
\]
Fix such a subsequence. Proposition~\ref{array:joint-representation}
gives the representation \eqref{array:joint-canonical-form}.
By Proposition~\ref{main:first-eliminations}, its parameters satisfy
$\sigma>0$ and $\kappa=0$.
Proposition~\ref{arr:rigidity} therefore gives
$(A,Q)\overset d=(\sigma Z,\sigma^2)$.
The median normalization determines the remaining variance.
Along the fixed subsequence, $\widehat V_n^2\pto\sigma^2$. Since
$b_n^2$ is a median,
\[
 \P(\widehat V_n^2\le1)\ge\tfrac12,
 \qquad \P(\widehat V_n^2\ge1)\ge\tfrac12.
\]
If $\sigma^2>1$, then $\P(\widehat V_n^2\le1)\to0$; if
$\sigma^2<1$, then $\P(\widehat V_n^2\ge1)\to0$.
Each case contradicts the corresponding lower bound $1/2$,
so $\sigma^2=1$. Every convergent subsequence therefore has the
same limit, and joint tightness gives
\begin{equation}\label{arr:denominator-finish}
 (A_n,\widehat V_n^2)\dto(Z,1),\qquad
 \widehat V_n^2\pto1,\qquad A_n\dto\N(0,1),
\end{equation}
which proves \eqref{main:normalization}.

To prove \eqref{main:tail}, we also use the vanishing
individual-coordinate limits in Proposition~\ref{arr:rigidity}.
Recall $w$ from \eqref{def:truncations}.
We first prove $\max_{i\le n}\E w(Y_{n,i})\to0$.
If this fails, choose deterministic indices $i_n$ attaining the
maximum and pass to a subsequence on which
\[
 \liminf_{n\to\infty}\E w(Y_{n,i_n})>0.
\]
Since $|Y_{n,i_n}|\le\widehat V_n$, the selected coordinates are
tight. Pass to a further subsequence on which $Y_{n,i_n}\dto Y$.
Apply the construction in the proof of
Proposition~\ref{array:joint-representation} with $Y_{n,i_n}$ in
the first position. It represents the joint limit $(Z,1)$ in
\eqref{arr:denominator-finish} with $Y_1\overset d=Y$.
Since $w$ is bounded and continuous,
\[
 \E w(Y_1)=\lim_{n\to\infty}\E w(Y_{n,i_n})>0.
\]
Propositions~\ref{main:first-eliminations} and~\ref{arr:rigidity}
force $Y_1=0$ in this representation, a contradiction.
Consequently,
\begin{equation}\label{main:derived-infinitesimality}
 \max_{i\le n}\E w(Y_{n,i})\longrightarrow0.
\end{equation}
Since
$1-e^{-y^2}\le\min\{y^2,1\}\le2y^2/(1+y^2)=2w(y)$, it follows from
\eqref{main:derived-infinitesimality} that
\[
 \max_{i\le n}\{1-\E e^{-Y_{n,i}^2}\}
 \le2\max_{i\le n}\E w(Y_{n,i})\longrightarrow0.
\]
Fix $\varepsilon>0$. For all sufficiently large $n$ and every
$i\le n$, it follows that
\begin{align} \label{arr:Eexponential}
      \E e^{-Y_{n,i}^2}\ge\frac{1+e^{-\varepsilon^2}}2.
\end{align}
By \eqref{arr:denominator-finish} and the boundedness of $e^{-v}$ and $e^{-2v}$,
\[
 \E e^{-\widehat V_n^2}\longrightarrow e^{-1},\qquad
 \E e^{-2\widehat V_n^2}\longrightarrow e^{-2}.
\]
Rowwise independence therefore gives
\begin{equation}\label{arr:exponential-variance-sum}
 \begin{aligned}
 1+\sum_{i=1}^n\operatorname{Var}(e^{-Y_{n,i}^2})
 &\le\prod_{i=1}^n
   \left(1+
   \frac{\operatorname{Var}(e^{-Y_{n,i}^2})}
        {(\E e^{-Y_{n,i}^2})^2}\right)
 =\prod_{i=1}^n
   \frac{\E e^{-2Y_{n,i}^2}}{(\E e^{-Y_{n,i}^2})^2}
 =\frac{\E e^{-2\widehat V_n^2}}
          {(\E e^{-\widehat V_n^2})^2}
 \longrightarrow1,
 \end{aligned}
\end{equation}
which implies $\sum_{i=1}^{n}\operatorname{Var}(e^{-Y_{n,i}^{2}})\to 0$.
Combining \eqref{arr:exponential-variance-sum} and \eqref{arr:Eexponential}
with Chebyshev's inequality gives
\begin{align} \label{arr:sumPY}
\begin{split}
 \sum_{i=1}^n\P(|Y_{n,i}|>\varepsilon)
 &=\sum_{i=1}^{n}\P\left(  e^{-Y_{n,i}^{2}}-\E e^{-Y_{n,i}^2}< e^{-\varepsilon^{2}}-\E e^{-Y_{n,i}^2} \right)\\
 &\le \sum_{i=1}^{n}\P\left(  \left|e^{-Y_{n,i}^{2}}-\E e^{-Y_{n,i}^2}\right|\ge \frac{1-e^{-\varepsilon^{2}}}{2}\right) \\
 &\le\frac{4}{(1-e^{-\varepsilon^2})^2}
       \sum_{i=1}^n\operatorname{Var}(e^{-Y_{n,i}^2})
 \longrightarrow0.
\end{split}
\end{align}
Substituting $Y_{n,i}=X_{n,i}/b_n$ proves \eqref{main:tail}.

It remains to prove
\eqref{main:maximal}. On $\{\widehat V_n^2\ge1/2\}$, the event
$\{\max_{i\le n}|X_{n,i}|>\varepsilon V_n\}$ implies
$\{\max_{i\le n}|Y_{n,i}|>\varepsilon/\sqrt2\}$.
Hence a union bound together with \eqref{arr:sumPY} and $\hat{V}_{n}^{2}\pto 1$ gives
\[
 \P\!\left(\max_{i\le n}|X_{n,i}|>\varepsilon V_n\right)
 \le\P(\widehat V_n^2<1/2)
      +\sum_{i=1}^n\P(|Y_{n,i}|>\varepsilon/\sqrt2)
 \longrightarrow0.
\]
When $V_n=0$, every coordinate is zero, so the event on the left
is empty under the stated convention.
This proves \eqref{main:maximal} and completes the theorem.
\end{proof}

\subsection{Proof of Theorem~\ref{main:shao-characterization}}
\label{sec:proof-characterization}

\begin{proof}[Proof of Theorem~\ref{main:shao-characterization}]
With the stated rowwise replacements,
\citet[Theorems~1.1 and~1.2, p.~1742]{Shao2018} show that, under
\eqref{main:shao-C1}, the normal limit \eqref{hist:cond:N} is equivalent
to \eqref{main:shao-C2} and~\eqref{main:shao-C3}.
This proves sufficiency. Conversely, \eqref{hist:cond:N} implies
\eqref{main:shao-C1} by Theorem~\ref{main:array}; Shao's results then
give \eqref{main:shao-C2} and~\eqref{main:shao-C3}.
\end{proof}

\subsection{Proof of Corollary~\ref{cor:fixed-sequence}}
\label{sec:proof-fixed-sequence}

\begin{proof}[Proof of Corollary~\ref{cor:fixed-sequence}]
Apply Theorem~\ref{main:array} to $X_{n,i}=X_i$, $1\le i\le n$.
This gives $\max_{i\le n}|X_i|/V_n\pto0$.

Fix $j\ge1$. For $n\ge j$, with $|X_j|/V_n=0$ when $V_n=0$,
\[
 0\le\frac{|X_j|}{V_n}
 \le\frac{\max_{i\le n}|X_i|}{V_n}\pto0.
\]
Since $V_n$ is nondecreasing, $|X_j|/V_n$ has an almost sure limit.
This limit is zero because $|X_j|/V_n\pto0$. Thus
\[
 V_n\longrightarrow\infty
 \quad\text{almost surely on }\{X_j\ne0\}.
\]
The assertions $V_n\to\infty$ on $\{X_j\ne0\}$ hold simultaneously
for all $j$ outside a single null set. Hence $V_n\to\infty$
almost surely whenever at least one coordinate is nonzero.

It remains to exclude the event that all coordinates vanish.
This event is contained in $\{T_n=0\}$ for every $n$.
By the Portmanteau theorem applied to the closed set $\{0\}$,
\[
 \P(X_i=0\text{ for every }i)
 \le\limsup_{n\to\infty}\P(T_n=0)=0,
\]
because the standard normal law has no atom at zero.
Therefore $\sum_{i=1}^\infty X_i^2=\infty$ almost surely.

To prove \eqref{cor:finite-series-ratio}, suppose that $\sum_{i=1}^\infty X_i$
converges almost surely and $0<\sum_{i=1}^\infty X_i^2<\infty$
almost surely. Then
\[
 T_n\longrightarrow
 \frac{\sum_{i=1}^\infty X_i}
      {\left(\sum_{i=1}^\infty X_i^2\right)^{1/2}}
 \quad\text{almost surely}.
\]
If $\lim_{n\to\infty}T_n$ had the standard normal law, then
$T_n\dto\N(0,1)$. The divergence argument for $V_n$ would give
$\sum_{i=1}^\infty X_i^2=\infty$ almost surely, contradicting
the assumed finiteness.
\end{proof}

\section{Proof of Proposition~\ref{main:first-eliminations}}
\label{sec:rigidity-start}
\label{main:auxiliary-facts}
\label{app:annuli}
Assume the hypotheses of Proposition~\ref{main:first-eliminations} and use
the representation in Proposition~\ref{array:joint-representation}.
Define the residual pair
\begin{equation}\label{def:residual}
 U=b+U_{\mathrm{ord}}+U_{\mathrm{poi}},\qquad
 R=R_{\mathrm{ord}}+R_{\mathrm{poi}}.
\end{equation}
\proofanchor
\label{main:residual-tool}

The following lemma shows that, for every fixed $d>0$,
$U/\sqrt{d+R}$ has lighter tails than any
nondegenerate Gaussian. By comparison, a standard normal variable $G$ satisfies
$\E e^{\lambda G^2}<\infty$ if and only if $\lambda<1/2$.
We use this difference in tail decay to obtain the contradictions below.
\begin{lemma}
\label{main:residual-integrability}
For the pair $(U,R)$ in \eqref{def:residual},
\begin{equation}\label{rigidity:residual-offset}
\E\exp\left\{\frac{\lambda U^2}{d+R}\right\}<\infty
 \qquad\text{for every }d>0\text{ and }\lambda>0.
\end{equation}
\end{lemma}

\label{supp:residual-integrability-proof}

\begin{proof}
Split $U$ and $R$ at magnitude one:
\[
 U=b+U_1+U_2,\qquad R=R_1+R_2,
\]
where, by \eqref{def:components},
\begin{equation}\label{def:residual-size-split}
\begin{aligned}
 U_1
 &=\sum_j\left\{Y_j\ind_{\{|Y_j|\le1\}}
       -\E[Y_j\ind_{\{|Y_j|\le1\}}]\right\}
   +\int_{|y|\le1}y(\mathcal N-\nu)(dy),\\
 U_2
 &=\sum_jY_j\ind_{\{|Y_j|>1\}}
   +\int_{|y|>1}y\mathcal N(dy),\\
 R_1
 &=\sum_jY_j^2\ind_{\{|Y_j|\le1\}}
   +\int_{|y|\le1}y^2\mathcal N(dy),\\
 R_2
 &=\sum_jY_j^2\ind_{\{|Y_j|>1\}}
   +\int_{|y|>1}y^2\mathcal N(dy).
\end{aligned}
\end{equation}

Let
\begin{equation}\label{arr:offset-large-count}
 N_2=\sum_j\ind_{\{|Y_j|>1\}}+\mathcal N(\{y:|y|>1\}).
\end{equation}
The conditions \eqref{array:component-summability} in
Proposition~\ref{array:joint-representation} give
\begin{equation}\label{arr:offset-large-count-mean}
 \E N_2=\sum_j\P(|Y_j|>1)+\nu(|y|>1)<\infty.
\end{equation}
{Thus $N_2<\infty$ almost surely. By
\eqref{def:residual-size-split} and \eqref{arr:offset-large-count},
the Cauchy--Schwarz inequality gives
}\begin{equation}\label{arr:offset-large-cauchy}
 U_2^2\le N_2R_2.
\end{equation}
{For $R_2>0$, the weighted Cauchy--Schwarz inequality gives
\[
 \frac{(b+U_1+U_2)^2}{(d+R_1)+R_2}
 \le\frac{(b+U_1)^2}{d+R_1}+\frac{U_2^2}{R_2}.
\]
If $R_2=0$, then $U_2=0$, so the same bound holds with the last
quotient set to zero. Thus \eqref{arr:offset-large-cauchy} gives
}\begin{equation}\label{arr:offset-pointwise-reduction}
 \frac{U^2}{d+R}
 \le\frac{(b+U_1)^2}{d+R_1}+N_2
 \le\frac{2b^2}{d}+\frac{2U_1^2}{d+R_1}+N_2.
\end{equation}
{Applying the Cauchy--Schwarz inequality to the exponential of
\eqref{arr:offset-pointwise-reduction}, we obtain
}\begin{equation}\label{arr:offset-moment-reduction}
 \E e^{\lambda U^2/(d+R)}
 \le e^{2\lambda b^2/d}
 \left\{\E e^{4\lambda U_1^2/(d+R_1)}\right\}^{1/2}
 \left\{\E e^{2\lambda N_2}\right\}^{1/2}
 \qquad(d,\lambda>0).
\end{equation}
It remains to bound the two expectations in \eqref{arr:offset-moment-reduction}.
Independence and the monotone convergence theorem applied to
\eqref{arr:offset-large-count}, together with
\eqref{arr:offset-large-count-mean}, give
\begin{equation}\label{arr:offset-large-integrability}
\begin{aligned}
 \E e^{vN_2}
 &=e^{(e^v-1)\nu(|y|>1)}
   \prod_j\bigl[1+(e^v-1)\P(|Y_j|>1)\bigr]\\
 &\le\exp\!\left\{(e^v-1)
   \left[\sum_j\P(|Y_j|>1)+\nu(|y|>1)\right]\right\}
 <\infty\qquad(v>0).
\end{aligned}
\end{equation}

To estimate $\E e^{4\lambda U_1^2/(d+R_1)}$, put $Z_j=Y_j\ind_{\{|Y_j|\le1\}}$.
By \eqref{array:component-summability},
\begin{equation}\label{arr:offset-small-energy}
 \E R_1=\sum_j\E Z_j^2+\int_{|y|\le1}y^2\nu(dy)<\infty.
\end{equation}
Also, $|Z_j|\le1$ implies, for every integer $m\ge1$,
\begin{equation}\label{arr:offset-mean-bound}
 \sum_{j\le m}|\E Z_j|\le\sum_{j\le m}\E|Z_j|\le m.
\end{equation}
Fix $a>0$. Since $tZ_j-at^2Z_j^2\le1/(4a)$,
\[
 \log\E e^{t(Z_j-\E Z_j)-at^2Z_j^2}
 \le\frac1{4a}-t\E Z_j\le\frac1{4a}+|t||\E Z_j|.
\]
Summing and using \eqref{arr:offset-mean-bound} gives
\begin{align}\label{arr:offset-ordinary-estimates-1}
     \sum_{j\le m}\log\E e^{t(Z_j-\E Z_j)-at^2Z_j^2}
 &\le\frac{m}{4a}+|t|\sum_{j\le m}|\E Z_j|
 \le\frac{m}{4a}+m|t|.
\end{align}
For every $a>0$, there is a finite constant $C_a$ such that
\begin{equation}\label{arr:offset-elementary}
\begin{aligned}
 e^{u-au^2}&\le1+u+C_au^2.
\end{aligned}
\end{equation}
Indeed, $(e^{u-au^2}-1-u)/u^2$ extends continuously to zero and is
bounded above on $\R$.
{It follows from \eqref{arr:offset-elementary} and the inequality $\log v\le v-1$ for $v>0$ that, for $j>m$,
}\begin{equation}\label{arr:offset-ordinary-estimates-2}
\begin{aligned}
 \log\E e^{t(Z_j-\E Z_j)-at^2Z_j^2}
 &\le-t\E Z_j+\E(e^{tZ_j-at^2Z_j^2}-1)\\
 &\le -t\E Z_j+t\E Z_j+C_{a}t^{2}\E Z_j^{2}\\
 &= C_at^2\E Z_j^2.
\end{aligned}
\end{equation}
For the small-jump contribution to $tU_1-at^2R_1$, first
restrict the jumps to $\epsilon<|y|\le1$, where $0<\epsilon<\eta<1$.
The intensity on this set is finite since
$\nu(\epsilon<|y|\le1)\le\epsilon^{-2}\int_{|y|\le1}y^2\nu(dy)$.
The Poisson exponential formula \citep[Proposition~3.8]{LP2018} gives
\begin{equation}\label{arr:offset-poisson-transform}
\begin{aligned}
 &\log\E\exp\Bigl\{
 t\int_{\epsilon<|y|\le1}y(\mathcal N-\nu)(dy)
 -at^2\int_{\epsilon<|y|\le1}y^2\mathcal N(dy)\Bigr\}\\
 &\quad=-t\int_{\epsilon<|y|\le1}y\nu(dy)
 +\log\E\exp\Bigl\{
 \int_{\epsilon<|y|\le1}(ty-at^2y^2)\mathcal N(dy)\Bigr\}\\
 &\quad=\int_{\epsilon<|y|\le1}
       (e^{ty-at^2y^2}-1-ty)\nu(dy).
\end{aligned}
\end{equation}
For $|y|\le\eta$, use \eqref{arr:offset-elementary}.
For $\eta<|y|\le1$, use $ty-at^2y^2\le1/(4a)$. This gives
\begin{equation}\label{arr:offset-poisson-estimate}
\begin{aligned}
 &\int_{\epsilon<|y|\le1}(e^{ty-at^2y^2}-1-ty)\nu(dy)\\
 &\quad\le C_at^2\int_{|y|\le\eta}y^2\nu(dy)
       +\int_{\eta<|y|\le1}(e^{1/(4a)}+1+|t||y|)\nu(dy)\\
 &\quad\le C_at^2\int_{|y|\le\eta}y^2\nu(dy)
       +C_{a,\eta}(1+|t|),
\end{aligned}
\end{equation}
Here $C_{a,\eta}$ is finite because, on $\eta<|y|\le1$,
$1\le\eta^{-2}y^2$ and $|y|\le\eta^{-1}y^2$.
It is independent of $\epsilon,t$.
For finite ordinary sums and Poisson jumps restricted away from zero,
independence makes the logarithmic exponential moments additive.
Add \eqref{arr:offset-ordinary-estimates-1},
\eqref{arr:offset-ordinary-estimates-2}, and \eqref{arr:offset-poisson-estimate},
then let the number of ordinary terms increase and the jump cutoff
decrease to zero. Almost sure convergence and the uniform bounds
allow us to apply Fatou's lemma:
\begin{equation}\label{arr:offset-transform}
 \log\E e^{tU_1-at^2R_1}
 \le C_{a,m,\eta}(1+|t|)
 +C_at^2\left\{\sum_{j>m}\E Z_j^2
                  +\int_{|y|\le\eta}y^2\nu(dy)\right\}.
\end{equation}
Fix $\delta>0$. By \eqref{arr:offset-small-energy}, choose $m$ large
and $\eta$ small so that the coefficient of $t^2$ in
\eqref{arr:offset-transform} is at most $\delta/2$.
The remaining linear term is bounded by
\begin{equation}\label{arr:offset-linear-absorption}
 C(1+|t|)\le\frac\delta2t^2+C+\frac{C^2}{2\delta}
 \qquad(C\ge0,\ t\in\mathbb R),
\end{equation}
Combining \eqref{arr:offset-transform} and
\eqref{arr:offset-linear-absorption} gives
\begin{equation}\label{arr:offset-subquadratic}
 \E e^{tU_1-at^2R_1}
 \le K_{a,\delta}e^{\delta t^2}
 \qquad(t\in\mathbb R),\qquad K_{a,\delta}<\infty.
\end{equation}

For arbitrary $d,\lambda_0>0$, take
$a=1/(4\lambda_0)$ and $\delta=d/(8\lambda_0)$ in
\eqref{arr:offset-subquadratic}.
For $u\in\R$ and $v\ge0$, completing the square gives
\[
 tu-\frac{(d+v)t^2}{4\lambda_0}
 =-\frac{d+v}{4\lambda_0}
       \left(t-\frac{2\lambda_0u}{d+v}\right)^2
   +\frac{\lambda_0u^2}{d+v}.
\]
Integrating in $t$ and then applying Tonelli's theorem yields
the Gaussian-mixture identity
\citep[proof of Theorem~2.1]{DKL2004},
\begin{equation}\label{arr:offset-gaussian-mixture}
\begin{aligned}
 \E[(d+R_1)^{-1/2}e^{\lambda_0U_1^2/(d+R_1)}]
 &\quad=\frac1{\sqrt{4\pi\lambda_0}}\int_{\mathbb R}
 e^{-dt^2/(4\lambda_0)}
 \E e^{tU_1-t^2R_1/(4\lambda_0)}\,dt\\
 &\quad\le\frac{K_{a,\delta}}{\sqrt{4\pi\lambda_0}}
 \int_{\mathbb R}e^{-dt^2/(8\lambda_0)}\,dt<\infty.
\end{aligned}
\end{equation}
For $d,\lambda>0$, the Cauchy--Schwarz inequality and
\eqref{arr:offset-gaussian-mixture}, with $\lambda_0=2\lambda$, give
\begin{equation}\label{arr:offset-small-integrability}
\begin{aligned}
 \E e^{\lambda U_1^2/(d+R_1)}
 \le\left\{\E[(d+R_1)^{-1/2}
       e^{2\lambda U_1^2/(d+R_1)}]\right\}^{1/2}
       \left\{\E\sqrt{d+R_1}\right\}^{1/2}
 <\infty,
\end{aligned}
\end{equation}
since $\E\sqrt{d+R_1}\le\sqrt{d+\E R_1}<\infty$.

Finally, substitute \eqref{arr:offset-small-integrability}, with
exponent $4\lambda$, and \eqref{arr:offset-large-integrability}, with
$v=2\lambda$, into \eqref{arr:offset-moment-reduction}.
This proves \eqref{rigidity:residual-offset} for every $d,\lambda>0$.
\end{proof}

The representation \eqref{array:joint-canonical-form} gives
\begin{equation}\label{rigidity:representation}
 A=\sigma Z+U,\qquad Q=\sigma^2+\kappa+R,\qquad
 T=\frac{A}{\sqrt Q},
\end{equation}
where $R\ge0$, $Z$ is independent of $(U,R)$, and $T\sim\N(0,1)$.
We first exclude the case $\sigma=0$, in which the ratio is $T=U/\sqrt Q$.
Lemma~\ref{main:residual-integrability} gives lighter-than-Gaussian
tails for $U/\sqrt{d+R}$ at every fixed $d>0$.
We therefore examine far-tail events to seek a contradiction with
the standard normal law of $T$.

\proofanchor
\label{main:radial-construction}
\proofanchor
\label{main:tail-energy-tools}
A standard approach to studying tail events is to change measure so
that they carry more probability. To study $Q$ on tail bands, we use quadratic exponential weights as in
\citet[(1.4) and the proof of Theorem~2.1]{DKL2004}.
We write $\P$ for the original law of this representation and $\E$
for expectation under $\P$. Expectations under other measures carry
the corresponding measure or tilt parameters as subscripts. For $x>2$ and
$\varsigma\in\{-1,1\}$, put $z_x=\sqrt{1-x^{-2}}$ and set
\begin{equation}\label{main:radial-definition}
 \Gamma_x^\varsigma(dr,d\omega)
 =x^{-2}rQe^{\varsigma z_xrA-r^2Q/2}\,dr\,\P(d\omega),\qquad r>0.
\end{equation}
In the tail-band estimates, a $+$ or $-$ subscript or superscript
abbreviates the choice $\varsigma=1$ or $\varsigma=-1$, respectively.
The variable $r$ combines different tilting scales.
With $u=r\sqrt Q$, we have $rQ\,dr=u\,du$, so integrating out $r$
leaves a density depending only on $T$.

We consider the bands
{\begin{equation}\label{main:tail-bands}
 I_x^\varsigma=\{x\le\varsigma T\le 5x/4\}.
\end{equation}
}

The next lemma bounds the sample marginal density and shows that
$\Gamma_x^\varsigma(I_x^\varsigma)$ has a positive limit. These bounds
transfer tail estimates from $\P$ to the conditional laws of $Q$ below.

\begin{lemma}
\label{main:radial-law}
For $x>2$ and either sign $\varsigma$, $\Gamma_x^\varsigma$ is a
probability law. Its sample marginal is
\begin{equation}\label{main:radial-likelihood}
 \ell_x(\varsigma T)\,d\P,\qquad
 \ell_x(v)=x^{-2}\int_0^\infty u e^{z_xvu-u^2/2}\,du,
\end{equation}
and there is a constant $C$ such that
\begin{equation}\label{main:radial-bound}
 0\le\ell_x(v)\le Cx^{-2}(1+|v|)e^{v^2/2}.
\end{equation}
Moreover, the band probability $p_x=\Gamma_x^\varsigma(I_x^\varsigma)>0$ is the
same for both signs and satisfies
\begin{equation}\label{main:band-positive-mass}
 p_x\longrightarrow p_*>0\qquad(x\to\infty).
\end{equation}
\end{lemma}

The proof is given in Supplementary Section~\ref{app:radial-law-proof}.

For Borel $D\subset(0,\infty)$, define
\begin{equation}\label{def:band-law}
 H_x^\varsigma(D)=\Gamma_x^\varsigma(Q\in D\mid I_x^\varsigma).
\end{equation}
In particular, $H_x^+=H_x^{1}$ and $H_x^-=H_x^{-1}$.
We write $H_x^\varsigma(Q\le u)$ for $H_x^\varsigma((0,u])$,
and similarly for other events expressed in terms of $Q$.

Intuitively, $H_x^\varsigma$ records which values of $Q$ contribute
to the $\Gamma_x^\varsigma$-mass of the tail band $I_x^\varsigma$.
If $\sigma=0$, then $T=U/\sqrt{\kappa+R}$.
On events where $Q=\kappa+R$ is bounded away from zero, the
exponential moment bound in \eqref{rigidity:residual-offset} gives
lighter-than-Gaussian tails for this ratio.
To account for the standard normal tails of $T$, the weighted band
mass should therefore come from increasingly small values of $Q$
as $x$ increases.
This motivates studying whether the medians of $H_x^\varsigma$
can tend to zero.

Choose a median $m_x$ of $H_x^+$ for each $x>2$; explicitly,
\begin{equation}\label{def:band-median}
 H_x^+((0,m_x])\ge\tfrac12,\qquad
 H_x^+([m_x,\infty))\ge\tfrac12.
\end{equation}
Since $0<Q<\infty$ almost surely, $0<m_x<\infty$.

The next lemma says that these medians cannot collapse to zero along
geometric tail levels. 

\begin{lemma}
\label{main:noncollapse}
{Under the assumptions of Proposition~\ref{main:first-eliminations},
let $m_x$ be any choice of medians satisfying \eqref{def:band-median}.
For every fixed sufficiently large $x_0$,
}{\[
 m_{(9/8)^jx_0}\longrightarrow m_*(x_0)\in(0,\infty)
 \qquad (j\to\infty).
\]
}\end{lemma}

The proof is given in Supplementary Section~\ref{supp:noncollapse-proof}.

In the proof below, $\sigma=0$ will instead force
$m_x\to0$, giving the required contradiction. Once $\sigma>0$ is known, conditioning
on $(U,R)$ and integrating over $Z$ excludes $\kappa>0$.

\proofanchor
\label{main:first-eliminations-proof}
\begin{proof}[Proof of Proposition~\ref{main:first-eliminations}]
\label{main:sigma-positive}
To prove $\sigma>0$, suppose that $\sigma=0$, so $A=U$ and $Q=\kappa+R$.
Fix $q_0>0$. On $\{Q\ge q_0\}$,
\[
 q_0+R\le Q+R=\kappa+2R\le2Q,
 \qquad T^2=\frac{U^2}{Q}\le\frac{2U^2}{q_0+R}.
\]
Since
$(1+|v|)e^{v^2/2}\le Ce^{v^2}$, Lemma~\ref{main:residual-integrability}
with $(d,\lambda)=(q_0,2)$ gives
\begin{equation}\label{rigidity:localized-integrability}
 \E[(1+|T|)e^{T^2/2}\ind_{\{Q\ge q_0\}}]<\infty.
\end{equation}
By \eqref{main:radial-bound} and
\eqref{rigidity:localized-integrability},
\begin{equation}\label{rigidity:weighted-energy-vanishes}
\begin{aligned}
 \Gamma_x^\varsigma(Q\ge q_0)
 &=\E[\ell_x(\varsigma T)\ind_{\{Q\ge q_0\}}]
 \le Cx^{-2}\E[(1+|T|)e^{T^2/2}\ind_{\{Q\ge q_0\}}]
 \longrightarrow0.
\end{aligned}
\end{equation}
Combining \eqref{rigidity:weighted-energy-vanishes} with
$p_x\to p_*>0$ gives
\[
 H_x^+([q_0,\infty))
 \le\Gamma_x^+(Q\ge q_0)/p_x\longrightarrow0.
\]
If $m_x\ge q_0$, then
$H_x^+([q_0,\infty))\ge H_x^+([m_x,\infty))\ge1/2$.
Thus $m_x<q_0$ eventually. Since $q_0>0$ is arbitrary,
$m_x\to0$. In particular, $m_{(9/8)^jx_0}\to0$ for every fixed
sufficiently large $x_0$, contradicting Lemma~\ref{main:noncollapse}.
Therefore $\sigma>0$.

\label{main:kappa-zero}
To prove $\kappa=0$, suppose that $\kappa>0$ and choose
{\[
 \lambda=\frac12+\frac{\kappa}{8\sigma^2}>\frac12.
\]
}{Since $Q-2\lambda\sigma^2=3\kappa/4+R>0$,
\[
 \frac{\lambda(\sigma z+U)^2}{Q}-\frac{z^2}{2}
 =-\frac{Q-2\lambda\sigma^2}{2Q}
   \left(z-\frac{2\lambda\sigma U}{Q-2\lambda\sigma^2}\right)^2
   +\frac{\lambda U^2}{Q-2\lambda\sigma^2}.
\]
Integrating against $(2\pi)^{-1/2}\,dz$ gives
}{\begin{equation}\label{rigidity:gaussian-square-integral}
 \E[e^{\lambda T^2}\mid U,R]
 =\sqrt{\frac{\sigma^2+\kappa+R}{3\kappa/4+R}}
   \exp\!\left\{\frac{\lambda U^2}{3\kappa/4+R}\right\}.
\end{equation}
}{Here $3\kappa/4+R\ge(\kappa+R)/2$, and the square-root factor is
at most $\sqrt{2(1+\sigma^2/\kappa)}$. Consequently,
}\[
 \E e^{\lambda T^2}
 \le \sqrt{2(1+\sigma^2/\kappa)}\,
       \E\exp\!\left\{\frac{2\lambda U^2}{\kappa+R}\right\}<\infty
\]
by Lemma~\ref{main:residual-integrability} with $d=\kappa$.
This contradicts $T\sim\N(0,1)$, for which this expectation is infinite
when $\lambda>1/2$. Hence $\kappa=0$.
\end{proof}

\section{Proof of Proposition~\ref{arr:rigidity}}
\label{array:gaussian-exclusion}

We use the representation from Proposition~\ref{array:joint-representation}.
Proposition~\ref{main:first-eliminations} gives $\sigma>0$ and $\kappa=0$, so
\begin{equation}\label{array:positive-gaussian-setup}
 A=\sigma Z+U,\qquad Q=\sigma^2+R,\qquad
 T=\frac{A}{\sqrt Q}\sim\N(0,1),
\end{equation}
{where $Z$ is independent of $(U,R)$ in \eqref{def:residual}.
Write $\varphi$ for the standard normal density. Conditional on $(U,R)$,
}{\begin{equation}\label{array:conditional-normal}
 \begin{gathered}
 T\mid(U,R)\sim\N\!\left(\frac U{\sqrt Q},\frac{\sigma^2}{Q}\right),\quad
 f_{T\mid U,R}(s)
 =\frac{\sqrt Q}{\sigma}
   \varphi\!\left(\frac{s\sqrt Q-U}{\sigma}\right).
 \end{gathered}
\end{equation}
}{We use the standard normal law of $T$ through the identity
\[
 \E f_{T\mid U,R}(s)=\varphi(s)
 \qquad\text{for Lebesgue-almost every }s\in\R.
\]
Since $\varphi(s)>0$, the density ratio
$f_{T\mid U,R}(s)/\varphi(s)$ has expectation one.
For $\varsigma\in\{-1,1\}$ and $t>0$, expansion gives
}{\begin{equation}\label{array:conditional-density-ratio}
 \begin{aligned}
 w_{\varsigma,t}
 &:=\frac{f_{T\mid U,R}(\varsigma t)}{\varphi(t)}
 =\frac{\sqrt Q}{\sigma}
   \exp\left\{\frac{t^2}{2}
        -\frac{(\varsigma t\sqrt Q-U)^2}{2\sigma^2}\right\}
=e^{\varsigma(t/\sigma)U-t^2R/(2\sigma^2)}
    \mathcal C_{\varsigma,t},
 \end{aligned}
\end{equation}
}where
{\begin{equation}\label{array:density-correction}
 \mathcal C_{\varsigma,t}
 =\sqrt{1+R/\sigma^2}\,
  \exp\left\{\frac{\varsigma tU}{\sigma}
          \bigl(\sqrt{1+R/\sigma^2}-1\bigr)
                    -\frac{U^2}{2\sigma^2}\right\}.
\end{equation}
}{The averaging identity above gives, for Lebesgue-almost every $t>0$,
}{\begin{equation}\label{array:normal-density-normalization}
 1=\E w_{\varsigma,t}
  =\E\!\left[
     e^{\varsigma(t/\sigma)U-t^2R/(2\sigma^2)}
     \mathcal C_{\varsigma,t}\right].
\end{equation}
}{We isolate the exponential factor in
\eqref{array:normal-density-normalization} to use independence.
Its exponent is additive over the coordinates and the Poisson component.
Writing $r=t/\sigma$, define the quadratic transforms
}\begin{equation}\label{def:critical-transform}
 L_\varsigma(r)=\E e^{\varsigma rU-r^2R/2}
              =\E e^{\varsigma rA-r^2Q/2},
 \qquad \varsigma\in\{-1,1\},\quad r>0.
\end{equation}
{The Gaussian part contributes exactly one, since
$\E e^{\varsigma r\sigma Z-r^2\sigma^2/2}=1$.
When $U=R=0$, both transforms equal one.
}

\proofanchor
\label{sec:residual-transform-deficit}

The following lemma shows that any surviving ordinary coordinate or
Poisson component leaves a strict asymptotic gap below one.

{
\begin{lemma}
\label{lem:strict-transform-deficit}
Let $(U,R)$ be the pair in \eqref{def:residual},
with components satisfying \eqref{array:component-summability}, and
let $L_\varsigma(r)=\E e^{\varsigma rU-r^2R/2}$, as in
\eqref{def:critical-transform}. Then
\[
 0<L_+(r)L_-(r)\le1\qquad(r>0).
\]
If $\P(Y_j\ne0)>0$ for some $j$, or if $\nu\ne0$, then
\begin{equation}\label{array:strict-transform-gap}
 \limsup_{r\to\infty}\sqrt{L_+(r)L_-(r)}<1.
\end{equation}
\end{lemma}
\begin{proof}
For $r>0$, put
\[
 g_r(u)=\frac{e^{ru-r^2u^2/2}+e^{-ru-r^2u^2/2}}2
       =e^{-r^2u^2/2}\cosh(ru).
\]
The identities $(\log\cosh 0)=(\log\cosh)'(0)=0$ and
$(\log\cosh v)''\le1$ give $\log\cosh v\le v^2/2$.
For $u\ne0$, $g_r(u)\le e^{r|u|-r^2u^2/2}\to0$. Hence
\begin{equation}\label{array:scalar-kernel-bounds}
 0<g_r(u)\le1,\qquad
 g_r(u)\longrightarrow\ind_{\{u=0\}}
 \quad\text{as }r\to\infty.
\end{equation}
The arithmetic--geometric mean inequality gives
\begin{equation}\label{array:ordinary-transform-bound}
 \sqrt{\E e^{rY_j-r^2Y_j^2/2}\,
              \E e^{-rY_j-r^2Y_j^2/2}}
 \le\frac{\E e^{rY_j-r^2Y_j^2/2}
               +\E e^{-rY_j-r^2Y_j^2/2}}2
 =\E g_r(Y_j)\le1.
\end{equation}
Both expectations are positive and finite, since
$e^{\pm ru-r^2u^2/2}\le e^{1/2}$.

To use independence, truncate to the first $N$ coordinates and
Poisson jumps with $|u|>N^{-1}$:
\begin{align*}
 U^{(N)}
 &=b+\sum_{j\le N}
       \left\{Y_j-\E[Y_j\ind_{\{|Y_j|\le1\}}]\right\}\\
 &\quad+\int_{|u|>N^{-1}}u\,\mathcal N(du)
       -\int_{N^{-1}<|u|\le1}u\,\nu(du),\\
 R^{(N)}
 &=\sum_{j\le N}Y_j^2
       +\int_{|u|>N^{-1}}u^2\,\mathcal N(du),\qquad N\ge1.
\end{align*}
By \eqref{array:component-summability}, the truncated jump region has
finite intensity and a finite compensating integral.
Proposition~\ref{array:joint-representation} therefore gives
\[
 U^{(N)}\longrightarrow U,\qquad R^{(N)}\uparrow R
 \quad\text{almost surely}.
\]
Fix $r>0$ and put
\[
 L_\varsigma^{(N)}(r)
 =\E e^{\varsigma rU^{(N)}-r^2R^{(N)}/2}.
\]
Independence and the Poisson exponential formula give
\[
 \begin{aligned}
 \log L_\varsigma^{(N)}(r)
 &=\varsigma rb+
   \sum_{j\le N}\left\{
    \log\E e^{\varsigma rY_j-r^2Y_j^2/2}
    -\varsigma r\E[Y_j\ind_{\{|Y_j|\le1\}}]\right\}\\
 &\quad+\int_{|u|>N^{-1}}
   \left\{e^{\varsigma ru-r^2u^2/2}-1
          -\varsigma ru\ind_{\{|u|\le1\}}\right\}\nu(du).
 \end{aligned}
\]
Adding the two signs cancels all terms linear in $\varsigma$.
By \eqref{array:ordinary-transform-bound} and
\eqref{array:scalar-kernel-bounds},
\begin{equation}\label{array:truncated-transform-product}
 \begin{aligned}
 \sqrt{L_+^{(N)}(r)L_-^{(N)}(r)}
 &=\prod_{j\le N}
    \sqrt{\E e^{rY_j-r^2Y_j^2/2}\,
                 \E e^{-rY_j-r^2Y_j^2/2}}\times\exp\left\{
       \int_{|u|>N^{-1}}(g_r(u)-1)\nu(du)\right\}
 \le1.
 \end{aligned}
\end{equation}
Fatou's lemma gives
$0<L_\varsigma(r)\le\liminf_N L_\varsigma^{(N)}(r)$ for each sign.
Consequently,
\[
 \begin{aligned}
 0<L_+(r)L_-(r)
 &\le\left(\liminf_N L_+^{(N)}(r)\right)
       \left(\liminf_N L_-^{(N)}(r)\right)
 \le\liminf_N\{L_+^{(N)}(r)L_-^{(N)}(r)\}\le1.
 \end{aligned}
\]

If $\P(Y_j\ne0)>0$ for a fixed $j$, retain its factor in
\eqref{array:truncated-transform-product} for every $N\ge j$.
All remaining factors are at most one. Fatou's lemma at fixed $r$ gives
\[
 \sqrt{L_+(r)L_-(r)}\le\E g_r(Y_j).
\]
The bounded convergence theorem and \eqref{array:scalar-kernel-bounds} yield
\[
 \limsup_{r\to\infty}\sqrt{L_+(r)L_-(r)}
 \le\P(Y_j=0)<1.
\]

If $\nu\ne0$, the sets $\{|u|>a\}$ increase to
$\R\setminus\{0\}$ as $a\downarrow0$, so choose a fixed $a>0$ with
$0<\nu(\{|u|>a\})<\infty$.
For $N^{-1}<a$, bound the ordinary factors in
\eqref{array:truncated-transform-product} by one and retain only
the Poisson integral over $|u|>a$, using $g_r-1\le0$.
Fatou's lemma then gives
\[
 \sqrt{L_+(r)L_-(r)}
 \le\exp\left\{\int_{|u|>a}(g_r(u)-1)\nu(du)\right\}.
\]
Applying the bounded convergence theorem on this fixed finite-intensity set gives
\[
 \limsup_{r\to\infty}\sqrt{L_+(r)L_-(r)}
 \le e^{-\nu(\{|u|>a\})}<1.
\]
This proves \eqref{array:strict-transform-gap}.
\end{proof}
}

\proofanchor
\label{sec:normal-transform-lower}

Lemma~\ref{lem:strict-transform-deficit} gives
$\limsup_{r\to\infty}\sqrt{L_+(r)L_-(r)}<1$ if an ordinary factor or
Poisson component survives. To contradict this bound, we seek lower
bounds approaching one for both transforms at the same arguments.
The expectation in \eqref{array:normal-density-normalization},
without its factor $\mathcal C_{\varsigma,t}$, is $L_\varsigma(t/\sigma)$.
We therefore need to control the effect of $\mathcal C_{\varsigma,t}$
in the tail-band averages.

With $u=r\sqrt Q$, the joint density of $(T,u)$ under
$\Gamma_x^\varsigma$ is
\[
 \begin{aligned}
 \frac{\Gamma_x^\varsigma(T\in dt,u\in du)}{dt\,du}
 &=x^{-2}u e^{\varsigma z_xut-u^2/2}\varphi(t)
 =x^{-2}u e^{-u^2/(2x^2)}\varphi(t-\varsigma z_xu),
 \qquad u>0,
 \end{aligned}
\]
so
\begin{equation}\label{app:radial-independent-laws}
 \mathsf E=T-\varsigma z_xr\sqrt Q\sim\N(0,1)
 \qquad\text{under }\Gamma_x^\varsigma.
\end{equation}
Under \eqref{array:positive-gaussian-setup}, define
\begin{equation}\label{def:residual-noise}
 v_x^\varsigma=\E_{\Gamma_x^\varsigma}\frac RQ,\qquad
 \xi_x^\varsigma=\frac{U-\varsigma z_xrR}{\sqrt Q}.
\end{equation}
Substituting
\eqref{array:positive-gaussian-setup} into
\eqref{app:radial-independent-laws} gives
\[
 \begin{aligned}
 \mathsf E
 &=\frac{\sigma Z+U-\varsigma z_xr(\sigma^2+R)}{\sqrt Q}
 =\xi_x^\varsigma+
   \frac{\sigma(Z-\varsigma\sigma z_xr)}{\sqrt Q}.
 \end{aligned}
\]
For fixed $(r,U,R)$, the conditional density of $Z$ under
$\Gamma_x^\varsigma$ is
\[
 f^{\Gamma_x^\varsigma}_{Z\mid r,U,R}(z)
 =\frac{\varphi(z)e^{\varsigma\sigma z_xrz}}
        {e^{\sigma^2z_x^2r^2/2}}
 =\varphi(z-\varsigma\sigma z_xr).
\]
Thus $\xi_x^\varsigma$ is the conditional mean of $\mathsf E$
given $(r,U,R)$, and the other term has conditional variance
$\sigma^2/Q$. The next lemma bounds this conditional mean in $L^2$
and quantifies the concentration of $Q$ near $\sigma^2$ on the bands.

\begin{lemma}
\label{lem:gaussian-energy-improvement}
Under \eqref{array:positive-gaussian-setup}, for every $x>2$ and either sign,
\begin{equation}\label{array:exact-residual-noise}
 \E_{\Gamma_x^\varsigma}(\xi_x^\varsigma)^2
 =\E_{\Gamma_x^\varsigma}\frac RQ=v_x^\varsigma.
\end{equation}
{There are constants $\beta\in(2/3,1)$, $C<\infty$ and $x_0>2$
such that, for every $x\ge x_0$ and either sign,
}{\begin{equation}\label{array:improved-rates}
 \begin{gathered}
 H_x^\varsigma\!\left([\sigma^2,\sigma^2+Cx^{-\beta}]\right)
       \ge1-Cx^{-\beta},
 v_x^\varsigma\le Cx^{-\beta}.
 \end{gathered}
\end{equation}
}\end{lemma}

The proof is given in Supplementary Section~\ref{app:positive-gaussian-energy}.

The first bound in \eqref{array:improved-rates} controls
$R=Q-\sigma^2$ on $I_x^\varsigma$. The second bound and
\eqref{array:exact-residual-noise} give
$\E_{\Gamma_x^\varsigma}(\xi_x^\varsigma)^2\le Cx^{-\beta}$.
Eliminating $\varsigma z_xr$ from \eqref{app:radial-independent-laws}
and \eqref{def:residual-noise} gives
\begin{equation}\label{array:posterior-noise}
 U-\frac{TR}{\sqrt Q}
 =\sqrt Q\xi_x^\varsigma-\frac{\mathsf E R}{\sqrt Q}.
\end{equation}
We now seek lower bounds for both transforms at a common argument
to contradict \eqref{array:strict-transform-gap}.

{
\begin{lemma}
\label{lem:common-transform-lower}
Under \eqref{array:positive-gaussian-setup}, the transforms
\eqref{def:critical-transform} satisfy $0<L_\varsigma(r)<\infty$ for
every $r>0$ and either sign. There is a deterministic choice
$t_x\in[x,5x/4]$ for every sufficiently large $x$ such that
\begin{equation}\label{array:two-transform-lower}
 \liminf_{x\to\infty}L_\varsigma(t_x/\sigma)\ge1,
 \qquad \varsigma=\pm1.
\end{equation}
\end{lemma}
\begin{proof}
Lemma~\ref{lem:strict-transform-deficit} gives finiteness and positivity.
We first express $L_\varsigma$ in terms of
$\mathcal C_{\varsigma,t}$.
Let $\Pi_{\varsigma,t}$ be the conditional law of $(U,R)$ given
$T=\varsigma t$ under $\P$.
Bayes' formula and \eqref{array:conditional-density-ratio} give density
$w_{\varsigma,t}$ relative to the original law of $(U,R)$.
Hence, for almost every $t>0$,
\begin{equation}\label{array:posterior-transform}
 L_\varsigma(t/\sigma)
 =\E[w_{\varsigma,t}\mathcal C_{\varsigma,t}^{-1}]
 =\E_{\Pi_{\varsigma,t}}\mathcal C_{\varsigma,t}^{-1}.
\end{equation}
For each sign, fix a Lebesgue-null set $N_\varsigma$ outside which
the density formula for $\Pi_{\varsigma,t}$ and
\eqref{array:posterior-transform} hold.
Write $a_+=\max\{a,0\}$. For $t\notin N_\varsigma$, put
\begin{equation}\label{array:transform-deficit}
 \begin{aligned}
 d_\varsigma(t)
 &:=\E_{\Pi_{\varsigma,t}}
      [(1-\mathcal C_{\varsigma,t}^{-1})_+]\in[0,1],\\
 1-L_\varsigma(t/\sigma)
 &=\E_{\Pi_{\varsigma,t}}[1-\mathcal C_{\varsigma,t}^{-1}]
 \le d_\varsigma(t).
 \end{aligned}
\end{equation}
We set $d_\varsigma=0$ on $N_\varsigma$.
It is enough to find $t_x$ such that $d_+(t_x)+d_-(t_x)\to0$.

Fix $x$ and a sign, and write
\begin{equation}\label{def:conditional-band-law}
 \widehat P=\Gamma_x^\varsigma(\,\cdot\mid I_x^\varsigma),
 \qquad
 \mu_x(dt)=\frac{\varphi(t)\ell_x(t)}{p_x}
       \ind_{[x,5x/4]}(t)\,dt.
\end{equation}
Under $\widehat P$, the marginals of $Q$ and $\varsigma T$ are
$H_x^\varsigma$ and $\mu_x$. The measure $\mu_x$ is the same for both
signs because $\varphi$ is even.
Because the weight $\ell_x(\varsigma T)$ depends only on $T$,
the conditional law of $(U,R)$ given $T=\varsigma t$ remains
$\Pi_{\varsigma,t}$.
Thus, for every nonnegative measurable $g$,
\begin{equation}\label{array:band-posterior-mixture}
 \E_{\widehat P}g(\varsigma T,U,R)
 =\int_x^{5x/4}
    \E_{\Pi_{\varsigma,t}}g(t,U,R)\,\mu_x(dt).
\end{equation}

It remains to bound the mean of
$(1-\mathcal C_{\varsigma,\varsigma T}^{-1})_+$ under $\widehat P$.
Let $\beta\in(2/3,1)$ be the exponent in Lemma~\ref{lem:gaussian-energy-improvement}.
Since $\mathsf E$ is standard normal under $\Gamma_x^\varsigma$ and
$p_x\to p_*>0$ by \eqref{main:band-positive-mass},
\eqref{array:exact-residual-noise} and \eqref{array:improved-rates} give
\begin{equation}\label{array:band-noise-moments}
 \E_{\widehat P}(\xi_x^\varsigma)^2
 \le\frac{v_x^\varsigma}{p_x}\le Cx^{-\beta},
 \qquad
 \E_{\widehat P}\mathsf E^2\le\frac1{p_x}\le C.
\end{equation}
The concentration bound in \eqref{array:improved-rates} gives
$K<\infty$ such that
\begin{equation}\label{array:band-small-residual}
 G_x:=\{R\le Kx^{-\beta}\},\qquad
 \widehat P(G_x^c)\le Cx^{-\beta}.
\end{equation}
Recall from \eqref{array:positive-gaussian-setup} that $Q=\sigma^2+R$.
Since $\sigma>0$ and $Kx^{-\beta}\to0$, on $G_x$ we have, for all
sufficiently large $x$,
\begin{equation}\label{array:band-energy-bounds}
 \sigma^2\le Q\le\sigma^2+Kx^{-\beta}\le2\sigma^2.
\end{equation}

To substitute for $U$ in \eqref{array:density-correction}, put
$t=\varsigma T\in[x,5x/4]$ and use \eqref{array:posterior-noise} to define
\begin{equation}\label{def:posterior-error}
 \begin{aligned}
 \zeta
 &=\frac1\sigma\left(U-\frac{TR}{\sqrt Q}\right)
 =\frac{\sqrt Q}{\sigma}\xi_x^\varsigma
   -\frac{R}{\sigma\sqrt Q}\mathsf E,\qquad
 \frac U\sigma&=\frac{\varsigma tR}{\sigma\sqrt Q}+\zeta.
 \end{aligned}
\end{equation}
On $G_x$, \eqref{def:posterior-error} and
\eqref{array:band-energy-bounds} imply
$\zeta^2\le C\{(\xi_x^\varsigma)^2+x^{-2\beta}\mathsf E^2\}$.
Equation~\eqref{array:band-noise-moments} and the Cauchy--Schwarz
inequality give
\begin{equation}\label{array:band-posterior-error-moments}
 \begin{aligned}
 \E_{\widehat P}[\zeta^2\ind_{G_x}]
 &\le C\E_{\widehat P}(\xi_x^\varsigma)^2
      +Cx^{-2\beta}\E_{\widehat P}\mathsf E^2
 \le Cx^{-\beta},\\
 \E_{\widehat P}[|\zeta|\ind_{G_x}]
 &\le\bigl(\E_{\widehat P}[\zeta^2\ind_{G_x}]\bigr)^{1/2}
 \le Cx^{-\beta/2}.
 \end{aligned}
\end{equation}

Put $q=\sqrt{1+R/\sigma^2}$. Substitution of
\eqref{def:posterior-error} into \eqref{array:density-correction} gives
\[
 \begin{aligned}
 \log\mathcal C_{\varsigma,t}
 &=\log q+\varsigma t
   \left(\frac{\varsigma tR}{\sigma^2q}+\zeta\right)(q-1)
   -\frac12\left(\frac{\varsigma tR}{\sigma^2q}+\zeta\right)^2\\
 &=\log q+\frac{t^2R(q-1)^2}{2\sigma^2q^2}
       -\frac{\varsigma t\zeta(q-1)}q-\frac{\zeta^2}{2}.
 \end{aligned}
\]
The coefficients of $t^2$ and $\varsigma t\zeta$ are simplified using
\[
 \frac R{\sigma^2}(1-q^{-1})-\frac{R^2}{2\sigma^4q^2}
 =\frac{R(q-1)^2}{2\sigma^2q^2},
 \qquad
 q-1-\frac R{\sigma^2q}=-\frac{q-1}{q}.
\]
On $G_x$, $q-1=R/(\sigma^2(q+1))\le CR$ and $0\le\log q\le CR$.
It follows that
\begin{equation}\label{array:density-correction-bound}
 |\log\mathcal C_{\varsigma,t}|
 \le C\{R+t^2R^3+tR|\zeta|+\zeta^2\}
 \qquad\text{on }G_x.
\end{equation}
Using $t\le5x/4$, \eqref{array:band-small-residual}, and
\eqref{array:band-posterior-error-moments}, we obtain
\begin{equation}\label{array:posterior-correction-small}
 \begin{aligned}
 \E_{\widehat P}
    [|\log\mathcal C_{\varsigma,\varsigma T}|\ind_{G_x}]
  &\le C\bigl\{x^{-\beta}+x^{2-3\beta}
       +x^{1-\beta}\E_{\widehat P}[|\zeta|\ind_{G_x}]
       +\E_{\widehat P}[\zeta^2\ind_{G_x}]\bigr\}\\
  &\le C\{x^{-\beta}+x^{2-3\beta}+x^{1-3\beta/2}\}
  \le Cx^{1-3\beta/2}.
 \end{aligned}
\end{equation}
Since $\beta>2/3$, $x^{1-3\beta/2}\to0$.

For every $c>0$,
\begin{equation}\label{array:inverse-correction-loss}
 0\le(1-c^{-1})_+\le\min\{1,|\log c|\}.
\end{equation}
By \eqref{array:band-small-residual},
\eqref{array:posterior-correction-small}, and
\eqref{array:inverse-correction-loss},
\begin{equation}\label{array:averaged-transform-deficit}
 \begin{aligned}
 \E_{\widehat P}
      [(1-\mathcal C_{\varsigma,\varsigma T}^{-1})_+]
 &\le\widehat P(G_x^c)
   +\E_{\widehat P}
      [|\log\mathcal C_{\varsigma,\varsigma T}|\ind_{G_x}]
 \le Cx^{-\beta}+Cx^{1-3\beta/2}
 \le Cx^{1-3\beta/2}\longrightarrow0.
 \end{aligned}
\end{equation}

Equations~\eqref{array:band-posterior-mixture} and
\eqref{array:averaged-transform-deficit} give
$\int d_\varsigma(t)\,\mu_x(dt)\le Cx^{1-3\beta/2}$.
The same measure $\mu_x$ appears for both signs, so
\[
 a_x:=\int[d_+(t)+d_-(t)]\,\mu_x(dt)
 \le Cx^{1-3\beta/2}\longrightarrow0.
\]
Markov's inequality gives
\[
 \mu_x\bigl(d_+(t)+d_-(t)\le a_x+x^{-1}\bigr)
 \ge\frac{x^{-1}}{a_x+x^{-1}}>0.
\]
Since $\mu_x(N_+\cup N_-)=0$, we can choose a deterministic
$t_x\in[x,5x/4]\setminus(N_+\cup N_-)$ satisfying
\[
 d_+(t_x)+d_-(t_x)\le a_x+x^{-1}.
\]
For each sign, \eqref{array:transform-deficit} now yields
\[
 L_\varsigma(t_x/\sigma)
 \ge1-d_\varsigma(t_x)\ge1-a_x-x^{-1}.
\]
Since $a_x\to0$, this proves \eqref{array:two-transform-lower}.
\end{proof}
}

\phantomsection
\label{array:rigidity-proof}
We can now compare the upper and lower bounds at the same parameter.
This excludes the ordinary coordinates and the Poisson component;
the standard normal law of $T$ then determines the remaining drift.

{
\begin{proof}[Proof of Proposition~\ref{arr:rigidity}]
Use the setup in \eqref{array:positive-gaussian-setup}.
By Lemma~\ref{lem:strict-transform-deficit},
$0<L_+(r)L_-(r)\le1$ for every $r>0$.
Lemma~\ref{lem:common-transform-lower} supplies a common choice
$t_x\in[x,5x/4]$ such that
\[
 \liminf_{x\to\infty}L_\varsigma(t_x/\sigma)\ge1,
 \qquad\varsigma=\pm1.
\]
Thus
\[
 \sqrt{L_+(t_x/\sigma)L_-(t_x/\sigma)}\longrightarrow1.
\]
If $\P(Y_j\ne0)>0$ for some $j$, or if $\nu\ne0$,
\eqref{array:strict-transform-gap} contradicts
$\sqrt{L_+(t_x/\sigma)L_-(t_x/\sigma)}\to1$,
since $t_x/\sigma\ge x/\sigma\to\infty$. Therefore
\[
 \P(Y_j=0)=1\quad\text{for every }j,\qquad \nu=0.
\]
By \eqref{def:components} and \eqref{def:residual}, $R=0$ and $U=b$
almost surely. Hence $T=Z+b/\sigma$.
Since both $T$ and $Z$ have the standard normal law,
$0=\E T=\E Z+b/\sigma=b/\sigma$, so $b=0$.
Substituting $U=R=0$ into \eqref{array:positive-gaussian-setup} gives
\[
 (A,Q)\overset d=(\sigma Z,\sigma^2),\qquad \sigma>0.
\]
\end{proof}
}

\section*{Acknowledgements}
Liu S.H. was partially supported by the Fundamental Research Funds for the Central Universities DUT25RC(3)133.
Shao Q.M. was partially supported by National Nature Science Foundation of China NSFC 12031005 and Shenzhen Outstanding Talents Training Fund, China.

\bmhead{AI Statement}
During the preparation of this manuscript, the authors used ChatGPT for English
editing, LaTeX formatting, and checks of notation, references and cross-references.
ChatGPT also helped organize the exposition and check intermediate calculations. These
checks covered the joint-limit representation, quadratic tilting and entropy estimates
in the necessity proof of the self-normalized central limit theorem. The authors
reviewed and revised the material and take full responsibility for the mathematical
arguments and conclusions.

\Needspace{8\baselineskip}
\phantomsection

\clearpage
\begin{center}
{\Large\bfseries Supplement to ``Characterizing standard normal limits of self-normalized sums''}
\end{center}
\medskip
\setcounter{section}{0}
\renewcommand{\thesection}{S.\arabic{section}}
\renewcommand{\theHsection}{supp.\arabic{section}}
\renewcommand{\theequation}{\thesection.\arabic{equation}}
\renewcommand{\thetheorem}{\thesection.\arabic{theorem}}

\section{Proof of Lemma~\ref{app:infinitesimal-remainder}}
\label{app:infinitesimal}

Let $\{Z_{n,i}\}$ be the rowwise independent infinitesimal array
in Lemma~\ref{app:infinitesimal-remainder}. Recall from
\eqref{def:truncations} that
\[
 w(x)=\frac{x^2}{1+x^2},\qquad
 \tau(x)=\frac{x}{1+x^2},\qquad
 \tau_0(x)=x\ind_{\{|x|\le1\}}.
\]
Recall also that
\[
 \mu_n(dx)=w(x)\sum_i\mathcal L(Z_{n,i})(dx),\qquad
 \kappa_n=\sum_i[\E\tau(Z_{n,i})]^2,
\]
where $\mathcal L(X)$ denotes the law of $X$.
The hypotheses give $\mu_n\Rightarrow\mu$ as finite measures and
$\kappa_n\to\kappa<\infty$. Write
\[
 \nu(dx)=\ind_{\{x\ne0\}}w(x)^{-1}\mu(dx).
\]

\begin{proof}
Weak convergence gives $\sup_n\mu_n(\R)<\infty$. Since
$(1\wedge x^2)/w(x)\le2$ for $x\ne0$,
\[
 \int(1\wedge x^2)\nu(dx)\le2\mu(\R)<\infty,
\]
so $\nu$ is a L\'evy measure. Write
$D_n=\sum_i\{Z_{n,i}-\E\tau(Z_{n,i})\}$ and $E_n=\sum_iZ_{n,i}^2$.
Fix $(s,t)\ne(0,0)$ and set
\[
 H(x)=sx+tx^2,\qquad c_n=s\sum_i\E\tau(Z_{n,i}).
\]
Since $H$ is continuous and $H(0)=0$, the independent row variables
$H(Z_{n,i})$ are infinitesimal
\citep[Chapter~2, Definition~3.3, p.~54]{AG1980}, and
\[
 \sum_iH(Z_{n,i})-c_n=sD_n+tE_n.
\]
To prove convergence of $sD_n+tE_n$, we will use the theorem of
\citet[Chapter~2, Theorem~4.7, p.~61]{AG1980}.
With the truncation used there \citep[p.~51]{AG1980}, write
\[
 H_{n,i,\delta}=H(Z_{n,i})\ind_{\{|H(Z_{n,i})|\le\delta\}}.
\]
Define a measure on
$\R\setminus\{0\}$ by
\[
 \lambda(B)=\nu(\{x:H(x)\in B\}),
 \qquad B\subset\R\setminus\{0\}.
\]
We exclude zero from the image measure because $H$ may have a nonzero root.
Since $1\wedge H(x)^2=O(x^2)$ near zero and is bounded elsewhere,
$\lambda$ is a L\'evy measure.
For this infinitesimal array, we verify the following sufficient
conditions of the theorem:
\begin{enumerate}[label=(\roman*),leftmargin=*,itemsep=2pt]
\item The truncated variance sums satisfy
\[
 \lim_{\delta\downarrow0}\limsup_n
       \sum_i\operatorname{Var}(H_{n,i,\delta})
 =\lim_{\delta\downarrow0}\liminf_n
       \sum_i\operatorname{Var}(H_{n,i,\delta})=v
\]
for some finite $v\ge0$.
\item For every bounded
continuous $g$ and each $\delta>0$ with $\lambda(\{-\delta,\delta\})=0$,
\[
 \sum_i\E\big[g(H(Z_{n,i}))\ind_{\{|H(Z_{n,i})|>\delta\}}\big]
 \longrightarrow\int_{|y|>\delta}g(y)\lambda(dy).
\]
\item For each such
$\delta$, there is $\alpha_\delta\in\R$ such that
\[
 \sum_i\E H_{n,i,\delta}-c_n\longrightarrow\alpha_\delta.
\]
\end{enumerate}

For (ii), fix $\delta>0$ with $\lambda(\{-\delta,\delta\})=0$
and a bounded continuous function $g$. Since the quotient
$g(H(x))\ind_{\{|H(x)|>\delta\}}/w(x)$, extended by zero at the origin,
is bounded, its possible discontinuities lie in $\{|H|=\delta\}$. The definition of $\lambda$ gives
$\nu(\{|H|=\delta\})=\lambda(\{-\delta,\delta\})=0$.
Since $H(0)=0$ and $\delta>0$, this set excludes zero, so
\[
 \mu(\{|H|=\delta\})
 =\int_{\{|H(x)|=\delta\}}w(x)\nu(dx)=0.
\]
The Portmanteau theorem for finite measures, applied to
$\mu_n\Rightarrow\mu$, then gives
\[
 \begin{aligned}
 \sum_i\E\big[g(H(Z_{n,i}))\ind_{\{|H(Z_{n,i})|>\delta\}}\big]
 &=\int\frac{g(H(x))\ind_{\{|H(x)|>\delta\}}}{w(x)}\,\mu_n(dx)\\
 &\longrightarrow\int\frac{g(H(x))\ind_{\{|H(x)|>\delta\}}}{w(x)}\,\mu(dx)
 =\int_{|y|>\delta}g(y)\lambda(dy).
 \end{aligned}
\]

For (iii), keep this $\delta$ and consider
\[
 q_\delta(x)=
 \frac{H(x)\ind_{\{|H(x)|\le\delta\}}-s \tau(x)}{w(x)}
 \quad(x\ne0),\qquad q_\delta(0)=t,
\]
which is bounded and continuous at zero. Its remaining possible discontinuities also lie
in the $\mu$-null set $\{|H|=\delta\}$,
so the same form of the Portmanteau theorem yields
\begin{equation}\label{app:projected-center}
 \begin{aligned}
 \sum_i\E H_{n,i,\delta}-c_n
 &=\int q_\delta\,d\mu_n\longrightarrow\alpha_\delta:=t\mu(\{0\})+\int\big\{H(x)\ind_{\{|H(x)|\le\delta\}}-s \tau(x)\big\}\nu(dx).
 \end{aligned}
\end{equation}
The integral is absolutely convergent, by boundedness of $q_\delta$.

For (i), first fix any $\delta>0$.
Boundedness of $q_\delta$ gives
\[
 \sum_i\big|\E H_{n,i,\delta}-s\E\tau(Z_{n,i})\big|
 \le\|q_\delta\|_\infty\mu_n(\R)=O(1).
\]
Both functions being averaged are bounded and tend to zero at the
origin, so infinitesimality gives
\[
 \max_i\bigl(|\E H_{n,i,\delta}|+|s\E\tau(Z_{n,i})|\bigr)
 \longrightarrow0.
\]
Factoring the difference of squares and using these two bounds gives
\begin{equation}\label{app:mean-square-match}
 \begin{aligned}
 \left|\sum_i(\E H_{n,i,\delta})^2-s^2\kappa_n\right|
 &\le\max_i\bigl(|\E H_{n,i,\delta}|+|s\E\tau(Z_{n,i})|\bigr)\times\sum_i\big|\E H_{n,i,\delta}-s\E\tau(Z_{n,i})\big|
 \longrightarrow0.
 \end{aligned}
\end{equation}

The bounded quotient
$H(x)^2\ind_{\{|H(x)|\le\delta\}}/w(x)$ extends continuously to zero
with value $s^2$. As in (ii), its discontinuities are $\mu$-null when
$\lambda(\{-\delta,\delta\})=0$, so the Portmanteau theorem gives
\begin{equation}\label{app:truncated-energy-match}
 \sum_i\E H_{n,i,\delta}^2
 \longrightarrow s^2\mu(\{0\})+
     \int_{|H(x)|\le\delta}H(x)^2\nu(dx).
\end{equation}
Put $\sigma^2=\mu(\{0\})-\kappa$, as in
Lemma~\ref{app:infinitesimal-remainder}.
Since $\lambda$ has at most countably many atoms, for each $\delta>0$
choose $\delta_-\in(\delta/2,\delta)$ and $\delta_+\in(\delta,2\delta)$
with $\lambda(\{-\delta_\pm,\delta_\pm\})=0$. Then
\[
 \sum_i\E H_{n,i,\delta_-}^2
 \le\sum_i\E H_{n,i,\delta}^2
 \le\sum_i\E H_{n,i,\delta_+}^2.
\]
By \eqref{app:truncated-energy-match} and
$\sum_i(\E H_{n,i,\delta})^2\to s^2\kappa$ from
\eqref{app:mean-square-match}, these inequalities yield
\[
 \begin{aligned}
 s^2\sigma^2
 &\le\liminf_n\sum_i\operatorname{Var}(H_{n,i,\delta})
 \le\limsup_n\sum_i\operatorname{Var}(H_{n,i,\delta})\\
 &\le s^2\sigma^2+\int_{|H(x)|\le2\delta}H(x)^2\nu(dx).
 \end{aligned}
\]
For $0<\delta\le1/2$, the last integrand is bounded by
$H(x)^2\ind_{\{|H(x)|\le1\}}\le C(1\wedge x^2)$.
This bound is $\nu$-integrable, so the dominated convergence theorem
and the preceding inequalities give
\[
 \lim_{\delta\downarrow0}\limsup_n
       \sum_i\operatorname{Var}(H_{n,i,\delta})
 =\lim_{\delta\downarrow0}\liminf_n
       \sum_i\operatorname{Var}(H_{n,i,\delta})
 =s^2\sigma^2.
\]
Taking $(s,t)=(1,0)$ gives
$\sigma^2=\mu(\{0\})-\kappa\ge0$.
Thus (i) holds with $v=s^2\sigma^2$.

For any $\delta>0$ with $\lambda(\{-\delta,\delta\})=0$, Theorem~4.7 of
\citet[Chapter~2, p.~61]{AG1980} now gives
\begin{equation}\label{app:projected-limit}
 sD_n+tE_n\dto
 \N(\alpha_\delta,s^2\sigma^2)*c_\delta\mathrm{Pois}\,\lambda,
\end{equation}
where $c_\delta\mathrm{Pois}\,\lambda$ is the centred generalized
Poisson law of \citet[Chapter~2, Definition~2.5 and (2.5), p.~49]{AG1980},
whose characteristic exponent at frequency $\theta\in\R$ can be
rewritten using the definition of $\lambda$ as
\begin{equation}\label{app:projected-poisson-exponent}
 \begin{aligned}
 &\int\bigl(e^{i\theta y}-1
     -i\theta y\ind_{\{|y|\le\delta\}}\bigr)\lambda(dy)=\int\bigl(e^{i\theta H(x)}-1
     -i\theta H(x)\ind_{\{|H(x)|\le\delta\}}\bigr)\nu(dx).
 \end{aligned}
\end{equation}
Write the joint characteristic function of $(D_n,E_n)$ as
\[
 \varphi_n(u,v):=\E e^{i(uD_n+vE_n)},\qquad (u,v)\in\R^2.
\]
By \eqref{app:projected-center}, \eqref{app:projected-poisson-exponent} and \eqref{app:projected-limit}, for every $\theta\in\R$, we have
\[
 \varphi_n(\theta s,\theta t)
 =\E e^{i\theta(sD_n+tE_n)}
 \longrightarrow e^{\Psi(\theta s,\theta t)},
\]
where
\begin{equation}\label{array:G10}
 \begin{aligned}
 \Psi(s,t)
 &:=-\tfrac12s^2\sigma^2+i\alpha_\delta
     +\int\bigl(e^{iH(x)}-1-iH(x)\ind_{\{|H(x)|\le\delta\}}\bigr)\nu(dx)\\
 &=-\tfrac12s^2\sigma^2+it\mu(\{0\})
     +\int\bigl(e^{isx+itx^2}-1-is\tau(x)\bigr)\nu(dx).
 \end{aligned}
\end{equation}
The integrand is
$O(x^2)$ near zero and bounded outside a neighborhood of zero, so the
integral is absolutely convergent, and we set $\Psi(0,0)=0$.
Since $(s,t)$ was arbitrary, taking $\theta=1$ yields
\[
 \varphi_n(s,t)\longrightarrow e^{\Psi(s,t)}
 \qquad\text{for every }(s,t)\in\R^2.
\]

It remains to verify that the vector in \eqref{app:remainder-limit}
is well defined and has joint characteristic function $e^{\Psi(s,t)}$.
Take a Poisson random measure $\mathcal N$ of intensity $\nu$ and an
independent standard normal $Z$. Recall the correction and Poisson components from
Lemma~\ref{app:infinitesimal-remainder} and \eqref{def:components}:
\[
 \begin{aligned}
 \gamma_\nu&=\int(\tau_0-\tau)\,d\nu,\\
 U_{\mathrm{poi}}
 &=\int_{|x|\le1}x(\mathcal N-\nu)(dx)
   +\int_{|x|>1}x\mathcal N(dx),\\
 R_{\mathrm{poi}}&=\int x^2\mathcal N(dx).
 \end{aligned}
\]
By definition, $|\tau_0(x)-\tau(x)|\le w(x)$ for all $x\in \R$, so
\[
 \begin{aligned}
 \int|\tau_0(x)-\tau(x)|\,\nu(dx)
 \le\int w(x)\,\nu(dx)
 =\mu(\R\setminus\{0\})<\infty.
 \end{aligned}
\]
Thus $\tau_0-\tau\in L^1(\nu)$, and $\gamma_\nu$ is well defined.
For the compensated small-jump term in $U_{\mathrm{poi}}$, note that
\[
 \nu(\epsilon<|x|\le1)
 \le\epsilon^{-2}\int_{|x|\le1}x^2\nu(dx)<\infty
 \qquad(0<\epsilon<1).
\]
The truncated integrals are therefore well defined. The moment formulas
for compensated Poisson integrals
\citep[equation~(12.5) and Lemma~12.2]{LP2018} give
\[
 \begin{aligned}
 \E\left[\int_{\epsilon<|x|\le1}x(\mathcal N-\nu)(dx)\right]&=0,\\
 \E\left|\int_{\epsilon<|x|\le1}x(\mathcal N-\nu)(dx)\right|^2
 &=\int_{\epsilon<|x|\le1}x^2\nu(dx)
 \le\int_{|x|\le1}x^2\nu(dx)<\infty.
 \end{aligned}
\]
Along $\epsilon=2^{-k}$, successive increments are compensated integrals
over disjoint annuli, so they are independent and centred. The second-moment
bound above makes this sequence an $L^2$-bounded martingale.
By the $L^2$ martingale convergence theorem, we may therefore define
\[
 \int_{|x|\le1}x(\mathcal N-\nu)(dx)
 :=\lim_{k\to\infty}\int_{2^{-k}<|x|\le1}x(\mathcal N-\nu)(dx)
 \qquad\text{in }L^2\text{ and almost surely}.
\]
For the small-jump part of $R_{\mathrm{poi}}$, Campbell's formula
\citep[Proposition~2.7]{LP2018} gives
\[
 \E\left[\int_{|x|\le1}x^2\mathcal N(dx)\right]
 =\int_{|x|\le1}x^2\nu(dx)<\infty.
\]
Since this random integral is nonnegative,
\[
 \int_{|x|\le1}x^2\mathcal N(dx)<\infty
 \qquad\text{almost surely}.
\]
For the large jumps, the same formula yields
\[
 \E\mathcal N(\{|x|>1\})
 =\nu(|x|>1)
 \le\int(1\wedge x^2)\nu(dx)<\infty.
\]
Hence $\mathcal N(\{|x|>1\})<\infty$ almost surely. The large-jump
integrals are therefore finite sums of finite terms, and
\[
 \int_{|x|>1}|x|\mathcal N(dx)
 +\int_{|x|>1}x^2\mathcal N(dx)<\infty
 \qquad\text{almost surely}.
\]
These estimates show that $U_{\mathrm{poi}}$ and $R_{\mathrm{poi}}$
are well defined and finite almost surely, and hence so is the vector
in \eqref{app:remainder-limit}.

We next compute the joint characteristic function of the Poisson components.
For $0<\epsilon<1$, write the truncated components as
\[
 \begin{aligned}
 U_\epsilon
 &:=\int_{|x|>\epsilon}x\mathcal N(dx)
     -\int_{\epsilon<|x|\le1}x\nu(dx),\\
 R_\epsilon&:=\int_{|x|>\epsilon}x^2\mathcal N(dx).
 \end{aligned}
\]
Since $\nu(|x|>\epsilon)<\infty$, conditioning on the Poisson number
of jumps \citep[Proposition~3.8]{LP2018} gives
\begin{equation}\label{app:truncated-poisson-characteristic}
 \begin{aligned}
 \E e^{isU_\epsilon+itR_\epsilon}
 &=\exp\left\{-is\int_{\epsilon<|x|\le1}x\nu(dx)\right\}
   \E\exp\left\{i\int_{|x|>\epsilon}(sx+tx^2)\mathcal N(dx)\right\}\\
 &=\exp\left\{-is\int_{\epsilon<|x|\le1}x\nu(dx)
       +\int_{|x|>\epsilon}(e^{isx+itx^2}-1)\nu(dx)\right\}\\
 &=\exp\left\{\int_{|x|>\epsilon}
       (e^{isx+itx^2}-1-is\tau_0(x))\nu(dx)\right\}.
 \end{aligned}
\end{equation}

Along $\epsilon=2^{-k}$, the preceding construction gives
\[
 U_\epsilon\longrightarrow U_{\mathrm{poi}},\qquad
 R_\epsilon\uparrow R_{\mathrm{poi}}
 \qquad\text{almost surely}.
\]
Since $|e^{isU_\epsilon+itR_\epsilon}|=1$, the bounded convergence theorem yields
\[
 \E e^{isU_\epsilon+itR_\epsilon}
 \longrightarrow\E e^{isU_{\mathrm{poi}}+itR_{\mathrm{poi}}}.
\]
To pass to the limit in the exponent, use
$|e^{iu}-1-iu|\le u^2/2$ for real $u$. For $|x|\le1$,
\[
 \begin{aligned}
 |e^{isx+itx^2}-1-is\tau_0(x)|
 &\le |e^{i(sx+tx^2)}-1-i(sx+tx^2)|+|t|x^2\\
 &\le\tfrac12(sx+tx^2)^2+|t|x^2\\
 &\le\left[\tfrac12(|s|+|t|)^2+|t|\right]x^2.
 \end{aligned}
\]
For $|x|>1$, $\tau_0(x)=0$ and $|e^{isx+itx^2}-1|\le2$.
These bounds are $\nu$-integrable because
$\int(1\wedge x^2)\nu(dx)<\infty$. The dominated convergence theorem therefore gives
\[
 \int_{|x|>\epsilon}(e^{isx+itx^2}-1-is\tau_0(x))\nu(dx)
 \longrightarrow
 \int(e^{isx+itx^2}-1-is\tau_0(x))\nu(dx).
\]
Taking limits on both sides of \eqref{app:truncated-poisson-characteristic} now yields
\[
 \E e^{isU_{\mathrm{poi}}+itR_{\mathrm{poi}}}
 =\exp\left\{\int
       (e^{isx+itx^2}-1-is \tau_0(x))\nu(dx)\right\}.
\]
Finally, independence of $Z$ and $\mathcal N$ shows that the vector
$(\sigma Z+U_{\mathrm{poi}}+\gamma_\nu,
\sigma^2+\kappa+R_{\mathrm{poi}})$ has joint characteristic exponent
\[
 \begin{aligned}
 &-\tfrac12s^2\sigma^2+it(\sigma^2+\kappa)+is\gamma_\nu+\int(e^{isx+itx^2}-1-is\tau_0(x))\,\nu(dx)
 =\Psi(s,t),
 \end{aligned}
\]
where we used $\sigma^2+\kappa=\mu(\{0\})$ and
$\gamma_\nu=\int(\tau_0-\tau)\,d\nu$. The multivariate L\'evy continuity theorem now
proves \eqref{app:remainder-limit}.
\end{proof}

\section{Proofs of Lemmas~\ref{main:radial-law} and~\ref{main:noncollapse}}
\label{supp:tail-band-tools}
\label{arr:product-interface}
\label{main:tilt-tools}

\subsection{Proof of Lemma~\ref{main:radial-law}}

\label{app:radial-law-proof}

Let $\P$ be the original law, with $0<Q<\infty$ almost surely and
$T=A/\sqrt Q\sim\N(0,1)$. As in the main text, $\E$ denotes
expectation under $\P$. Write $\varphi$ for the standard normal density.
For $x>2$ and $\varsigma\in\{-1,1\}$, recall
$z_x=\sqrt{1-x^{-2}}$ and the measure in \eqref{main:radial-definition}:
\[
 \Gamma_x^\varsigma(dr,d\omega)
 =x^{-2}rQe^{\varsigma z_xrA-r^2Q/2}\,dr\,\P(d\omega),
 \qquad r>0.
\]
The bands and the positive-band mass are
\[
 I_x^\varsigma=\{x\le\varsigma T\le5x/4\},\qquad
 p_x=\Gamma_x^+(I_x^+).
\]
Superscripts $+$ and $-$ denote the two choices of $\varsigma$.

The following lemma includes all the conclusions of
Lemma~\ref{main:radial-law}.

\begin{lemma}
\label{app:radial-identities}
Assume that $0<Q<\infty$ almost surely and $T=A/\sqrt Q\sim\N(0,1)$.
For $x>2$ and either sign $\varsigma$, the measure
$\Gamma_x^\varsigma$ in \eqref{main:radial-definition} is a probability law
with the following properties.
\begin{enumerate}[label=(\roman*),leftmargin=*,itemsep=2pt]
\item
The sample marginal of $\Gamma_x^\varsigma$ is
$\ell_x(\varsigma T)\,d\P$, where
\begin{equation}\label{app:radial-likelihood-explicit}
 \ell_x(v)=x^{-2}\int_0^\infty u e^{z_xvu-u^2/2}\,du=x^{-2}
 \{1+z_xv\sqrt{2\pi}\,e^{z_x^2v^2/2}\Phi(z_xv)\},
\end{equation}
and $\Phi$ is the standard normal distribution function.
There is a constant $C$, independent of $x$ and $v$, such that
\[
 0\le\ell_x(v)\le Cx^{-2}(1+|v|)e^{v^2/2}.
\]
\item
Under $\Gamma_x^\varsigma$, the residual
\begin{equation}\label{app:radial-residual}
 \mathsf E=T-\varsigma z_xr\sqrt Q
\end{equation}
has the standard normal law.

\item
The sample marginal in (i) yields the same positive band probability
for both signs:
\[
 p_x=\Gamma_x^\varsigma(I_x^\varsigma)
     =\int_1^{5/4}x\varphi(xw)\ell_x(xw)\,dw>0.
\]
As $x\to\infty$,
\[
 p_x\longrightarrow p_*:=\int_1^{5/4}w e^{-w^2/2}\,dw>0.
\]
\end{enumerate}
\end{lemma}

\begin{proof}
We first prove that $\Gamma_x^\varsigma$ is a probability law and
establish the standard normal law in (ii).
Changing from $r$ to $u=r\sqrt Q$, we obtain
\[
 \Gamma_x^\varsigma(du,d\omega)
 =
 x^{-2}u e^{\varsigma z_xuT(\omega)-u^2/2}\,du\,\P(d\omega),
 \qquad u>0.
\]
Since $T$ has density $\varphi$ under $\P$, Tonelli's theorem gives,
for every nonnegative measurable $f$ on $(0,\infty)\times\R$,
\[
 \begin{aligned}
 \int f(u,T(\omega))\,\Gamma_x^\varsigma(du,d\omega)
 &=x^{-2}\int_0^\infty u e^{-u^2/2}
       \E\big[f(u,T)e^{\varsigma z_xuT}\big]\,du\\
 &=x^{-2}\int_0^\infty\int_{\R}
       f(u,t)u e^{\varsigma z_xut-u^2/2}\varphi(t)\,dt\,du.
 \end{aligned}
\]
Thus the joint density of $(u,T)$ under $\Gamma_x^\varsigma$ is
\[
 x^{-2}u e^{\varsigma z_xut-u^2/2}\varphi(t)
 =x^{-2}u e^{-u^2/(2x^2)}\varphi(t-\varsigma z_xu),\qquad u>0,
\]
where the equality follows by completing the square and using
$1-z_x^2=x^{-2}$.
With $w=u/x$ and $y=t-\varsigma z_xu$, we have $du\,dt=x\,dw\,dy$.
The joint density of $(u/x,\mathsf E)$ at $(w,y)$ is therefore
\[
 w e^{-w^2/2}\ind_{\{w>0\}}\,\varphi(y).
\]
Both factors integrate to one. Thus $\Gamma_x^\varsigma$ is a
probability law and $\mathsf E\sim\N(0,1)$, proving (ii).

Integration by parts and completion of the square give
\[
 \int_0^\infty u e^{vu-u^2/2}\,du
 =1+v\int_0^\infty e^{vu-u^2/2}\,du
 =1+v\sqrt{2\pi}e^{v^2/2}\Phi(v).
\]
Take $v=\varsigma z_xT$ to obtain the sample marginal
$\ell_x(\varsigma T)\,d\P$ and the expression
\eqref{app:radial-likelihood-explicit}.
For $v\le0$, the integral formula gives $0\le\ell_x(v)\le x^{-2}$.
For $v>0$, we have $z_xv\le v$,
$e^{z_x^2v^2/2}\le e^{v^2/2}$, and $\Phi(z_xv)\le1$.
Thus \eqref{app:radial-likelihood-explicit} gives
\[
 \begin{aligned}
 \ell_x(v)
 &\le x^{-2}\{1+\sqrt{2\pi}\,v e^{v^2/2}\}\le\sqrt{2\pi}\,x^{-2}(1+v)e^{v^2/2}.
 \end{aligned}
\]
 Together with the case $v\le0$, this completes (i).

It remains to prove (iii).
Since $\varsigma T$ has the same $\P$-law for both signs, the band
masses agree. The positivity of $\ell_x$ and $\varphi$ gives $p_x>0$.
Multiplying \eqref{app:radial-likelihood-explicit}, with $v=xw$, by
$x\varphi(xw)$ and using $1-z_x^2=x^{-2}$ gives
\begin{equation}\label{app:radial-local-limit}
 x\varphi(xw)\ell_x(xw)
 =x^{-1}\varphi(xw)+z_xw e^{-w^2/2}\Phi(z_xxw)
 \longrightarrow w e^{-w^2/2}
\end{equation}
as $x\to\infty$, uniformly for $w$ in compact subsets of $(0,\infty)$.
Indeed, on every $[a,b]\subset(0,\infty)$,
$z_x\to1$, $\Phi(z_xxw)\to1$, and $x^{-1}\varphi(xw)\to0$ uniformly.
Integrating over $[1,5/4]$ therefore gives
\[
 p_x=\int_1^{5/4}x\varphi(xw)\ell_x(xw)\,dw
 \longrightarrow
 p_*:=\int_1^{5/4}w e^{-w^2/2}\,dw>0.
\]
This proves (iii).
\end{proof}

\subsection{Proof of Lemma~\ref{main:noncollapse}}
\label{supp:noncollapse-proof}

For Borel $D\subset(0,\infty)$, recall the conditional law
\eqref{def:band-law}:
\[
 H_x^\varsigma(D)=\Gamma_x^\varsigma(Q\in D\mid I_x^\varsigma).
\]
Choose a median $m_x$ of $H_x^+$ as in \eqref{def:band-median}, so that
\[
 H_x^+((0,m_x])\ge\tfrac12,\qquad
 H_x^+([m_x,\infty))\ge\tfrac12.
\]
Constants below are independent of $r$, the tilted mean of $Q$,
and the sign, unless stated otherwise.

Let $\{Y_{n,i}\}$ be the rowwise independent array from
Proposition~\ref{array:joint-representation}. Its sums satisfy
\[
 A_n=\sum_{i=1}^nY_{n,i},\qquad
 \widehat V_n^2=\sum_{i=1}^nY_{n,i}^2,\qquad
 (A_n,\widehat V_n^2)\dto(A,Q).
\]

For $s\in\R$ and $t>0$, define
\begin{equation}\label{def:product-tilt}
 M(s,t)=\E e^{sA-tQ},\qquad
 dP_{s,t}=M(s,t)^{-1}e^{sA-tQ}\,d\P,
\end{equation}
and write $\E_{s,t}$ for expectation under $P_{s,t}$.

We first record convergence and a variance bound for these tilted laws.

\begin{lemma}
\label{arr:transform-ui}
For every $s\in\R$ and $t>0$,
\[
 M_n(s,t):=\E e^{sA_n-t\widehat V_n^2}
 \longrightarrow
 M(s,t)\in(0,\infty).
\]
The convergence remains valid after multiplying the integrands by any
fixed polynomial and any bounded continuous function of the two coordinates.
For the tilted law in \eqref{def:product-tilt},
\begin{equation}\label{arr:variance-interface}
 0<\E_{s,t}Q<\infty,\qquad \operatorname{Var}_{s,t}A\le\E_{s,t}Q.
\end{equation}
\end{lemma}

The proof is given in Section~\ref{supp:transform-ui-proof}.

We use this variance bound to compare the square sum with its tilted mean.

\begin{lemma}\label{core:good-state}
Use $\Gamma_x^\varsigma$ from \eqref{main:radial-definition}, and let
$\Omega$ be the sample space of $\P$.
For $r>0$, set $s=\varsigma z_xr$ and $t=r^2/2$, and define
\[
 a:=\frac{\varsigma\E_{s,t}A}{x\sqrt{\E_{s,t}Q}},\qquad
 p:=\frac{r\sqrt{\E_{s,t}Q}}x.
\]
There exist constants $K_0,K_1<\infty$ and $x_0>2$ such that, for every
$x\ge x_0$ and either sign $\varsigma$, the set
\begin{equation}\label{core:admissible-radii}
 C_x^\varsigma
 :=\left\{r>0:\begin{aligned}
 &a,p\in[1-K_0x^{-1/4},5/4+K_0x^{-1/4}],|a-p|\le K_0x^{-1/4}
 \end{aligned}\right\}
\end{equation}
is measurable, and there is a measurable set
$G_x^\varsigma\subset I_x^\varsigma\cap(C_x^\varsigma\times\Omega)$ with
\[
 \Gamma_x^\varsigma(I_x^\varsigma\setminus G_x^\varsigma)
 \le K_1x^{-3/2},
\]
such that, for every $(r,\omega)\in G_x^\varsigma$,
\begin{equation}\label{core:good-properties}
 \left|\frac{Q(\omega)}{\E_{s,t}Q}-1\right|\le K_0x^{-1/4}.
\end{equation}
\end{lemma}

\begin{proof}
\proofanchor\label{supp:good-state-proof}
Fix $x>2$ and a sign $\varsigma$.
Define $C_x^\varsigma$ by \eqref{core:admissible-radii}, with $K_0$
to be chosen below.
By Lemma~\ref{arr:transform-ui}, the tilted moments are finite and
$0<\E_{s,t}Q<\infty$. These moments are measurable functions of $r$,
so $a,p$ and $C_x^\varsigma$ are measurable. Put
\begin{equation}\label{notation:radial-errors}
 d_r=\frac{{\E_{s,t}A}-\varsigma z_xrQ}{\sqrt Q},\qquad
 e_r=\frac{{\E_{s,t}A}-\varsigma z_xr{\E_{s,t}Q}}{\sqrt {\E_{s,t}Q}}.
\end{equation}
Integrating \eqref{main:radial-definition} over the sample identifies the
marginal law of $r$ under $\Gamma_x^\varsigma$:
\begin{equation}\label{core:rho}
 \begin{aligned}
 \rho_x^\varsigma(dr)
 &=x^{-2}r\left(\int Qe^{sA-tQ}\,d\P\right)dr
 =x^{-2}rM(s,t)\E_{s,t}Q\,dr.
 \end{aligned}
\end{equation}
Using the conditional-density formula and the tilted density in \eqref{def:product-tilt}, we obtain
\[
 \begin{aligned}
 \Gamma_x^\varsigma(d\omega\mid r)
 =\frac{x^{-2}rQ(\omega)e^{sA(\omega)-tQ(\omega)}}
         {x^{-2}rM(s,t)\E_{s,t}Q}\,\P(d\omega)
 =\frac{Q(\omega)}{\E_{s,t}Q}
 \frac{e^{sA(\omega)-tQ(\omega)}}{M(s,t)}\,\P(d\omega)
 =\frac{Q(\omega)}{\E_{s,t}Q}\,P_{s,t}(d\omega).
 \end{aligned}
\]
The tower property of conditional expectation and \eqref{arr:variance-interface}
yield
\begin{align*}
 \E_{\Gamma_x^\varsigma}\frac{(A-\E_{s,t}A)^2}{Q}
 &=\int_0^\infty\left[\int
     \frac{(A(\omega)-\E_{s,t}A)^2}{Q(\omega)}
     \frac{Q(\omega)}{\E_{s,t}Q}\,P_{s,t}(d\omega)
     \right]\rho_x^\varsigma(dr)\\
 &=\int_0^\infty\frac{1}{\E_{s,t}Q}
     \left[\int(A(\omega)-\E_{s,t}A)^2\,P_{s,t}(d\omega)\right]
     \rho_x^\varsigma(dr)\\
 &=\int_0^\infty\frac{\operatorname{Var}_{s,t}A}{\E_{s,t}Q}\,\rho_x^\varsigma(dr)
 \le\int_0^\infty\rho_x^\varsigma(dr)=1.
\end{align*}
By Lemma~\ref{app:radial-identities}(ii), $\mathsf E=T-\varsigma z_xr\sqrt Q$ has the standard
normal law under $\Gamma_x^\varsigma$.
Combining this fact with the preceding estimate and the
Cauchy--Schwarz inequality, we find
\begin{equation}\label{core:d-error-moment}
 \begin{aligned}
 \E_{\Gamma_x^\varsigma}d_r^2
 &=\E_{\Gamma_x^\varsigma}
     \left(\frac{\E_{s,t}A-A}{\sqrt Q}+\mathsf E\right)^2\le2\E_{\Gamma_x^\varsigma}\frac{(A-\E_{s,t}A)^2}{Q}
     +2\E_{\Gamma_x^\varsigma}\mathsf E^2\le4.
 \end{aligned}
\end{equation}

For $e_r$, Jensen's inequality under each fixed tilt, followed by
Tonelli's theorem, gives
\begin{equation}\label{core:e-error-moment}
 \begin{aligned}
 \E_{\Gamma_x^\varsigma}e_r^2
 =\int e_r^2\,d\rho_x^\varsigma
 &=x^{-2}\int_0^\infty rM(\varsigma z_xr,r^2/2)
       ({\E_{s,t}A}-\varsigma z_xr{\E_{s,t}Q})^2\,dr\\
 &\le x^{-2}\int_0^\infty rM(\varsigma z_xr,r^2/2)
       \E_{s,t}(A-\varsigma z_xrQ)^2\,dr\\
 &=\E_{\Gamma_x^\varsigma}
       \frac{(A-\varsigma z_xrQ)^2}{Q}
 =\E_{\Gamma_x^\varsigma}{\mathsf E}^2=1.
 \end{aligned}
\end{equation}

Taking absolute values in the identity
\[
 \varsigma z_xr(\E_{s,t}Q-Q)
 =\sqrt Q\,d_r-\sqrt{\E_{s,t}Q}\,e_r
\]
and factoring $\E_{s,t}Q-Q$, we obtain
\begin{equation}\label{core:radial-root-difference}
 r\left|\sqrt{\E_{s,t}Q}-\sqrt Q\right|
 \le\frac{\sqrt Q\,|d_r|+\sqrt{\E_{s,t}Q}\,|e_r|}
         {z_x(\sqrt Q+\sqrt{\E_{s,t}Q})}
 \le\frac{|d_r|+|e_r|}{z_x}.
\end{equation}
Suppose that $(r,\omega)\in I_x^\varsigma$ satisfies
\begin{equation}\label{core:small-errors}
 \max(|\mathsf E|,|d_r|,|e_r|)\le x^{3/4}.
\end{equation}
Since $z_x\ge\sqrt3/2$ for $x>2$, \eqref{app:radial-residual} and
\eqref{core:small-errors} give
\begin{equation}\label{core:radial-sample-approx}
 r\sqrt Q=\frac{\varsigma T-\varsigma {\mathsf E}}{z_x}
   =\frac{\varsigma T}{z_x}+O(x^{3/4}).
\end{equation}
The corresponding approximation for the tilted mean follows from
\eqref{core:radial-root-difference} and \eqref{core:small-errors}:
\begin{equation}\label{core:radial-mean-approx}
 r\sqrt{\E_{s,t}Q}=r\sqrt Q+O(x^{3/4}).
\end{equation}
For this pair, $x\le\varsigma T\le5x/4$ by the definition of
$I_x^\varsigma$, while $z_x=\sqrt{1-x^{-2}}=1+O(x^{-2})$.
Together with \eqref{core:radial-sample-approx} and
\eqref{core:radial-mean-approx}, these bounds yield
\begin{equation}\label{core:q-p-estimates}
 \begin{aligned}
 \left|\frac Q{\E_{s,t}Q}-1\right|
 &=\frac{|\sqrt Q-\sqrt{\E_{s,t}Q}|
          (\sqrt Q+\sqrt{\E_{s,t}Q})}{\E_{s,t}Q}
 =O(x^{-1/4}),\\
 p&=\frac{r\sqrt{\E_{s,t}Q}}{x}
 \in[1-O(x^{-1/4}),5/4+O(x^{-1/4})].
 \end{aligned}
\end{equation}
The definition of $e_r$ relates $a$ to $p$:
\[
 \frac{\varsigma {\E_{s,t}A}}{\sqrt {\E_{s,t}Q}}
 =z_xr\sqrt{\E_{s,t}Q}+\varsigma e_r,
 \qquad a=z_xp+\frac{\varsigma e_r}{x}.
\]
In particular,
\begin{equation}\label{core:a-estimates}
 \begin{gathered}
 |a-p|\le|1-z_x||p|+|e_r|/x=O(x^{-1/4}),\text{ and }
 a\in[1-O(x^{-1/4}),5/4+O(x^{-1/4})].
 \end{gathered}
\end{equation}
Choose $K_0$ large enough to cover \eqref{core:q-p-estimates} and
\eqref{core:a-estimates}. For all sufficiently large $x$, these estimates
and \eqref{core:admissible-radii} show that \eqref{core:good-properties}
holds for every
\[
 (r,\omega)\in I_x^\varsigma\cap
 \{\max(|\mathsf E|,|d_r|,|e_r|)\le x^{3/4}\}
 \subseteq C_x^\varsigma\times\Omega.
\]
Set
\begin{equation}\label{core:good-set}
 G_x^\varsigma
 =I_x^\varsigma\cap(C_x^\varsigma\times\Omega)\cap
   \{\max(|\mathsf E|,|d_r|,|e_r|)\le x^{3/4}\}.
\end{equation}
We then have
\[
 \begin{aligned}
 I_x^\varsigma\setminus G_x^\varsigma
 &\subseteq
\{\max(|\mathsf E|,|d_r|,|e_r|)>x^{3/4}\}.
 \end{aligned}
\]
Combining this inclusion with Lemma~\ref{app:radial-identities}(ii) and
\eqref{core:d-error-moment}--\eqref{core:e-error-moment}, we apply
a union bound and Markov's inequality to obtain
\[
 \begin{aligned}
 \Gamma_x^\varsigma(I_x^\varsigma\setminus G_x^\varsigma)
 &\le\Gamma_x^\varsigma\!\left(
       \max(|\mathsf E|,|d_r|,|e_r|)>x^{3/4}\right)\le\frac{\E_{\Gamma_x^\varsigma}{\mathsf E}^2+
       \E_{\Gamma_x^\varsigma}d_r^2+
       \E_{\Gamma_x^\varsigma}e_r^2}{x^{3/2}}
 \le6x^{-3/2}.
 \end{aligned}
\]
Taking $K_1=6$ and $x_0$ sufficiently large completes the proof.

\end{proof}

For probability measures $\mathsf R\ll P$, define the relative entropy by
\[
 D(\mathsf R\Vert P)
 =\int\log\!\left(\frac{d\mathsf R}{dP}\right)\,d\mathsf R
 \in[0,\infty].
\]
Let
\begin{equation}\label{arr:j-definition}
 j(d)=\int_0^1\frac{1-u}{1-u+ud}\,du
     =\frac{1-d+d\log d}{(1-d)^2}
       \quad(d>0,\ d\ne1),\qquad
 j(1)=\frac12.
\end{equation}
The integral shows that $j$ is decreasing and extends continuously
to $j(0)=1$.

\begin{lemma}
\label{main:tilted-estimates}
Under the hypotheses of Proposition~\ref{main:first-eliminations},
fix $s\in\R$ and $t>0$.
For every probability law $\mathsf R$ of the pair $(A,Q)$ with bounded
density relative to its law under $P_{s,t}$,
\begin{equation}\label{main:tilted-entropy}
 D(\mathsf R\Vert P_{s,t})
 \ge\frac{(\E_{s,t}A-\E_{\mathsf R}A)^2}{\E_{s,t}Q}
 j\!\left(\frac{\E_{\mathsf R}Q}{\E_{s,t}Q}\right).
\end{equation}
Here $j$ is defined in \eqref{arr:j-definition}.
The entropy is computed from the laws of $(A,Q)$.
\end{lemma}

The proof is given in Section~\ref{supp:tilted-estimates-proof}.

For conditional laws of the original sample, the following lemma
provides an upper bound to compare with \eqref{main:tilted-entropy}.

\begin{lemma}\label{core:conditional-entropy}
Under the hypotheses of Proposition~\ref{main:first-eliminations},
let $x>2$, $\varsigma\in\{-1,1\}$, and $r>0$, and set
$s=\varsigma z_xr$ and $t=r^2/2$.
Take an event $E_0\in\sigma(A,Q)$ with $\P(E_0)>0$ on which
$A$ and $Q$ are bounded.
Then $\mathsf R=\P(\,\cdot\mid E_0)$ has bounded density relative
to $P_{s,t}$ and satisfies
\begin{equation}\label{core:conditional-entropy-upper}
 D(\mathsf R\Vert P_{s,t})
 \le\log\frac{x}{\P(E_0)}-\varsigma z_xr\E_{\mathsf R}A
       +\frac{r^2}{2}\E_{\mathsf R}Q.
\end{equation}
\end{lemma}

\begin{proof}
By the definition of $\mathsf R$ and \eqref{def:product-tilt},
\[
 \frac{d\mathsf R}{d\P}
 =\frac{\ind_{E_0}}{\P(E_0)},\qquad
 \frac{dP_{s,t}}{d\P}=M(s,t)^{-1}e^{sA-tQ}>0.
\]
The chain rule for Radon--Nikodym derivatives therefore gives
\[
 \frac{d\mathsf R}{dP_{s,t}}
 =\frac{d\mathsf R/d\P}{dP_{s,t}/d\P}
 =\frac{M(s,t)}{\P(E_0)}\ind_{E_0}
       \exp\{-\varsigma z_xrA+r^2Q/2\}.
\]
Since $A,Q$ are bounded on $E_0$ and $E_0\in\sigma(A,Q)$, this
density is a bounded function of $(A,Q)$. We can therefore compute
its entropy from the pair laws. Taking its logarithm and expectation
under $\mathsf R$ gives
\[
 D(\mathsf R\Vert P_{s,t})
 =\log\frac1{\P(E_0)}-\varsigma z_xr\E_{\mathsf R}A
       +\frac{r^2}{2}\E_{\mathsf R}Q+\log M(s,t).
\]
Completing the square gives
\[
 \varsigma z_xrA-r^2Q/2
 =z_x^2T^2/2-(r\sqrt Q-\varsigma z_xT)^2/2
 \le z_x^2T^2/2.
\]
Since $T\sim\N(0,1)$ and $1-z_x^2=x^{-2}$,
\[
 M(s,t)\le\E e^{z_x^2T^2/2}=x.
\]
Substituting this bound into the entropy identity proves
\eqref{core:conditional-entropy-upper}.
\end{proof}

\begin{lemma}
\label{core:gap}
Fix $0<a_*<a^*<\infty$ with $a^*/a_*<\sqrt2$.  There exists $g>0$,
depending only on these two bounds, such that for every
$a,b\in[a_*,a^*]$ and every $0<d<1$,
\begin{equation}\label{core:linear-gap}
 (a+b\sqrt d)^2j(d)-\frac12(b+a\sqrt d)^2
 \ge g(1-\sqrt d).
\end{equation}
\end{lemma}

The proof is given in Section~\ref{supp:gap-proof}.

\begin{lemma}
\label{core:no-drop}
Use $H_x^\varsigma$ from \eqref{def:band-law} and take $s,t,a,p$ as in
Lemma~\ref{core:good-state}.
There exists $x_0$ such that, for every
$x\ge x_0$, every $\varsigma\in\{-1,1\}$, and every
$r\in C_x^\varsigma$ defined in \eqref{core:admissible-radii},
\begin{equation}\label{core:no-drop-form}
 H_x^{-\varsigma}\big((0,(1-x^{-1/8}){\E_{s,t}Q}]\big)<x^{-1/8}.
\end{equation}
\end{lemma}
\begin{proof}
We argue by contradiction, first for $\varsigma=1$.
Recall from Lemma~\ref{core:good-state} that $s=z_xr$, $t=r^2/2$, and
\[
 a=\frac{\E_{s,t}A}{x\sqrt{\E_{s,t}Q}},\qquad
 p=\frac{r\sqrt{\E_{s,t}Q}}x.
\]
Suppose that, for arbitrarily large $x$, there is an $r\in C_x^+$ such that
\begin{equation}\label{core:no-drop-contradiction}
 H_x^-(Q\le(1-x^{-1/8})\E_{s,t}Q)\ge x^{-1/8}.
\end{equation}
By \eqref{core:admissible-radii},
\[
 a,p\in[1-K_0x^{-1/4},5/4+K_0x^{-1/4}],\qquad
 |a-p|\le K_0x^{-1/4}.
\]

Since $j(d)\to1$ as $d\downarrow0$ by \eqref{arr:j-definition},
we can choose a fixed $d_0\in(0,1/2)$ such that
\begin{equation}\label{core:d0choice}
 j(d_0)>\frac{25}{32}+\frac{25}{16}\sqrt{d_0}+\frac{25d_0}{32}.
\end{equation}
First consider the case
$H_x^-(Q\le d_0\E_{s,t}Q)\ge x^{-1/8}/2$ and take
\[
 E_0=\{-5x/4\le T\le-x,\ Q\le d_0\E_{s,t}Q\}.
\]
On $E_0$, we have
\begin{equation}\label{core:coarse-event-bounds}
 |A|\le\tfrac54x\sqrt{d_0\E_{s,t}Q},\qquad
 Q\le d_0\E_{s,t}Q.
\end{equation}
Since $p_x=\Gamma_x^-(I_x^-)$, the conditional-probability formula
in \eqref{def:band-law} and the sample marginal in
Lemma~\ref{app:radial-identities}(i) give
\begin{equation}\label{core:coarse-band-mass}
 \begin{aligned}
 p_xH_x^-(Q\le d_0\E_{s,t}Q)
 &=\Gamma_x^-\bigl(\{Q\le d_0\E_{s,t}Q\}\cap I_x^-\bigr)
 =\E[\ell_x(-T)\ind_{E_0}].
 \end{aligned}
\end{equation}
It follows from Lemma~\ref{app:radial-identities}(i)
with $v=-T$ and $x\le-T\le5x/4$ that
\begin{equation}\label{core:coarse-density-bound}
 \begin{aligned}
 \ell_x(-T)
 &\le Cx^{-2}(1+|T|)e^{T^2/2}
 \le Cx^{-2}(1+5x/4)e^{25x^2/32}
 \le Cx e^{25x^2/32}.
 \end{aligned}
\end{equation}
By Lemma~\ref{app:radial-identities}(iii), $p_x$ is bounded below
for all sufficiently large $x$. Combining this with
$H_x^-(Q\le d_0\E_{s,t}Q)\ge x^{-1/8}/2$ and
\eqref{core:coarse-band-mass}--\eqref{core:coarse-density-bound}, we obtain
\begin{equation}\label{core:coarse-event-cost}
 \P(E_0)\ge c x^{-9/8}e^{-25x^2/32},\qquad
 \log(1/\P(E_0))\le25x^2/32+O(\log x),
\end{equation}
for some constant $c>0$.
Set $\mathsf R=\P(\,\cdot\mid E_0)$.
Since $A=T\sqrt Q\le-x\sqrt Q<0$ on $E_0$, taking expectations
under $\mathsf R$ and using \eqref{core:coarse-event-bounds} gives
\begin{equation}\label{core:coarse-conditioned-moments}
 \E_{\mathsf R}A\le-x\E_{\mathsf R}\sqrt Q<0,\qquad
 \E_{\mathsf R}Q\le d_0\E_{s,t}Q.
\end{equation}
By \eqref{core:admissible-radii} and the definition of $a$, for large $x$,
\begin{equation}\label{core:coarse-tilted-mean}
 \E_{s,t}A=a x\sqrt{\E_{s,t}Q}
 \ge(1-K_0x^{-1/4})x\sqrt{\E_{s,t}Q}>0.
\end{equation}
Lemma~\ref{core:conditional-entropy} verifies the bounded-density
hypothesis of \eqref{main:tilted-entropy}.
The estimates \eqref{core:coarse-conditioned-moments}--\eqref{core:coarse-tilted-mean},
together with the nonnegativity and monotonicity of $j$ from
\eqref{arr:j-definition}, therefore give
\begin{equation}\label{core:coarse-kl}
 \begin{aligned}
 D(\mathsf R\Vert P_{s,t})
 &\ge\frac{(\E_{s,t}A-\E_{\mathsf R}A)^2}{\E_{s,t}Q}
       j\!\left(\frac{\E_{\mathsf R}Q}{\E_{s,t}Q}\right)
 \ge\frac{(\E_{s,t}A)^2}{\E_{s,t}Q}\,j(d_0)
 =a^2x^2j(d_0).
 \end{aligned}
\end{equation}
Combining \eqref{core:coarse-kl} with \eqref{core:conditional-entropy-upper}
yields
\begin{equation}\label{core:coarse-entropy-comparison}
 a^2x^2j(d_0)
 \le\log\frac{x}{\P(E_0)}-z_xr\E_{\mathsf R}A
       +\frac{r^2}{2}\E_{\mathsf R}Q.
\end{equation}
The probability bound \eqref{core:coarse-event-cost}, together with
\eqref{core:coarse-event-bounds}, $z_x\le1$, and the definition of $p$,
controls the three terms on the right:
\[
 \begin{aligned}
 \log\frac{x}{\P(E_0)}&\le\frac{25x^2}{32}+O(\log x),\quad
 -z_xr\E_{\mathsf R}A\le\tfrac54p\sqrt{d_0}\,x^2,\text{ and }
 \frac{r^2}{2}\E_{\mathsf R}Q\le\frac{p^2d_0x^2}{2}.
 \end{aligned}
\]
Inserting these estimates into \eqref{core:coarse-entropy-comparison} and applying
the bounds on $a,p$ in \eqref{core:admissible-radii} give
\[
 (1-K_0x^{-1/4})^2j(d_0)
 \le\frac{25}{32}+(5/4+K_0x^{-1/4})\tfrac54\sqrt{d_0}
 +\frac{(5/4+K_0x^{-1/4})^2d_0}{2}+O(x^{-2}\log x)
\]
which contradicts \eqref{core:d0choice} for large $x$.
Consequently, $H_x^-(Q\le d_0\E_{s,t}Q)<x^{-1/8}/2$.
Together with \eqref{core:no-drop-contradiction}, this gives
\[
 H_x^-\!\left(d_0\E_{s,t}Q<Q\le(1-x^{-1/8})\E_{s,t}Q\right)
 \ge\frac12x^{-1/8}.
\]
The definition \eqref{def:band-law}, with $p_x=\Gamma_x^-(I_x^-)$,
therefore converts this to the following bound under $\Gamma_x^-$:
\begin{equation}\label{core:remaining-band-mass}
 \begin{aligned}
 &\Gamma_x^-\!\left(1\le-T/x\le5/4,\quad
       d_0\le Q/(\E_{s,t}Q)\le1-x^{-1/8}\right)\\
 &=p_xH_x^-\!\left(d_0\E_{s,t}Q\le Q
                       \le(1-x^{-1/8})\E_{s,t}Q\right)\\
 &\ge p_x\left[H_x^-\!\left(Q\le(1-x^{-1/8})\E_{s,t}Q\right)
                   -H_x^-\!\left(Q\le d_0\E_{s,t}Q\right)\right]
       >\frac{p_x}{2}x^{-1/8}.
 \end{aligned}
\end{equation}

Partition each coordinate interval of
$[1,5/4]\times[d_0,1-x^{-1/8}]$ into $m=\lfloor x^{1/4}\rfloor$
equal parts, assigning each shared boundary to one box.
The resulting $m^2\le x^{1/2}$ disjoint boxes have side lengths
\[
 \frac1{4m}\le x^{-1/4},\qquad
 \frac{1-x^{-1/8}-d_0}{m}\le x^{-1/4}
\]
for all sufficiently large $x$.
By \eqref{core:remaining-band-mass} and the pigeonhole principle,
one box has preimage $E\in\sigma(A,Q)$ satisfying
\begin{equation}\label{core:box-mass}
 \Gamma_x^-(E)\ge\frac{p_x x^{-1/8}}{2m^2}
                \ge\frac{p_x}{2}x^{-5/8}.
\end{equation}
Write $(t_E,q_E)$ for the lower left corner of this box. Then, on $E$,
\begin{equation}\label{core:box-coordinate-bounds}
 \begin{gathered}
 t_E\le-T/x\le t_E+x^{-1/4},\qquad
 q_E\le Q/(\E_{s,t}Q)\le q_E+x^{-1/4},\\
 t_E\in[1,5/4],\qquad q_E\in[d_0,1-x^{-1/8}].
 \end{gathered}
\end{equation}
Every point of $E$ still has $Q/(\E_{s,t}Q)\le1-x^{-1/8}$.
On $E$, the bound in Lemma~\ref{app:radial-identities}(i), applied
with $v=-T$, gives
\begin{equation}\label{core:box-likelihood}
 \ell_x(-T)\le C_3x\exp\{(t_E+x^{-1/4})^2x^2/2\}.
\end{equation}
Lemma~\ref{app:radial-identities}(i), \eqref{core:box-mass}, and
\eqref{core:box-likelihood} give
\[
 \frac{p_x}{2}x^{-5/8}\le\Gamma_x^-(E)
 =\E[\ell_x(-T)\ind_E]
 \le C_3x e^{(t_E+x^{-1/4})^2x^2/2}\P(E).
\]
Since $p_x$ is bounded below for large $x$ by
Lemma~\ref{app:radial-identities}(iii),
\[
 \P(E)\ge c_3x^{-13/8}e^{-(t_E+x^{-1/4})^2x^2/2}.
\]
Hence
\begin{equation}\label{core:box-cost}
 \log\frac1{\P(E)}
 \le\frac{t_E^2x^2}{2}+K_1x^{7/4}+K_1\log x.
\end{equation}
Put $\mathsf R=\P(\,\cdot\mid E)$.
Since $A/(x\sqrt{\E_{s,t}Q})=-(-T/x)\sqrt{Q/(\E_{s,t}Q)}$ on $E$
and the square root is Lipschitz on $[d_0,1]$,
\eqref{core:box-coordinate-bounds} gives, after taking expectations
under $\mathsf R$,
\begin{equation}\label{core:box-moments}
 \frac{\E_{\mathsf R}A}{x\sqrt {\E_{s,t}Q}}=-t_E\sqrt{q_E}+O(x^{-1/4}),\qquad
 \frac{\E_{\mathsf R}Q}{{\E_{s,t}Q}}=q_E+O(x^{-1/4}).
\end{equation}
The variables $A,Q$ are bounded on $E$ for fixed $r,x$, so
Lemmas~\ref{core:conditional-entropy} and~\ref{main:tilted-estimates}
apply to $\mathsf R$.
On $[d_0,1]$, $j$ has bounded derivatives; inserting
\eqref{core:box-moments} into \eqref{main:tilted-entropy} therefore yields
\begin{equation}\label{core:box-lower}
 D({\mathsf R}\Vert P_{s,t})\ge x^2(a+t_E\sqrt{q_E})^2j(q_E)-K_2x^{7/4}.
\end{equation}
In the upper bound \eqref{core:conditional-entropy-upper}, take $E_0=E$
and use \eqref{core:box-cost} and~\eqref{core:box-moments} to obtain
\begin{align*}
 D(\mathsf R\Vert P_{s,t})/x^2
 &\le t_E^2/2+z_xp t_E\sqrt{q_E}+p^2q_E/2
       +O(x^{-1/4})+O(x^{-2}\log x)\\
 &=\tfrac12(t_E+a\sqrt{q_E})^2+O(x^{-1/4})+O(x^{-2}\log x),
\end{align*}
where $p=a+O(x^{-1/4})$ and $z_x=1+O(x^{-2})$. Thus
\begin{equation}\label{core:box-upper}
 D({\mathsf R}\Vert P_{s,t})\le\frac{x^2}{2}(t_E+a\sqrt{q_E})^2
                    +K_3x^{7/4}+K_3\log x,
\end{equation}
for some constant $K_{3}>0$.
Choose fixed $a_*<1<5/4<a^*$ with $a^*/a_*<\sqrt2$.
For all large $x$, $a,t_E\in[a_*,a^*]$.
Since $1-\sqrt{q_E}\ge(1-q_E)/2\ge x^{-1/8}/2$, applying
Lemma~\ref{core:gap} with $(b,d)=(t_E,q_E)$ to
\eqref{core:box-lower}--\eqref{core:box-upper} yields
\[
 \frac g2x^{15/8}
 \le g(1-\sqrt{q_E})x^2
 \le K_4x^{7/4}+K_4\log x.
\]
After division by $x^{15/8}$, the right-hand side tends to zero,
contradicting $g>0$.
This proves \eqref{core:no-drop-form} for $\varsigma=1$, with constants
and thresholds uniform in $r$ and $\E_{s,t}Q$.

Replacing $A$ by $-A$ preserves the assumptions and interchanges
the positive and negative tilts and bands. Applying the result to
$(-A,Q)$ gives the case $\varsigma=-1$ with the same constants and
thresholds.
\end{proof}

\begin{lemma}\label{arr:overlap}
Suppose that, for some $\alpha>0$, both signs satisfy
\begin{equation}\label{supp:overlap-hypothesis}
 H_x^\varsigma\!\left(
 [m_x(1-4x^{-\alpha}),m_x(1+4x^{-\alpha})]\right)
 \ge1-3x^{-\alpha}
\end{equation}
for all sufficiently large $x$. There is a constant $K$ such that,
for every sufficiently large real $x$,
\[
 |\log m_{9x/8}-\log m_x|\le Kx^{-\alpha}.
\]
\end{lemma}
\begin{proof}
Define
\[
 \widehat H_x(d\omega)
 =\frac{\ind_{\{x\le T\le 5x/4\}}\ell_x(T)}{p_x}\,\P(d\omega),
\]
which is the sample marginal of the law conditional on $I_x^+$, and its
$Q$-marginal is $H_x^+$. Let
$J_x=I_x^+\cap I_{9x/8}^+=\{9x/8\le T\le 5x/4\}$
be the intersection of the two conditioning events.
The uniform asymptotic formula \eqref{app:radial-local-limit}
implies
\begin{equation}\label{core:overlap-mass}
 \widehat H_x(J_x)\longrightarrow
 \frac1{p_*}\int_{9/8}^{5/4} w e^{-w^2/2}\,dw>0.
\end{equation}
By \eqref{app:radial-likelihood-explicit}, the ratio of the two
conditional densities on $J_x$ is
\begin{equation}\label{core:overlap-density-ratio}
 \frac{\ell_{9x/8}(T)/p_{9x/8}}{\ell_x(T)/p_x}
 =\frac{64}{81}\frac{p_x}{p_{9x/8}}
 \frac{1+\sqrt{2\pi}\,z_{9x/8}T
       e^{z_{9x/8}^2T^2/2}\Phi(z_{9x/8}T)}
      {1+\sqrt{2\pi}\,z_xT
       e^{z_x^2T^2/2}\Phi(z_xT)}.
\end{equation}
Recall that $p_x\to p_*>0$ by \eqref{main:band-positive-mass}
and $z_x=\sqrt{1-x^{-2}}$. Hence, as $x\to\infty$,
\begin{equation}\label{core:overlap-estimates}
 \frac{p_x}{p_{9x/8}}\longrightarrow1,
 \qquad \frac{z_{9x/8}}{z_x}\longrightarrow1,
 \qquad z_{9x/8}^2-z_x^2=\frac{17}{81x^2}.
\end{equation}
On $J_x$, we have $T\ge9x/8$. Since $z_{9x/8}\ge z_x$
and $\Phi(v)\to1$ as $v\to\infty$,
\begin{equation}\label{core:overlap-uniform-estimates}
 z_{9x/8}T\ge z_xT\ge\frac98\sqrt{x^2-1}\longrightarrow\infty,
 \qquad \frac{\Phi(z_{9x/8}T)}{\Phi(z_xT)}\longrightarrow1,
\end{equation}
uniformly on $J_x$. Substituting
\eqref{core:overlap-estimates}--\eqref{core:overlap-uniform-estimates}
into \eqref{core:overlap-density-ratio} gives, uniformly on $J_x$,
\begin{equation}\label{core:overlap-ratio-asymptotic}
 \frac{\ell_{9x/8}(T)/p_{9x/8}}{\ell_x(T)/p_x}
 =\frac{64}{81}
   \exp\left\{\frac{17T^2}{162x^2}\right\}(1+o(1)).
\end{equation}
Since $9/8\le T/x\le5/4$ on $J_x$,
\eqref{core:overlap-mass} and \eqref{core:overlap-ratio-asymptotic}
yield constants $\gamma>0$ and $L<\infty$, independent of $x$,
such that, for all sufficiently large $x$,
\begin{equation}\label{core:overlap-comparison}
 \widehat H_x(J_x)\ge\gamma,
 \qquad
 L^{-1}\le
 \frac{\ell_{9x/8}(T)/p_{9x/8}}{\ell_x(T)/p_x}
 \le L\quad\hbox{on }J_x.
\end{equation}

Set
\[
 B_x=[m_x(1-4x^{-\alpha}),m_x(1+4x^{-\alpha})].
\]
By \eqref{core:overlap-comparison} and \eqref{supp:overlap-hypothesis},
\begin{align*}
 \widehat H_x(J_x\cap\{Q\in B_x\cap B_{9x/8}\})
 &\ge\widehat H_x(J_x)
      -\widehat H_x(Q\notin B_x)
      -\widehat H_x(J_x\cap\{Q\notin B_{9x/8}\})\\
 &\ge\gamma-3x^{-\alpha}-3L(9x/8)^{-\alpha}>0
\end{align*}
for all sufficiently large $x$. Hence the two deterministic intervals
$B_x$ and $B_{9x/8}$ intersect. Comparing their endpoints gives
\[
 \frac{1-4x^{-\alpha}}{1+4(9x/8)^{-\alpha}}
 \le\frac{m_{9x/8}}{m_x}
 \le\frac{1+4x^{-\alpha}}{1-4(9x/8)^{-\alpha}}.
\]
For large $x$, the denominators are positive. Using
$|\log(1+u)|\le2|u|$ for $|u|\le1/2$ gives
\begin{equation}\label{core:summable-scale}
 |\log m_{9x/8}-\log m_x|\le Cx^{-\alpha}
\end{equation}
for every sufficiently large real $x$.

\end{proof}

\begin{proof}[Proof of Lemma~\ref{main:noncollapse}]
Let $K_0,K_1$ be the constants in Lemma~\ref{core:good-state}.
Take $x$ sufficiently large that the conclusions of
Lemmas~\ref{core:good-state} and~\ref{core:no-drop} hold.
Recall $C_x^\varsigma$ from \eqref{core:admissible-radii} and
$G_x^\varsigma$ from \eqref{core:good-set}.
Since $p_x\to p_*>0$ by \eqref{main:band-positive-mass},
Lemma~\ref{core:good-state} gives, for all sufficiently large $x$,
\[
 \Gamma_x^\varsigma\!\left((G_x^\varsigma)^c\mid I_x^\varsigma\right)
 =\frac{\Gamma_x^\varsigma(I_x^\varsigma\setminus G_x^\varsigma)}{p_x}
 \le\frac{K_1x^{-3/2}}{p_x}\le x^{-1/8}.
\]
Fix $\varsigma\in\{-1,1\}$ and $u>0$ such that
$H_x^\varsigma(Q\ge u)\ge2x^{-1/8}$.
By \eqref{def:band-law},
\[
 \Gamma_x^\varsigma\!\left(
 G_x^\varsigma\cap\{Q\ge u\}\mid I_x^\varsigma\right)
 \ge H_x^\varsigma(Q\ge u)-x^{-1/8}
 \ge x^{-1/8}>0.
\]
Thus there is a pair $(r,\omega)\in G_x^\varsigma$ with $Q(\omega)\ge u$ and \eqref{core:good-set} ensures that $r\in C_x^\varsigma$.
Set $h=(1-x^{-1/8})/(1+K_0x^{-1/4})$.
By \eqref{core:good-properties},
\[
 hu=\frac{(1-x^{-1/8})u}{1+K_0x^{-1/4}}
 \le(1-x^{-1/8})\E_{\varsigma z_xr,r^2/2}Q.
\]
Applying \eqref{core:no-drop-form} at this $r$ therefore gives
\begin{equation}\label{core:tail-transfer}
 H_x^{-\varsigma}(Q\le hu)
 \le H_x^{-\varsigma}\!\left(
 Q\le(1-x^{-1/8})\E_{\varsigma z_xr,r^2/2}Q\right)
 <x^{-1/8}.
\end{equation}

Since $m_x$ is a median of $H_x^+$, we have
$H_x^+(Q\ge m_x)\ge1/2>2x^{-1/8}$ and
$H_x^+(Q\le m_x)\ge1/2$.
Since $h=1-x^{-1/8}+O(x^{-1/4})$, we have
$h\ge1-2x^{-1/8}$ for all sufficiently large $x$.
Applying \eqref{core:tail-transfer} with $(\varsigma,u)=(1,m_x)$ yields
\begin{equation}\label{core:minus-lower-tail}
 H_x^-(Q\le m_x(1-2x^{-1/8}))<x^{-1/8}.
\end{equation}
Suppose that $H_x^-(Q\ge m_x(1+2x^{-1/8}))\ge2x^{-1/8}$.
Since $h(1+2x^{-1/8})\ge1$, applying
\eqref{core:tail-transfer} with
$(\varsigma,u)=(-1,m_x(1+2x^{-1/8}))$ would give
$H_x^+(Q\le m_x)<x^{-1/8}$, contrary to the median inequality.
Thus
\begin{equation}\label{core:minus-upper-tail}
 H_x^-(Q\ge m_x(1+2x^{-1/8}))<2x^{-1/8}.
\end{equation}
Combining \eqref{core:minus-lower-tail} and
\eqref{core:minus-upper-tail} gives
\begin{equation}\label{core:minus-band}
 H_x^-\big([m_x(1-2x^{-1/8}),m_x(1+2x^{-1/8})]\big)
 \ge1-3x^{-1/8}.
\end{equation}
Using \eqref{core:minus-lower-tail} and \eqref{core:minus-upper-tail},
we apply \eqref{core:tail-transfer} and its contrapositive, respectively,
to obtain
\begin{equation}\label{core:plus-tails}
 H_x^+(Q\le hm_x(1-2x^{-1/8}))<x^{-1/8},
 \qquad
 H_x^+(Q\ge m_x(1+2x^{-1/8})/h)<2x^{-1/8}.
\end{equation}
For sufficiently large $x$, both thresholds lie in
$[m_x(1-4x^{-1/8}),m_x(1+4x^{-1/8})]$.
Combining \eqref{core:plus-tails} with \eqref{core:minus-band} gives,
for both signs and all sufficiently large real $x$,
\begin{equation}\label{main:common-energy}
 H_x^\varsigma\!\left(
 [m_x(1-4x^{-1/8}),m_x(1+4x^{-1/8})]\right)
 \ge1-3x^{-1/8}.
\end{equation}
Lemma~\ref{arr:overlap}, with $\alpha=1/8$, now gives a constant
$K<\infty$, independent of $x$, such that, for all sufficiently large real $x$,
\begin{equation}\label{main:overlap-energy}
 |\log m_{9x/8}-\log m_x|\le Kx^{-1/8}.
\end{equation}
Fix $x_0$ large enough that \eqref{main:overlap-energy} holds
for every $x\ge x_0$, and put $x_j=(9/8)^jx_0$.
For $l>j$, \eqref{main:overlap-energy} gives
\[
 |\log m_{x_l}-\log m_{x_j}|
 \le K\sum_{r=j}^{l-1}x_r^{-1/8}
 \le\frac{Kx_j^{-1/8}}{1-(8/9)^{1/8}}\longrightarrow0
 \qquad(j\to\infty).
\]
Thus $(\log m_{x_j})$ is Cauchy in $\R$.
Its finite limit gives a positive finite limit for $m_{x_j}$ by
exponentiation.
\end{proof}

\section{Proof of Lemma~\ref{lem:gaussian-energy-improvement}}
\label{app:positive-gaussian-energy}

We use \eqref{array:positive-gaussian-setup}, with $\sigma>0$ and $\kappa=0$:
\[
 A=\sigma Z+U,\qquad Q=\sigma^2+R,\qquad
 T=\frac{A}{\sqrt Q}\sim\N(0,1).
\]
Here $Z\sim\N(0,1)$ is independent of $(U,R)$, and $R\ge0$.
Write $\varphi$ for the standard normal density.
The pair $(U,R)$ is defined in \eqref{def:residual}.
As above, $\P$ and $\E$ denote the original law and its expectation.
For $x>2$ and $\varsigma\in\{-1,1\}$, recall $z_x=\sqrt{1-x^{-2}}$ and
\[
 \Gamma_x^\varsigma(dr,d\omega)
 =x^{-2}rQe^{\varsigma z_xrA-r^2Q/2}\,dr\,\P(d\omega),
 \qquad r>0.
\]
Recall the bands and conditional laws from
\eqref{main:tail-bands} and \eqref{def:band-law}:
\[
 \begin{gathered}
 I_x^\varsigma=\{x\le\varsigma T\le5x/4\},\qquad
 p_x=\Gamma_x^\varsigma(I_x^\varsigma),\\
 H_x^\varsigma(D)=\Gamma_x^\varsigma(Q\in D\mid I_x^\varsigma),
 \qquad D\subset(0,\infty)\text{ Borel}.
 \end{gathered}
\]
As in \eqref{def:band-median}, $m_x$ is a median of $H_x^+$:
$H_x^+((0,m_x])\ge1/2$ and $H_x^+([m_x,\infty))\ge1/2$.
We write $v_x^\varsigma=\E_{\Gamma_x^\varsigma}(R/Q)$, as in
\eqref{def:residual-noise}.

\begin{lemma}
\label{lem:tail-energy-initial}
There are constants $C<\infty$ and $x_0>2$ such that, for all
$x\ge x_0$,
\begin{equation}\label{array:initial-median-rate}
 |m_x-\sigma^2|\le Cx^{-1/8},
\end{equation}
and, for both signs $\varsigma$,
\begin{equation}\label{array:initial-band-rate}
 \E_{\Gamma_x^\varsigma}\!\left[\frac RQ\,\middle|\,I_x^\varsigma\right]
 \le Cx^{-1/8}.
\end{equation}
\end{lemma}

\begin{proof}
For every $\eta>0$, we first show that some $\lambda>1/2$ satisfies
\begin{equation}\label{array:local-supergaussian}
 \E\!\left[e^{\lambda T^2}\ind_{\{R\ge\eta\}}\right]<\infty.
\end{equation}
For every $\lambda>0$, conditional Gaussian integration gives,
on $\{R>(2\lambda-1)\sigma^2\}$,
\begin{equation}\label{array:conditional-gaussian-integral}
 \E\!\left[e^{\lambda T^2}\mid U,R\right]
 =\sqrt{\frac{\sigma^2+R}{R-(2\lambda-1)\sigma^2}}
   \exp\left\{\frac{\lambda U^2}
                    {R-(2\lambda-1)\sigma^2}\right\}.
\end{equation}
Now fix $\eta>0$ and choose
\begin{equation}\label{array:local-lambda}
 \lambda=\frac12+\frac{\eta}{8\sigma^2}>\frac12,
 \qquad (2\lambda-1)\sigma^2=\frac\eta4<\frac\eta2.
\end{equation}
By \eqref{array:local-lambda}, on $\{R\ge\eta\}$ we have
\begin{equation}\label{array:gaussian-prefactor-bound}
 \sqrt{\frac{\sigma^2+R}{R-(2\lambda-1)\sigma^2}}
 \le\sqrt{2\left(1+\frac{\sigma^2}{R}\right)}
 \le\sqrt{2\left(1+\frac{\sigma^2}{\eta}\right)}
\end{equation}
and
\begin{equation}\label{array:gaussian-denominator-bound}
 R-(2\lambda-1)\sigma^2
 \ge R-\eta/2\ge\frac{\eta+R}{4}.
\end{equation}
Substituting \eqref{array:gaussian-prefactor-bound} and
\eqref{array:gaussian-denominator-bound} into
\eqref{array:conditional-gaussian-integral}, and using
\eqref{rigidity:residual-offset} with $d=\eta$, gives
\[
 \E[e^{\lambda T^2}\ind_{\{R\ge\eta\}}]
 \le\sqrt{2(1+\sigma^2/\eta)}\,
      \E e^{4\lambda U^2/(\eta+R)}<\infty.
\]
This proves \eqref{array:local-supergaussian} for the fixed $\eta$.

The marginal law of $\Gamma_x^\varsigma$ on the original sample space is
$\ell_x(\varsigma T)\,d\P$, where
\[
 \ell_x(v)=x^{-2}\int_0^\infty u e^{z_xvu-u^2/2}\,du.
\]
Lemma~\ref{main:radial-law} gives
\[
 \ell_x(\varsigma T)\le Cx^{-2}(1+|T|)e^{T^2/2}.
\]
For $\lambda$ in \eqref{array:local-lambda},
$(1+|T|)e^{T^2/2}\le C_\lambda e^{\lambda T^2}$.
Thus \eqref{array:local-supergaussian} gives a constant
$C_\eta<\infty$, independent of $x$ and the sign, such that
\begin{equation}\label{array:residual-vanishes}
  \Gamma_x^\varsigma(R\ge\eta)
 =\E\!\left[\ell_x(\varsigma T)\ind_{\{R\ge\eta\}}\right]
 \le C_\eta x^{-2}\longrightarrow0,
 \qquad(\eta>0).
\end{equation}
Since $p_x\to p_*>0$, for every fixed $\eta>0$,
\[
 H_x^\varsigma(Q\ge\sigma^2+\eta)
 \le\frac{C_\eta x^{-2}}{p_x}\longrightarrow0.
\]
Since $Q\ge\sigma^2$, the median inequalities in
\eqref{def:band-median} now give $m_x\to\sigma^2$.

Starting from any sufficiently large $x$, sum
\eqref{main:overlap-energy} along $x,9x/8,(9/8)^2x,\ldots$.
Since $m_{(9/8)^jx}\to\sigma^2>0$,
\begin{equation}\label{array:initial-log-median-rate}
 |\log m_x-\log\sigma^2|
 \le C\sum_{j\ge0}((9/8)^jx)^{-1/8}\le C'x^{-1/8}.
\end{equation}
By the mean value theorem, $|e^t-1|\le e|t|$ whenever $|t|\le1$.
Thus \eqref{array:initial-log-median-rate} gives, for all sufficiently
large $x$,
\[
 \begin{aligned}
 |m_x-\sigma^2|
 &=\sigma^2\left|e^{\log m_x-\log\sigma^2}-1\right|
 \le e\sigma^2|\log m_x-\log\sigma^2|
 \le Cx^{-1/8}.
 \end{aligned}
\]
This proves \eqref{array:initial-median-rate}, with a constant
independent of $x$.

Define the event
\[
 E_x=\{m_x(1-4x^{-1/8})\le Q\le m_x(1+4x^{-1/8})\}.
\]
By the definition of $H_x^\varsigma$ in \eqref{def:band-law},
\eqref{main:common-energy} gives, for both signs,
\begin{equation}\label{array:initial-event-probability}
 \begin{aligned}
 \Gamma_x^\varsigma(E_x^c\mid I_x^\varsigma)
 &=1-H_x^\varsigma\!\left(
 [m_x(1-4x^{-1/8}),m_x(1+4x^{-1/8})]\right)
 \le3x^{-1/8}.
 \end{aligned}
\end{equation}
Since $m_x\to\sigma^2$, \eqref{array:initial-median-rate} implies
that, on $E_x$,
\begin{equation}\label{array:initial-event-residual}
 0\le R=Q-\sigma^2
 \le |m_x-\sigma^2|+4m_xx^{-1/8}
 \le Cx^{-1/8}.
\end{equation}
Using $Q\ge\sigma^2$, $0\le R/Q\le1$, and
\eqref{array:initial-event-probability}--\eqref{array:initial-event-residual},
we obtain
\[
 \begin{aligned}
 \E_{\Gamma_x^\varsigma}\!\left[\frac RQ\,\middle|\,I_x^\varsigma\right]
 &=\E_{\Gamma_x^\varsigma}\!\left[
   \frac RQ\ind_{E_x}\,\middle|\,I_x^\varsigma\right]
  +\E_{\Gamma_x^\varsigma}\!\left[
   \frac RQ\ind_{E_x^c}\,\middle|\,I_x^\varsigma\right]\\
 &\le\frac{C}{\sigma^2}x^{-1/8}
      \Gamma_x^\varsigma(E_x\mid I_x^\varsigma)
      +\Gamma_x^\varsigma(E_x^c\mid I_x^\varsigma)\\
 &\le\left(\frac{C}{\sigma^2}+3\right)x^{-1/8}.
 \end{aligned}
\]
This proves \eqref{array:initial-band-rate}.
\end{proof}

\begin{lemma}
\label{lem:band-to-global}
Use $v_x^\varsigma$ from \eqref{def:residual-noise}.
Suppose $0<\beta<2$ and there are constants $C<\infty$ and $y_0>2$
such that, for every $y\ge y_0$ and both signs,
\begin{equation}\label{array:annular-energy-rate}
 \E_{\Gamma_y^\varsigma}\!\left[\frac RQ\,\middle|\,I_y^\varsigma\right]
 \le Cy^{-\beta}.
\end{equation}
Then there is $C'<\infty$ such that, for all sufficiently large $x$
and both signs,
\begin{equation}\label{array:global-energy-rate}
 v_x^\varsigma\le C'x^{-\beta}.
\end{equation}
\end{lemma}

\begin{proof}
Fix a sign, put $q=5/4$, and write $y_j=q^jy_0$.
The bands $[y_j,qy_j)$ partition $\{\varsigma T\ge y_0\}$ up to
null endpoints. Formula~\eqref{app:radial-likelihood-explicit} gives,
for $t\in[y,5y/4]$ and uniformly over $x,y\ge y_0$,
\begin{equation}\label{array:likelihood-transfer}
 \frac{\ell_x(t)}{\ell_y(t)}
 \le C(y/x)^2
       \exp\{(y^{-2}-x^{-2})t^2/2\}
 \le C'(y/x)^2e^{-y^2/(2x^2)}.
\end{equation}
Indeed, $z_x,z_y$ are bounded away from zero. For $t\ge y_0$, the
positive exponential term in \eqref{app:radial-likelihood-explicit}
uniformly dominates the additive constant $1$. The ratios $z_x/z_y$
and $\Phi(z_xt)/\Phi(z_yt)$ are uniformly bounded above and away
from zero. Here $\Phi$ is the standard normal distribution function.
For $y\le t\le5y/4$,
$(y^{-2}-x^{-2})t^2/2\le25/32-y^2/(2x^2)$.
Comparing the exponential factors gives \eqref{array:likelihood-transfer}.

By \eqref{array:likelihood-transfer} and
\eqref{array:annular-energy-rate},
\begin{align*}
 \E_{\Gamma_x^\varsigma}\!\left[
 \frac RQ\ind_{\{y\le\varsigma T<5y/4\}}\right]
 &\le C(y/x)^2e^{-y^2/(2x^2)}
       \E_{\Gamma_y^\varsigma}\!\left[\frac RQ\ind_{I_y^\varsigma}\right]\\
 &\le C'(y/x)^2e^{-y^2/(2x^2)}y^{-\beta},
\end{align*}
where $p_y\le1$ was used in the last step. Summing over $y=(5/4)^jy_0$
gives
\[
 Cx^{-\beta}\sum_{j\ge0}
       (y_j/x)^{2-\beta}e^{-(y_j/x)^2/2}=O(x^{-\beta}).
\]
Let $J=\max\{j:y_j\le x\}$. Since $y_J/x\le1<y_{J+1}/x\le q$,
\[
 \begin{aligned}
 \sum_{j\le J}(y_j/x)^{2-\beta}
 &\le\sum_{k\ge0}q^{-k(2-\beta)}
   =\frac1{1-q^{-(2-\beta)}},\\
 \sum_{j>J}(y_j/x)^{2-\beta}e^{-(y_j/x)^2/2}
 &\le\sum_{k\ge1}q^{k(2-\beta)}e^{-q^{2(k-1)}/2}<\infty.
 \end{aligned}
\]
Both bounds are independent of $x$.

On $0\le t<y_0$, $\ell_x(t)\le C_{y_0}x^{-2}$.
For $t<0$,
\[
 \ell_x(t)=x^{-2}\int_0^\infty u e^{z_xut-u^2/2}\,du
 \le x^{-2}\int_0^\infty u e^{-u^2/2}\,du=x^{-2}.
\]
Consequently,
\[
 \E_{\Gamma_x^\varsigma}\!\left[
       \frac RQ\ind_{\{\varsigma T<y_0\}}\right]
 \le C_{y_0}x^{-2}=O(x^{-\beta}).
\]
Adding the band contributions proves \eqref{array:global-energy-rate}.
The same constants cover both signs.
\end{proof}

For $s\in\R$ and $t>0$, recall the tilted law in \eqref{def:product-tilt}:
\[
 M(s,t)=\E e^{sA-tQ},\qquad
 dP_{s,t}=M(s,t)^{-1}e^{sA-tQ}\,d\P.
\]
Here $\E_{s,t}$ denotes expectation under $P_{s,t}$.

\label{supp:canonical-tilt-proof}

\begin{lemma}
\label{arr:canonical-tilt}
Use the representation
$A=\sigma Z+U$, $Q=\sigma^2+\kappa+R$ from
Proposition~\ref{array:joint-representation}. For every $s\in\R$ and $t>0$,
under $P_{s,t}$ the variable $Z$ has law $\N(\sigma s,1)$ and is
independent of $(U,R)$. Moreover, $U$ is square integrable and
\begin{equation}\label{arr:residual-tilted-variance}
 \operatorname{Var}_{s,t}U
 \le\E_{s,t}Q-\sigma^2
 =\kappa+\E_{s,t}R.
\end{equation}
\end{lemma}

\begin{proof}
The quadratic weight factors as
\[
 e^{sA-tQ}
 =e^{-t(\sigma^2+\kappa)}e^{s\sigma Z}e^{sU-tR}.
\]
Independence under $\P$ and Tonelli's theorem give
\[
 M(s,t)
 =e^{-t(\sigma^2+\kappa)+s^2\sigma^2/2}\,
     \E e^{sU-tR}\in(0,\infty),
\]
where finiteness follows from Lemma~\ref{arr:transform-ui}.
The normalized density is therefore
\[
 \frac{e^{sA-tQ}}{M(s,t)}
 =e^{s\sigma Z-s^2\sigma^2/2}
    \frac{e^{sU-tR}}{\E e^{sU-tR}}.
\]
The factors depend separately on $Z$ and $(U,R)$.
Thus independence is preserved and $Z\sim\N(\sigma s,1)$ under $P_{s,t}$.
Lemma~\ref{arr:transform-ui} gives $A\in L^2(P_{s,t})$, so
$U=A-\sigma Z\in L^2(P_{s,t})$.
By \eqref{arr:variance-interface},
\[
 \sigma^2+\operatorname{Var}_{s,t}U
 =\operatorname{Var}_{s,t}A\le\E_{s,t}Q,
\]
which proves \eqref{arr:residual-tilted-variance}.
\end{proof}

For fixed $x$, $\varsigma$, and $r$, set
$s=\varsigma z_xr$ and $t=r^2/2$. Recall the marginal in \eqref{core:rho}
and the quantities $a,p$ defined in Lemma~\ref{core:good-state}:
\[
 \begin{gathered}
 \rho_x^\varsigma(dr)
 =x^{-2}rM(\varsigma z_xr,r^2/2)\E_{s,t}Q\,dr,\\
 a=\frac{\varsigma\E_{s,t}A}{x\sqrt{\E_{s,t}Q}},\qquad
 p=\frac{r\sqrt{\E_{s,t}Q}}x.
 \end{gathered}
\]
For $\mathsf R\ll P$, write
$D(\mathsf R\Vert P)=\int\log(d\mathsf R/dP)\,d\mathsf R$
for relative entropy. The function in \eqref{arr:j-definition} is
\[
 j(d)=\int_0^1\frac{1-u}{1-u+ud}\,du,\qquad d>0.
\]

\begin{lemma}
\label{arr:parameter-reuse}
Use $\Gamma_x^\varsigma$ from \eqref{main:radial-definition} and its marginal
\eqref{core:rho}, together with the entropy bounds
\eqref{core:coarse-kl} and~\eqref{core:conditional-entropy-upper}.
Replace the error tolerance $x^{-1/4}$ in
\eqref{core:admissible-radii} and \eqref{core:good-properties} by $\eta_x$,
with one fixed constant $K_0$, and write $C_x^\varsigma$ for the corresponding
set of radii. Replace $x^{-1/8}$ in the proof of Lemma~\ref{core:no-drop}
and in the argument from \eqref{core:tail-transfer} to
\eqref{main:common-energy} in the proof of Lemma~\ref{main:noncollapse}
by $\delta_x$.
Let these deterministic tolerances satisfy
$\eta_x,\delta_x\to0$ and
\[
 \eta_x=o(\delta_x),\qquad \log x=o(\delta_x x^2),\qquad
 \log(1/(\delta_x \eta_x^2))=O(\log x).
\]
Suppose also that, for each sign, the pairs $(r,\omega)\in I_x^\varsigma$
with $r\notin C_x^\varsigma$ or failing \eqref{core:good-properties}
with the new tolerance
have total $\Gamma_x^\varsigma$-mass
$o(\delta_x)$. Then Lemma~\ref{core:no-drop} and the concentration
estimate \eqref{main:common-energy} remain valid with this choice of
$\eta_x,\delta_x$. If $\delta_x=x^{-\alpha}$ for some $\alpha>0$, then
Lemma~\ref{arr:overlap} applies as well.
\end{lemma}
\begin{proof}
Write $\eta=\eta_x>0$, $\delta=\delta_x>0$, and
\[
 \Lambda_x=\log x+\log\frac1{\delta\eta^2},\qquad
 \eta x^2+\Lambda_x=o(\delta x^2).
\]
Use $d_0$ from \eqref{core:d0choice} and $a,p$ as defined in
Lemma~\ref{core:good-state}.
Fix $r\in C_x^\varsigma$ and
reduce to $\varsigma=1$ by reflection. Put $P=P_{z_xr,r^2/2}$.
Suppose that
$H_x^-(Q\le(1-\delta)\E_PQ)\ge\delta$.

If $H_x^-(Q\le d_0\E_PQ)\ge\delta/2$, the event
$E_0=\{-5x/4\le T\le-x,\ Q\le d_0\E_PQ\}$ satisfies
\[
 \P(E_0)\ge c\delta x^{-1}e^{-25x^2/32},\qquad
 \log\frac1{\P(E_0)}\le\frac{25x^2}{32}+\log\frac{Cx}{\delta}.
\]
For $\mathsf R=\P(\,\cdot\mid E_0)$, equations
\eqref{core:coarse-kl} and~\eqref{core:conditional-entropy-upper} give
\[
 \begin{aligned}
 0<&\ j(d_0)-\frac{25}{32}-\frac{25}{16}\sqrt{d_0}
                   -\frac{25d_0}{32}\\
 \le&\ C\eta+Cx^{-2}\log\frac{Cx}{\delta}\longrightarrow0,
 \end{aligned}
\]
a contradiction.

The rectangle $1\le-T/x\le5/4$,
$d_0\le Q/(\E_PQ)\le1-\delta$ therefore has conditional mass at least
$\delta/2$ under $\Gamma_x^-(\,\cdot\mid I_x^-)$. Partition it into
at most $C\eta^{-2}$ boxes of side length at most $\eta$.
Some box $E$, with lower endpoints $t_E,q_E$, satisfies
\[
 \Gamma_x^-(E\mid I_x^-)\ge c\delta\eta^2.
\]
On $E$, $-T/x\le t_E+\eta$. The sample marginal gives
\[
 p_x\Gamma_x^-(E\mid I_x^-)
 =\E[\ell_x(-T)\ind_E]
 \le Cx e^{(t_E+\eta)^2x^2/2}\P(E).
\]
Since $p_x$ is bounded below,
\[
 \log\frac1{\P(E)}
 \le\frac{(t_E+\eta)^2x^2}{2}+\log\frac{Cx}{\delta\eta^2}
 \le\frac{t_E^2x^2}{2}+C(\eta x^2+\Lambda_x).
\]
Put $\mathsf R=\P(\,\cdot\mid E)$. The box bounds and $q_E\ge d_0$ give
\[
 \frac{\E_{\mathsf R}A}{x\sqrt{\E_PQ}}
 =-t_E\sqrt{q_E}+O(\eta),\qquad
 \frac{\E_{\mathsf R}Q}{\E_PQ}=q_E+O(\eta).
\]
The errors are uniform in $r$ and $\E_PQ$.
Substitution into \eqref{main:tilted-entropy} and
\eqref{core:conditional-entropy-upper} gives
\[
 \begin{aligned}
 D(\mathsf R\Vert P)&\ge x^2(a+t_E\sqrt{q_E})^2j(q_E)-C\eta x^2,\\
 D(\mathsf R\Vert P)&\le\frac{x^2}{2}(t_E+a\sqrt{q_E})^2
                              +C(\eta x^2+\Lambda_x).
 \end{aligned}
\]
Replacing $z_x$ by $1$ adds an $O(1)$ term, bounded by $C\Lambda_x$.
Since $q_E\le1-\delta$, Lemma~\ref{core:gap}, with
$(b,d)=(t_E,q_E)$, yields
\[
 \frac g2\delta x^2
 \le g(1-\sqrt{q_E})x^2
 \le C(\eta x^2+\Lambda_x)=o(\delta x^2),
\]
again a contradiction. This proves \eqref{core:no-drop-form} with
$x^{-1/8}$ replaced by $\delta$, uniformly in $r$ and $\E_PQ$.

Since $p_x\to p_*>0$, the exceptional conditional mass is $o(\delta)$.
If $H_x^\varsigma(Q\ge u)\ge2\delta$, we can therefore choose a pair
$(r,\omega)$ in $\{Q\ge u\}$ with $r\in C_x^\varsigma$ satisfying
\eqref{core:good-properties} with tolerance $\eta$.
Applying \eqref{core:no-drop-form} with $x^{-1/8}$ replaced by $\delta$
gives, for $u>0$,
\begin{equation}\label{core:reused-tail-transfer}
 H_x^\varsigma(Q\ge u)\ge2\delta
 \quad\Longrightarrow\quad
 H_x^{-\varsigma}(Q\le h_xu)<\delta,
 \qquad h_x=\frac{1-\delta}{1+K_0\eta}.
\end{equation}
Since $\eta=o(\delta)$ and $\delta\to0$, for all sufficiently large $x$,
\begin{equation}\label{core:reused-endpoints}
 \begin{gathered}
 h_x\ge1-2\delta,\qquad h_x(1+2\delta)\ge1,\\
 h_x(1-2\delta)\ge1-4\delta,\qquad
 h_x(1+4\delta)\ge1+2\delta.
 \end{gathered}
\end{equation}
Applying \eqref{def:band-median} as in the derivation of
\eqref{main:common-energy}, with \eqref{core:reused-tail-transfer}
and \eqref{core:reused-endpoints}, gives
\[
 H_x^\varsigma([m_x(1-4\delta),m_x(1+4\delta)])\ge1-3\delta.
\]
If $\delta_x=x^{-\alpha}$, Lemma~\ref{arr:overlap} now gives
$|\log m_{9x/8}-\log m_x|\le Kx^{-\alpha}$.
\end{proof}

Recall the variables in \eqref{app:radial-independent-laws}
and \eqref{def:residual-noise}:
\[
 \mathsf E=T-\varsigma z_xr\sqrt Q,\qquad
 \xi_x^\varsigma=\frac{U-\varsigma z_xrR}{\sqrt Q}.
\]

\begin{proof}[Proof of Lemma~\ref{lem:gaussian-energy-improvement}]
Lemmas~\ref{lem:tail-energy-initial} and~\ref{lem:band-to-global} give
\begin{equation}\label{array:initial-v}
 v_x^\varsigma=O(x^{-1/8}),\qquad\varsigma=\pm1.
\end{equation}
We begin with $\varsigma=1$; the other sign follows by replacing
$(A,U)$ with $(-A,-U)$. For a fixed $r>0$, write
$P_r=P_{z_xr,r^2/2}$.
Lemma~\ref{arr:canonical-tilt} then gives
$\E_{P_r}A=\sigma^2z_xr+\E_{P_r}U$.
By \eqref{core:rho}, the marginal law of $r$ is $\rho_x^+$,
and the conditional law on the original sample space given $r$ is
$(Q/(\E_{P_r}Q))P_r$.

Under $\Gamma_x^+$, condition on $(r,U,R)$.
Completing the square in $\varphi(Z)e^{z_xr\sigma Z}$ gives
$Z=\sigma z_xr+Z_0$. Here $Z_0$ is standard normal and independent
of $(r,U,R)$. The size-biasing factor $Q$ depends only
on $R$, so this conditional Gaussian law is unchanged. Thus
\[
 \begin{aligned}
 \mathsf E
 &=\frac{\sigma Z+U-z_xr(\sigma^2+R)}{\sqrt Q}
 =\frac{U-z_xrR}{\sqrt Q}
   +\frac{\sigma(Z-\sigma z_xr)}{\sqrt Q}
 =\xi_x^++\frac{\sigma Z_0}{\sqrt Q}.
 \end{aligned}
\]
By Lemma~\ref{app:radial-identities}, $\mathsf E\sim\N(0,1)$.
Since $\E_{\Gamma_x^+}\mathsf E^2=1$ and $\sigma^2/Q\le1$,
both $\sigma Z_0/\sqrt Q$ and $\xi_x^+$ are square integrable.
Conditional on $(r,U,R)$, $Z_0$ has mean zero and second moment one.
The cross term vanishes, giving
\[
 \E_{\Gamma_x^+}(\xi_x^+)^2
 =1-\E_{\Gamma_x^+}\frac{\sigma^2}{Q}
 =\E_{\Gamma_x^+}\frac RQ=v_x^+.
\]
Reflection proves \eqref{array:exact-residual-noise} for the other
sign.

To prove \eqref{array:improved-rates}, use the errors in
\eqref{notation:radial-errors}, which here take the form
\begin{equation}\label{def:radial-errors}
 d_r=\frac{\E_{P_r}U-z_xrR}{\sqrt Q},\qquad
 e_r=\frac{\E_{P_r}U-z_xr\E_{P_r}R}{\sqrt{\E_{P_r}Q}}.
\end{equation}
Since $\kappa=0$, \eqref{arr:residual-tilted-variance} gives
$\operatorname{Var}_{P_r}U\le\E_{P_r}R$. Under the conditional law $(Q/(\E_{P_r}Q))P_r$, the factor $Q$ cancels
the denominator in $(U-\E_{P_r}U)^2/Q$, giving
\begin{align*}
 \E_{\Gamma_x^+}\frac{(U-\E_{P_r}U)^2}{Q}
 &=\int\frac{\operatorname{Var}_{P_r}U}{\E_{P_r}Q}\,\rho_x^+(dr)
 \le\int\frac{\E_{P_r}R}{\E_{P_r}Q}\,\rho_x^+(dr)
 =\E_{\Gamma_x^+}\frac RQ=v_x^+.
\end{align*}
Using $d_r=\xi_x^+-(U-\E_{P_r}U)/\sqrt Q$ and
\eqref{array:exact-residual-noise},
\[
 \E_{\Gamma_x^+}d_r^2
 \le2\E_{\Gamma_x^+}(\xi_x^+)^2
   +2\E_{\Gamma_x^+}\frac{(U-\E_{P_r}U)^2}{Q}
 \le4v_x^+.
\]
Jensen's inequality gives
\[
 e_r^2=\frac{\bigl(\E_{P_r}(U-z_xrR)\bigr)^2}{\E_{P_r}Q}
 \le\frac{\E_{P_r}(U-z_xrR)^2}{\E_{P_r}Q}.
\]
Integrating against $\rho_x^+$ and canceling the factor
$Q/(\E_{P_r}Q)$ again gives
\[
 \int e_r^2\,d\rho_x^+
 \le\E_{\Gamma_x^+}(\xi_x^+)^2=v_x^+.
\]
The identity
\[
 z_xr(\E_{P_r}Q-Q)=\sqrt Q\,d_r-\sqrt{\E_{P_r}Q}\,e_r
\]
and $Q,\E_{P_r}Q>0$ imply
\[
 r\left|\sqrt{\E_{P_r}Q}-\sqrt Q\right|
 \le\frac{|d_r|+|e_r|}{z_x}.
\]
For $\eta_x>0$, let
\[
 B_x^+=\{\max(|d_r|,|e_r|,|\mathsf E|)>\eta_x x\}.
\]
Markov's inequality gives
\[
 \Gamma_x^+(B_x^+)
 \le\frac{\E_{\Gamma_x^+}d_r^2+\E_{\Gamma_x^+}e_r^2}
          {\eta_x^2x^2}
       +\Gamma_x^+(|\mathsf E|>\eta_x x).
\]
The Gaussian tail bound shows that $\Gamma_x^+(B_x^+)$ is at most
\begin{equation}\label{array:refined-bad-state}
 K\frac{v_x^+}{\eta_x^2x^2}+2e^{-\eta_x^2x^2/2},
\end{equation}
and $|d_r|,|e_r|,|\mathsf E|\le\eta_x x$ on $(B_x^+)^c$.

On $I_x^+\cap(B_x^+)^c$,
\[
 \begin{aligned}
 r\sqrt Q&=(T-\mathsf E)/z_x=T/z_x+O(\eta_x x),\\
 r\sqrt{\E_{P_r}Q}&=r\sqrt Q+O(\eta_x x).
 \end{aligned}
\]
For $\eta_x\to0$ with $x^{-2}=O(\eta_x)$,
$z_x=1+O(\eta_x)$, so both
$r\sqrt{\E_{P_r}Q}/x$ and $r\sqrt Q/x$ lie in
$[1-O(\eta_x),5/4+O(\eta_x)]$. Consequently,
\begin{equation}\label{array:refined-good-state}
 \begin{gathered}
 \frac{Q}{\E_{P_r}Q}=1+O(\eta_x),\qquad
 a:=\frac{\E_{P_r}A}{x\sqrt{\E_{P_r}Q}},\quad p:=\frac{r\sqrt{\E_{P_r}Q}}{x},\\
 a,p\in[1-O(\eta_x),5/4+O(\eta_x)],\qquad a-p=O(\eta_x).
 \end{gathered}
\end{equation}
The bound for $a-p$ follows from $a=z_xp+e_r/x$.
All constants are independent of $r$ and $\E_{P_r}Q$.
Reflection gives the negative-sign bounds with
\[
 a=\frac{\varsigma\E_{\varsigma z_xr,r^2/2}A}
        {x\sqrt{\E_{\varsigma z_xr,r^2/2}Q}}.
\]
Let $B_x^-$ be the corresponding reflected exceptional set.
Choose one fixed $K_0$ so that the bounds on $a,p$ in
\eqref{core:admissible-radii} and the estimate \eqref{core:good-properties}
hold with $x^{-1/4}$ replaced by $\eta_x$ on
$I_x^\varsigma\cap(B_x^\varsigma)^c$, for both signs.

For \eqref{array:improved-rates}, fix $\beta\in(2/3,7/10)$ and put
\[
 \eta_x=x^{-7/10},\qquad\delta_x=x^{-\beta}.
\]
By \eqref{array:initial-v}, the bound in
\eqref{array:refined-bad-state} is
\[
 O(x^{-1/8+7/5-2})+2e^{-x^{3/5}/2}
 =O(x^{-29/40})=o(\delta_x).
\]
Here $\beta<7/10<29/40$. In particular,
\[
 \begin{aligned}
 \frac{\Gamma_x^\varsigma(B_x^\varsigma)}{\delta_x}
 &\le Cx^{\beta-29/40}
      +2x^\beta e^{-x^{3/5}/2}\longrightarrow0,\\
 \frac{\eta_x}{\delta_x}
 &=x^{\beta-7/10}\longrightarrow0,\qquad
 \frac{\log x}{\delta_xx^2}\longrightarrow0,\\
 \log\frac1{\delta_x\eta_x^2}
 &=(\beta+7/5)\log x.
 \end{aligned}
\]
Since $p_x$ is bounded below, the conditional probability of $B_x^\varsigma$
given $I_x^\varsigma$ is also $o(\delta_x)$.
Lemma~\ref{arr:parameter-reuse} therefore gives
\[
 \begin{gathered}
 H_x^\varsigma([m_x(1-4x^{-\beta}),m_x(1+4x^{-\beta})])
       \ge1-3x^{-\beta},\\
 |\log m_{9x/8}-\log m_x|\le Cx^{-\beta}.
 \end{gathered}
\]
By Lemma~\ref{lem:tail-energy-initial}, $m_y\to\sigma^2>0$.
Summing the logarithmic increments gives
\[
 |\log(m_x/\sigma^2)|
 \le C\sum_{j\ge0}((9/8)^jx)^{-\beta}
 =\frac{Cx^{-\beta}}{1-(8/9)^\beta}.
\]
Exponentiating yields
\begin{equation}\label{array:improved-median-rate}
 |m_x-\sigma^2|\le Cx^{-\beta}.
\end{equation}
Since $Q\ge\sigma^2$ and
\[
 m_x(1+4x^{-\beta})
 \le(\sigma^2+Cx^{-\beta})(1+4x^{-\beta})
 \le\sigma^2+C'x^{-\beta},
\]
the concentration estimate in Lemma~\ref{arr:parameter-reuse}
proves the bound for $H_x^\varsigma$ in \eqref{array:improved-rates}.
Since $R/Q\le1$,
\[
 \E_{\Gamma_x^\varsigma}\!\left[\frac RQ\,\middle|\,I_x^\varsigma\right]
 \le\frac{C'x^{-\beta}}{\sigma^2}+3x^{-\beta}
 \le C''x^{-\beta}.
\]
Lemma~\ref{lem:band-to-global}, with $\beta\in(0,2)$, now gives
$v_x^\varsigma\le Cx^{-\beta}$.
\end{proof}

\section{Proofs of auxiliary lemmas}
\label{supp:auxiliary-proofs}

\subsection{Proof of Lemma~\ref{arr:transform-ui}}
\label{supp:transform-ui-proof}

\label{app:uniform-offset}

In the i.i.d.\ case, \citet[Theorem~2.5, p.~1520]{GGM1997} obtained uniform subgaussianity
from stochastic boundedness of the self-normalized sums.

\begin{theorem}
\label{new:offset}
For each $n$, let $Y_{n,1},\ldots,Y_{n,n}$ be independent,
almost surely finite real random variables.
Set $A_n=\sum_{i=1}^nY_{n,i}$ and
$\widehat V_n^2=\sum_{i=1}^nY_{n,i}^2$.
If $(A_n,\widehat V_n^2)$ is tight, then for every $\lambda>0$ there exists a
deterministic $d=d(\lambda)>0$ such that
\begin{equation}\label{new:offset-bound}
 \sup_n\E\exp\left\{\frac{\lambda A_n^2}{d+\widehat V_n^2}\right\}<\infty.
\end{equation}
The offset $d$ is independent of $n$ and may depend on $\lambda$
and the tight family.
\end{theorem}

\begin{proof}
For \eqref{new:offset-bound}, tightness of $\widehat V_n^2$ gives a finite $q\ge0$ such that
$\inf_n\P(\widehat V_n^2\le q)\ge1/2$. Thus
\[
 c:=\inf_n\E e^{-\widehat V_n^2}\ge e^{-q}/2>0.
\]
By independence, $-\log v\ge1-v$ for $0<v\le1$, and
\eqref{hist:eq:laplace-truncation},
\begin{align*}
 -\log\E e^{-\widehat V_n^2}
 &=\sum_i-\log\E e^{-Y_{n,i}^2}
 \ge\sum_i\E(1-e^{-Y_{n,i}^2})
 \ge(1-e^{-1})\sum_i\E\min(Y_{n,i}^2,1).
\end{align*}
Since $-\log\E e^{-\widehat V_n^2}\le-\log c$, there is
a finite deterministic constant $C$ such that
\begin{equation}\label{new:C}
 \sup_n\sum_i\E\min(Y_{n,i}^2,1)\le C.
\end{equation}

Truncate the coordinates at magnitude one by setting
\[
 Z_{n,i}=Y_{n,i}\ind_{\{|Y_{n,i}|\le1\}},
\]
and write
\[
 \begin{gathered}
 H_n=\sum_iZ_{n,i},\qquad K_n=\sum_iZ_{n,i}^2.
 \end{gathered}
\]
The truncated coordinates are bounded. Rowwise independence and
\eqref{new:C} imply
\[
 \E K_n=\sum_i\E Z_{n,i}^2\le C,
 \qquad
 \operatorname{Var}(H_n)
 =\sum_i\operatorname{Var}(Z_{n,i})\le C.
\]
Put $N_n=\sum_i\ind_{\{|Y_{n,i}|>1\}}$. For every $r>0$,
independence of these indicators gives
\begin{equation}\label{new:count}
 \begin{aligned}
 \E e^{rN_n}
 &=\prod_i\left[1+(e^r-1)\P(|Y_{n,i}|>1)\right]
 \le\exp\left\{(e^r-1)\sum_i\P(|Y_{n,i}|>1)\right\}
 \le\exp\{(e^r-1)C\}.
 \end{aligned}
\end{equation}
In particular, $N_n$ is tight. Define $L_n=A_n-H_n$ and
$R_n=\widehat V_n^2-K_n$. The Cauchy--Schwarz inequality gives
\[
 |L_n|\le\sqrt{N_nR_n}\le\sqrt{N_n\widehat V_n^2}.
\]
 Hence
$L_n$ is tight, and so is $H_n=A_n-L_n$.

The variance bound makes $H_n-\E H_n$ tight.
Choose $M<\infty$ so that, for every $n$,
\[
 \P(|H_n|\le M/2)\ge\tfrac34,\qquad
 \P(|H_n-\E H_n|\le M/2)\ge\tfrac34.
\]
The intersection has positive probability, so $\sup_n|\E H_n|\le M$.

To bound $\E e^{tH_n-at^2K_n}$ for $t\in\R$, fix $a>0$ and put
\[
 C_a=\max\left(0,\sup_{u\ne0}
       \frac{e^{u-au^2}-1-u}{u^2}\right).
\]
The quotient extends continuously at zero, with value $1/2-a$,
and tends to zero at both ends of the real line. Hence $C_a<\infty$ and
\begin{equation}\label{new:scalar-exponential-bound}
 e^{u-au^2}\le1+u+C_au^2\qquad(u\in\R).
\end{equation}
For each coordinate and each real $t$, apply \eqref{new:scalar-exponential-bound}
with $u=tZ_{n,i}$. Using $\log v\le v-1$ then gives
\[
 \log\E e^{tZ_{n,i}-at^2Z_{n,i}^2}
 \le t\E Z_{n,i}+C_at^2\E Z_{n,i}^2.
\]
These transforms are positive and finite because
$tz-at^2z^2\le1/(4a)$. Taking their product and using independence gives, for every $t\in\R$,
\begin{equation}\label{new:quadratic}
 \begin{aligned}
 \log\E e^{tH_n-at^2K_n}
 &\le t\E H_n+C_at^2\sum_i\E Z_{n,i}^2
 \le M|t|+C_aCt^2,
 \end{aligned}
\end{equation}

Fix $\lambda>0$ and choose
\[
 \ell=4\lambda,\qquad a=(4\ell)^{-1},\qquad
 d=4\ell(C_aC+1),\qquad D_n=d+K_n.
\]
Gaussian integration, Tonelli's theorem and \eqref{new:quadratic} give
\begin{align*}
 \E[D_n^{-1/2}e^{\ell H_n^2/D_n}]
 &=\frac1{\sqrt{4\pi\ell}}\int_\R
       e^{-dt^2/(4\ell)}
       \E e^{tH_n-t^2K_n/(4\ell)}\,dt\\
 &\le\frac1{\sqrt{4\pi\ell}}\int_\R
       \exp\left\{\left(-\frac d{4\ell}+C_aC\right)t^2
                         +M|t|\right\}\,dt\\
 &=\frac1{\sqrt{4\pi\ell}}\int_\R e^{-t^2+M|t|}\,dt
 =:K<\infty.
\end{align*}
Apply the Cauchy--Schwarz inequality to
$D_n^{-1/4}e^{2\lambda H_n^2/D_n}$ and $D_n^{1/4}$, then use
concavity of the square root:
\begin{equation}\label{new:remove-weight}
 \begin{aligned}
 \E e^{2\lambda H_n^2/D_n}
 &\le
 \{\E[D_n^{-1/2}e^{4\lambda H_n^2/D_n}]\}^{1/2}
 \{\E\sqrt{D_n}\}^{1/2}
 \le K^{1/2}(d+C)^{1/4}.
 \end{aligned}
\end{equation}
Here $4\lambda=\ell$ and
$\E\sqrt{D_n}\le\sqrt{d+\E K_n}\le\sqrt{d+C}$.

Use $L_n^2\le N_nR_n$ to include the large-coordinate sum $L_n$
in the exponential bound. For $R_n>0$,
the weighted Cauchy--Schwarz inequality gives
\[
 \frac{A_n^2}{d+\widehat V_n^2}
 =\frac{(H_n+L_n)^2}{(d+K_n)+R_n}
 \le\frac{H_n^2}{d+K_n}+\frac{L_n^2}{R_n}
 \le\frac{H_n^2}{d+K_n}+N_n.
\]
If $R_n=0$, then $L_n=N_n=0$, and the same inequality holds
with $L_n^2/R_n$ set to zero. Apply the Cauchy--Schwarz inequality
to the two exponentials. By \eqref{new:count} and \eqref{new:remove-weight},
\begin{align*}
 \E e^{\lambda A_n^2/(d+\widehat V_n^2)}
 &\le
 \{\E e^{2\lambda H_n^2/(d+K_n)}\}^{1/2}
 \{\E e^{2\lambda N_n}\}^{1/2}
 \le K^{1/4}(d+C)^{1/8}
       \exp\{(e^{2\lambda}-1)C/2\}<\infty.
\end{align*}
The bound is independent of $n$, proving
\eqref{new:offset-bound}.
\end{proof}

For a prescribed offset, the bound can fail. Take
$Y_{n,i}=\varepsilon_{n,i}/\sqrt n$, $1\le i\le n$, where the signs
in each row are independent and take the values $-1$ and $1$ with
equal probabilities.  Then $\widehat V_n^2=1$, $\E A_n=0$, and
$\E A_n^2=1$, so $(A_n,\widehat V_n^2)$ is tight by Chebyshev's inequality.
For any prescribed $d>0$, take $\lambda=d+1$.
The event that all signs are positive gives
\[
 \E\exp\left\{\frac{\lambda A_n^2}{d+\widehat V_n^2}\right\}
 =\E e^{A_n^2}
 \ge 2^{-n}e^n=(e/2)^n\longrightarrow\infty.
\]
Thus the offset must be allowed to depend on the exponent.

By Theorem~\ref{new:offset},
\begin{equation}\label{arr:offset-input}
 \forall\lambda>0\quad\exists d>0,\ C<\infty:
 \qquad
 \sup_n\E\exp\left\{\frac{\lambda A_n^2}{d+\widehat V_n^2}\right\}\le C.
\end{equation}
\begin{proof}[Proof of Lemma~\ref{arr:transform-ui}]
Recall that $(A_n,\widehat V_n^2)\dto(A,Q)$ and
$0<Q<\infty$ almost surely. In the tilted density below,
$\P$ denotes the original law of $(A,Q)$.
For fixed $s\in\R$ and $t>0$, the transforms and tilted law are
\[
 \begin{gathered}
 M_n(s,t)=\E e^{sA_n-t\widehat V_n^2},\qquad
 M(s,t)=\E e^{sA-tQ},\\
 dP_{s,t}=M(s,t)^{-1}e^{sA-tQ}\,d\P.
 \end{gathered}
\]
For $u\in\R$, $v\ge0$, and $d,\lambda>0$, completing the square gives
\begin{equation}\label{arr:tilt-square-bound}
 su-tv
 \le \frac{\lambda u^2}{d+v}
      +\frac{s^2d}{4\lambda}
      +\left(\frac{s^2}{4\lambda}-t\right)v.
\end{equation}
Choose $\lambda>s^2/(2t)$ and $d=d(\lambda)$ from
\eqref{arr:offset-input}. Applying
\eqref{arr:tilt-square-bound} with $(2s,2t)$ yields
\begin{equation}\label{arr:tilt-l2-bound}
 \sup_n\E e^{2sA_n-2t\widehat V_n^2}
 \le e^{s^2d/\lambda}
      \sup_n\E\exp\left\{\frac{\lambda A_n^2}{d+\widehat V_n^2}\right\}
 <\infty.
\end{equation}

For integers $m,k\ge0$, the bounds $|u|^m\le C_m e^{|u|}$ and
$v^k\le C_{k,t}e^{tv/2}$ give
\begin{equation}\label{arr:polynomial-absorption}
 |u|^m v^k e^{su-tv}
 \le C_{m,k,t}\left\{
    e^{(s+1)u-tv/2}+e^{(s-1)u-tv/2}\right\}.
\end{equation}
Apply \eqref{arr:tilt-l2-bound} at $(s+1,t/2)$ and $(s-1,t/2)$.
Equation~\eqref{arr:polynomial-absorption} then gives uniform $L^2$ bounds
for every polynomial times $e^{sA_n-t\widehat V_n^2}$.
Bounded continuous multipliers preserve uniform integrability.
Joint weak convergence gives the weighted expectation limits in
Lemma~\ref{arr:transform-ui}.
Their limits are finite, and $M(s,t)>0$ because $e^{sA-tQ}>0$.

To prove \eqref{arr:variance-interface}, factor the tilted density of
each finite row:
\[
 \frac{e^{sA_n-t\widehat V_n^2}}{M_n(s,t)}
 =\prod_{i\le n}
   \frac{e^{sY_{n,i}-tY_{n,i}^2}}
        {\E e^{sY_{n,i}-tY_{n,i}^2}}.
\]
The tilted coordinates remain independent and have finite second
moments, since $y^2e^{sy-ty^2}$ is bounded. Thus
\[
 \operatorname{Var}_{n;s,t}A_n
 =\sum_{i\le n}\operatorname{Var}_{n;s,t}Y_{n,i}
 \le\sum_{i\le n}\E_{n;s,t}Y_{n,i}^2
 =\E_{n;s,t}\widehat V_n^2.
\]
The weighted moment limits and
$M_n(s,t)\to M(s,t)>0$ give \eqref{arr:variance-interface}.
The tilted density is strictly positive, so $Q>0$ also gives $\E_{s,t}Q>0$.
\end{proof}

\subsection{Proof of Lemma~\ref{main:tilted-estimates}}
\label{supp:tilted-estimates-proof}

For $\mathsf R\ll P$, recall the relative entropy and the function
in \eqref{arr:j-definition}:
\[
 \begin{gathered}
 D(\mathsf R\Vert P)=\int\log\!\left(\frac{d\mathsf R}{dP}\right)
                         \,d\mathsf R,\\
 j(d)=\int_0^1\frac{1-u}{1-u+ud}\,du\quad(d>0),\qquad j(0)=1.
 \end{gathered}
\]
Logarithms are natural.

\begin{lemma}
\label{arr:finite-sharp-entropy}
Let $P$ be a product law of finitely many real coordinates
$(X_i)_{i\le m}$, each square integrable under $P$.
Let $\mathsf R\ll P$ be another probability law under which the coordinates are
also square integrable. Write the sums of the first and second moments as
\[
 \begin{gathered}
 m^P=\sum_i\E_PX_i,\qquad m^{\mathsf R}=\sum_i\E_{\mathsf R}X_i,\\
 q^P=\sum_i\E_PX_i^2>0,\qquad q^{\mathsf R}=\sum_i\E_{\mathsf R}X_i^2.
 \end{gathered}
\]
Then
\begin{equation}\label{arr:finite-sharp-bound}
 D(\mathsf R\Vert P)\ge\frac{(m^P-m^{\mathsf R})^2}{q^P}\,
 j\!\left(\frac{q^{\mathsf R}}{q^P}\right).
\end{equation}
\end{lemma}

\begin{proof}
It suffices to consider $D(\mathsf R\Vert P)<\infty$.
Write $P_i,\mathsf R_i$ for the marginal laws and
$L_i=d\mathsf R_i/dP_i$. For each coordinate, set
\[
 \begin{gathered}
 m_i^P=\E_PX_i,\qquad m_i^{\mathsf R}=\E_{\mathsf R}X_i,\\
 q_i^P=\E_PX_i^2,\qquad q_i^{\mathsf R}=\E_{\mathsf R}X_i^2.
 \end{gathered}
\]
The quantities in the lemma are the corresponding sums over $i$.
We omit coordinates with $q_i^P=0$, since they vanish under both $P$ and $\mathsf R$.

Apply the entropy representation of
\citet[Corollary~2, equation~(28)]{NS2020} to $(\mathsf R_i,P_i)$ and
reverse the integration parameter:
\begin{equation}\label{arr:coordinate-entropy-mixture}
 D(\mathsf R_i\Vert P_i)
 =\int_0^1(1-u)\,
     \E_{P_i}\frac{(L_i-1)^2}{1-u+uL_i}\,du.
\end{equation}
For $0<u<1$, the Cauchy--Schwarz inequality applied to
$\E_{P_i}[X_i(L_i-1)]=m_i^{\mathsf R}-m_i^P$ yields
\begin{equation}\label{arr:coordinate-entropy-cs}
 \begin{aligned}
 (m_i^P-m_i^{\mathsf R})^2
 &\le
 \E_{P_i}\frac{(L_i-1)^2}{1-u+uL_i}\,
 \E_{P_i}\!\left[X_i^2(1-u+uL_i)\right]\\
 &=
 \E_{P_i}\frac{(L_i-1)^2}{1-u+uL_i}\,
 \bigl((1-u)q_i^P+uq_i^{\mathsf R}\bigr).
 \end{aligned}
\end{equation}
To sum the coordinate bounds, let $L=d\mathsf R/dP$ and
$\mathsf Q=\bigotimes_i\mathsf R_i$.
Since $L_i(X_i)=\E_P[L\mid X_i]$, the conditional form of Jensen's inequality gives
\[
 D(\mathsf R_i\Vert P_i)\le D(\mathsf R\Vert P)<\infty.
\]
The negative part of $x\log x$ is bounded. Hence $\log L$ and all
$\log L_i(X_i)$ are integrable under $\mathsf R$.
Also, $\mathsf R(L_i(X_i)=0)=0$ for every $i$. Therefore
$\mathsf R\ll\mathsf Q$, and
\[
 \frac{d\mathsf Q}{dP}=\prod_iL_i(X_i),\qquad
 \frac{d\mathsf R}{d\mathsf Q}
 =\frac{L}{\prod_iL_i(X_i)}\quad\mathsf R\text{-almost surely}.
\]
Taking logarithms and expectations gives
\begin{equation}\label{arr:entropy-tensorization}
 D(\mathsf R\Vert P)
 =D\!\left(\mathsf R\,\middle\Vert\,\bigotimes_i\mathsf R_i\right)
       +\sum_iD(\mathsf R_i\Vert P_i)
 \ge\sum_iD(\mathsf R_i\Vert P_i).
\end{equation}

Combine \eqref{arr:entropy-tensorization},
\eqref{arr:coordinate-entropy-mixture}, and \eqref{arr:coordinate-entropy-cs}.
The Cauchy--Schwarz inequality across coordinates then gives
\begin{align*}
 D(\mathsf R\Vert P)
 &\ge\int_0^1(1-u)
       \sum_i\frac{(m_i^P-m_i^{\mathsf R})^2}{(1-u)q_i^P+uq_i^{\mathsf R}}\,du\\
 &\ge(m^P-m^{\mathsf R})^2\int_0^1
       \frac{1-u}{(1-u)q^P+uq^{\mathsf R}}\,du\\
 &=\frac{(m^P-m^{\mathsf R})^2}{q^P}\,
       j\!\left(\frac{q^{\mathsf R}}{q^P}\right).
\end{align*}
\end{proof}

\begin{proof}[Proof of Lemma~\ref{main:tilted-estimates}]
Fix $s\in\R$ and $t>0$, and put $P=P_{s,t}$.

First suppose
\[
 \frac{d\mathsf R}{dP}=\frac{g(A,Q)}{\E_Pg},
\]
where $g$ is nonnegative, bounded and continuous, and $\E_Pg>0$.
Let $P_n$ be the law of the full $n$th row tilted by
$e^{sA_n-t\widehat V_n^2}$. The weight factors across coordinates,
so $P_n$ is a product law. Lemma~\ref{arr:transform-ui} gives, with
$z_n=\E_{P_n}g(A_n,\widehat V_n^2)$,
\[
 z_n\longrightarrow\E_Pg>0,\qquad
 \E_{P_n}\widehat V_n^2\longrightarrow\E_PQ>0.
\]
For all sufficiently large $n$, both denominators are positive. Define
a law $\mathsf R_n$ on the full row by
\[
 \frac{d\mathsf R_n}{dP_n}
   =\frac{g(A_n,\widehat V_n^2)}{z_n}.
\]
The coordinates have finite second moments under both $P_n$ and
$\mathsf R_n$. Apply Lemma~\ref{arr:finite-sharp-entropy} on the full row:
\[
 D(\mathsf R_n\Vert P_n)
 \ge
 \frac{(\E_{P_n}A_n-\E_{\mathsf R_n}A_n)^2}
      {\E_{P_n}\widehat V_n^2}
 j\left(\frac{\E_{\mathsf R_n}\widehat V_n^2}{\E_{P_n}\widehat V_n^2}\right).
\]

The density depends only on $(A_n,\widehat V_n^2)$, so its entropy
can be computed from the pair:
\[
 D(\mathsf R_n\Vert P_n)
 =\frac{\E_{P_n}[g(A_n,\widehat V_n^2)\log g(A_n,\widehat V_n^2)]}{z_n}
       -\log z_n.
\]
The function $g\log g$ is bounded and continuous, with $0\log0=0$.
Lemma~\ref{arr:transform-ui} gives
$D(\mathsf R_n\Vert P_n)\to D(\mathsf R\Vert P)$.
The multipliers $1,A_n,\widehat V_n^2$ also give convergence of the means
of $A_n$ and $\widehat V_n^2$ under $P_n$ and $\mathsf R_n$.
Since $j$ is continuous, taking limits proves
\eqref{main:tilted-entropy} for bounded continuous densities.

Now let $f=d\mathsf R/dP$ be an arbitrary bounded density, say $0\le f\le K$.
The measure
\[
 d\mu=(1+|A|+Q)\,dP
\]
is a finite Borel measure on $\R^2$.
Approximate $f$ by bounded simple functions.
Regularity of $\mu$ gives continuous approximations to their indicators.
Truncating to $[0,K]$ yields continuous $g_m$ with
$0\le g_m\le K$ and $g_m\to f$ in $L^1(\mu)$.
Pass to a subsequence with $g_m\to f$ $P$-almost surely.

Let $z_m=\E_Pg_m$.  Then $z_m\to\E_Pf=1$, so $z_m>0$ eventually.
The weighted $L^1$ convergence gives
\[
 \frac{\E_P[Ag_m]}{z_m}\to\E_{\mathsf R}A,\qquad
 \frac{\E_P[Qg_m]}{z_m}\to\E_{\mathsf R}Q.
\]
The function $x\log x$ is bounded and continuous on $[0,K]$.
The dominated convergence theorem therefore also gives
\[
 \frac{\E_P[g_m\log g_m]}{z_m}-\log z_m
 \longrightarrow \E_P[f\log f]=D(\mathsf R\Vert P).
\]
Apply the continuous-density case to $g_m/z_m$ and let $m\to\infty$.
This proves \eqref{main:tilted-entropy}.
Boundedness of $f$ gives finite means, and $Q>0$ $P$-almost surely
implies $\E_{\mathsf R}Q>0$.
\end{proof}

\Needspace{20\baselineskip}
\subsection{Proof of Lemma~\ref{core:gap}}
\label{supp:gap-proof}

\begin{proof}
Let $a,b\in[a_*,a^*]$, where
$0<a_*<a^*<\infty$ and $a^*/a_*<\sqrt2$.
Put $y=\sqrt d\in(0,1)$. By \eqref{arr:j-definition},
\[
 j(y^2)-\frac1{1+y}
 =\frac{y(1-y^2+2y\log y)}{(1-y^2)^2}\ge0.
\]
Indeed, $f(y)=1-y^2+2y\log y$ satisfies
\[
 f'(y)=2(1-y+\log y)\le0,
 \qquad f(1)=0,
\]
so $f(y)\ge0$ for $0<y<1$.

Set $v=b/a$ and $v_*=a^*/a_*<\sqrt2$.  A direct expansion yields
\begin{align*}
 &2(1+vy)^2-(1+y)(v+y)^2=(1-y)\{2-v^2+(2+2v-2v^2)y+y^2\}.
\end{align*}
Since $0<v\le v_*<\sqrt2$,
$2+2v-2v^2>0$ and $2-v^2\ge2-v_*^2>0$. Thus
\begin{align*}
 (a+by)^2j(y^2)-\frac12(b+ay)^2
 &\ge
 a^2\left\{\frac{(1+vy)^2}{1+y}-\frac12(v+y)^2\right\}
 \ge\frac{a_*^2(2-v_*^2)}4(1-y).
\end{align*}
This proves \eqref{core:linear-gap} with $g=a_*^2(2-v_*^2)/4$.
\end{proof}

\end{document}